\documentclass[10pt]{article}
\usepackage[dvips, hmargin=1.6in, vmargin={1.4in, 1.9in}]{geometry}
\usepackage{amssymb, amsmath, amsthm, enumerate, microtype, stmaryrd, sectsty}
\usepackage[latin1]{inputenc}
\usepackage[dvips, arrow, matrix, tips, curve]{xy}
\SelectTips{cm}{10}

\allsectionsfont{\normalsize}

 \DeclareMathOperator{\colim}{colim}

\newcommand{\cat}[1]{\mathbf{#1}}
\newcommand{\op}{\mathrm{op}}
\newcommand{\id}{\mathrm{id}}
\newcommand{\thg}{{\mathord{\text{--}}}}

\newcommand{\abs}[1]{{\left|{#1}\right|}}

\newcommand{\res}[2]{\left.{#1}\right|_{#2}}

\newcommand{\spn}[1]{{\left<{#1}\right>}}

\newcommand{\defn}[1]{\textbf{#1}}
\newcommand{\cd}[2][]{\vcenter{\hbox{\xymatrix#1{#2}}}}
\newcommand{\cdl}[2][]{\xymatrix@1#1{#2}}

\renewcommand{\b}[1]{\boldsymbol{#1}}
\newcommand{\A}{{\mathcal A}}
\newcommand{\B}{{\mathcal B}}
\newcommand{\C}{{\mathcal C}}
\newcommand{\D}{{\mathcal D}}
\newcommand{\E}{{\mathcal E}}

\newcommand{\G}{{\mathcal G}}
\renewcommand{\H}{{\mathcal H}}
\newcommand{\I}{{\mathcal I}}

\newcommand{\K}{{\mathcal K}}
\renewcommand{\L}{{\mathcal L}}
\newcommand{\ELL}{{\mathcal L}}
\newcommand{\M}{{\mathcal M}}

\newcommand{\R}{{\mathcal R}}
\renewcommand{\S}{{\mathcal S}}

\newcommand{\V}{{\mathcal V}}
\newcommand{\W}{{\mathcal W}}

\newcommand{\xtor}[1]{\cdl[@1]{{} \ar[r]|-{\object@{|}}^{#1} & {}}}

\newcommand{\twocong}[2][0.5]{\ar@{}[#2] \save ?(#1)*{\cong}\restore}
\newcommand{\rtwocell}[3][0.5]{\ar@{}[#2] \ar@{=>}?(#1)+/l 0.2cm/;?(#1)+/r 0.2cm/^{#3}}
\newcommand{\ltwocell}[3][0.5]{\ar@{}[#2] \ar@{=>}?(#1)+/r 0.2cm/;?(#1)+/l 0.2cm/^{#3}}
\newcommand{\ltwocello}[3][0.5]{\ar@{}[#2] \ar@{=>}?(#1)+/r 0.2cm/;?(#1)+/l 0.2cm/_{#3}}
\newcommand{\dtwocell}[3][0.5]{\ar@{}[#2] \ar@{=>}?(#1)+/u  0.2cm/;?(#1)+/d 0.2cm/^{#3}}
\newcommand{\dthreecell}[3][0.5]{\ar@{}[#2] \ar@3{->}?(#1)+/u  0.2cm/;?(#1)+/d 0.2cm/^{#3}}
\newcommand{\utwocell}[3][0.5]{\ar@{}[#2] \ar@{=>}?(#1)+/d 0.2cm/;?(#1)+/u 0.2cm/_{#3}}
\newcommand{\dtwocelltarg}[3][0.5]{\ar@{}#2 \ar@{=>}?(#1)+/u  0.2cm/;?(#1)+/d 0.2cm/^{#3}}
\newcommand{\utwocelltarg}[3][0.5]{\ar@{}#2 \ar@{=>}?(#1)+/d  0.2cm/;?(#1)+/u 0.2cm/_{#3}}

\newcommand{\pullbackcorner}[1][dr]{\save*!/#1+1.2pc/#1:(1,-1)@^{|-}\restore}
\newcommand{\pushoutcorner}[1][dr]{\save*!/#1-1.2pc/#1:(-1,1)@^{|-}\restore}

\newdir{ >}{{}*!/-5pt/\dir{>}}

\newtheorem{Prop}{Proposition}
\newtheorem{Lem}[Prop]{Lemma}
\newtheorem{Cor}[Prop]{Corollary}

\theoremstyle{definition}
\newtheorem{Defn}[Prop]{Definition}
\newtheorem{Rk}[Prop]{Remark}

\newtheorem{Exs}[Prop]{Examples}
\newtheorem{Ex}[Prop]{Example}

\makeatletter \@namedef{itemize*}{\itemize\parsep\z@ \parskip\z@}
\@namedef{enditemize*}{\enditemize} \@namedef{enumerate*}{\enumerate\parsep\z@
\parskip\z@} \@namedef{endenumerate*}{\endenumerate} \makeatother

\begin{document}

\title{Cofibrantly generated natural weak factorisation systems}
\author{Richard Garner\thanks{Supported by a Research Fellowship of St John's College, Cambridge and a European Union Marie Curie Fellowship}\\Department of Mathematics, Uppsala University,\\Box 480, S-751 06 Uppsala, Sweden}
\maketitle
\begin{abstract}
\noindent There is an ``algebraisation'' of the notion of weak
factorisation system (w.f.s.) known as a \emph{natural weak
factorisation system}. In it, the two classes of maps of a w.f.s.
are replaced by two categories of maps-with-structure, where the
extra structure on a map now encodes a \emph{choice} of liftings
with respect to the other class. This extra structure has pleasant
consequences: for example, a natural w.f.s. on $\C$ induces a
canonical natural w.f.s. structure on any functor category $[\A,
\C]$.

In this paper, we define cofibrantly generated natural weak factorisation systems by analogy with
cofibrantly generated w.f.s.'s. We then construct them by a method which is reminiscent of
Quillen's small object argument but produces factorisations which are much smaller and easier to
handle, and show that the resultant natural w.f.s. is, in a suitable sense, \emph{freely} generated
by its generating cofibrations. Finally, we show that the two categories of maps-with-structure for
a natural w.f.s. are closed under all the constructions we would expect of them: (co)limits,
pushouts / pullbacks, transfinite composition, and so on.
\end{abstract}
\newcommand{\inj}{\textsf{in}}
\newcommand{\catc}{\K}
\newcommand{\comp}{\cat{Mon}(\K)}
\newcommand{\orth}{\mathop{\boxempty}}
\newcommand{\morth}{\mathop{\boxtimes}}
\newcommand{\lp}{\cat{Sq}}
\newcommand{\zig}{\cat{ZZ}}
\newcommand{\sorth}{\mathop{\bot}}
\newcommand{\lo}[1]{{}^{\orth} {#1}}
\newcommand{\ro}[1]{{#1}^{\orth}}
\newcommand{\mlo}[1]{{\vphantom{#1}}^{\morth} {#1}}
\newcommand{\mro}[1]{{#1}^{\morth}}
\newcommand{\Ar}[1]{{{#1}^\mathbf{2}}}
\newcommand{\wfs}{w.f.s.}
\newcommand{\nwfs}{n.w.f.s.}
\newcommand{\dom}{\mathrm{dom}}
\newcommand{\cod}{\mathrm{cod}}
\newcommand{\map}{\mathbf}
\newcommand{\Ll}{\mathsf L}
\newcommand{\Rr}{\mathsf R}
\newcommand{\Lmap}{\mathsf L\text-\cat{Map}}
\newcommand{\Rmap}{\mathsf R\text-\cat{Map}}
\newcommand{\Ff}{\cat{Ff}}
\let\L\undefined

\section{Introduction}
A \emph{weak factorisation system} on a category is given by two
classes of maps $\ELL$ and $\R$ which are related by a
``lifting-extension'' property guaranteeing the existence of
fill-ins for certain commuting squares, along with a way of
factorising an arbitrary map as $f = pi$, where the maps $i$ and $p$
lie in the respective classes $\ELL$ and $\R$. In typical examples,
these two classes of maps have distinctive feels to them: an
$\ELL$-map is given by freely ``glueing'' structure onto the source
of the map to obtain the target, whilst an $\R$-map allows one to
lift structure from the target to its source.

The most common place where weak factorisation systems (henceforth
\wfs's) arise is in \emph{Quillen model structures}
\cite{Quillen:homotopical} on a category: here one has two weak
factorisation systems (trivial cofibration, fibration) and
(cofibration, trivial fibration) which interact in a pleasant way,
providing a powerful framework within which one can do a lot of
abstract homotopy theory: one obtains formal notions of homotopy
category, homotopy equivalence, homotopy limits and colimits,
simplicial resolution, and so on.

The definition of weak factorisation system is, as the name suggests, a weakening of the older
categorical notion of \emph{orthogonal factorisation system} \cite{FreydKelly:continuous}, in which
the ``existence'' in the lifting-extension property becomes \emph{unique} existence. This makes the
theory of orthogonal factorisation systems cleaner than that of their weak cousins: for example,
their factorisations can always be chosen in a functorial way; their ``$\R$-maps'' are closed under
limits and their ``$\ELL$-maps'' under colimits; and they can be lifted with no effort to functor
categories. However, one also has that the factorisations themselves must be (essentially) unique,
making them ill-suited to homotopy-theoretic ends.

However, it turns out that if one is slightly more subtle about the way in which one defines a
\wfs, one can have both non-uniqueness of factorisations and also many of the pleasant properties
of orthogonal factorisation systems. The \emph{natural weak factorisation systems} (henceforth
\nwfs's) of \cite{nwfs} are an ``algebraisation'' of the concept of w.f.s. We review their formal
definition in Section \ref{Sec:nwfs}, but the intuition can be quickly illustrated by an analogy
with the notion of \emph{Grothendieck fibration}.

Given a functor $F \colon \C \to \D$, it may or may not have the \emph{property} of being a
Grothendieck fibration, namely that every arrow of $\D$ should have a cartesian lifting to $\C$.
However, if one asks for $F \colon \C \to \D$ to be a \emph{cloven fibration}, that is, to be
equipped with an explicitly given choice of cartesian lifting for every arrow of $\D$, then this is
no longer a property of $F$ but \emph{extra structure} borne by $F$. And in fact, this extra
structure is algebraic: there is a monad on $\Ar{\cat{Cat}}$~--~the category of arrows in
$\cat{Cat}$~--~whose algebras are precisely the cloven fibrations.

Likewise, for a \wfs\ on $\C$, we speak of a map having the \emph{property} of being an $\ELL$-map
or a $\R$-map, whilst in a \nwfs\ on $\C$, we speak instead of equipping a map with the
\emph{structure} of an $\ELL$-map or a $\R$-map. And again, this extra structure is (co)algebraic:
there is a monad $\Rr$ on $\Ar \C$ whose algebras are precisely the $\R$-maps in this new sense;
and dually, there is a comonad $\Ll$ on $\Ar \C$ whose coalgebras are the $\ELL$-maps.

In the language of the first paragraph, an $\ELL$-map now becomes an arrow \emph{together with an
explicit description} of how one should glue structure onto the source to obtain the target; and a
$\R$-map becomes an arrow \emph{together with an explicit description} of how one should lift
structure from the target to the source. These explicit descriptions conspire to give one a
\emph{canonical} choice of fill-ins for the ``lifting-extension'' property, whilst the one
remaining ingredient in a \wfs, namely factorisation, is already encoded in the comonad-monad pair
on $\Ar \C$: the functor parts of $\Ll$ and $\Rr$ simply send an arrow to the first and second
halves of its factorisation.

Natural \wfs's have certain advantages over plain \wfs's: for instance, the category of $\ELL$-maps
for a \nwfs\ is closed under \emph{all} colimits, and the category of $\R$-maps under \emph{all}
limits; moreover, \nwfs\ structures on $\C$ induce \nwfs\ structures on each functor category $[\A,
\C]$ in a completely canonical way. However, with this greater power comes greater complexity, and
thus one needs to do a good deal of groundwork to obtain a useful computational tool.

For example, one knows that the $\ELL$-maps (or dually, the $\R$-maps) for a \wfs\ are closed under
constructions like pushout, retracts, fibre coproducts and transfinite composition: we would
obviously like the same to be true for \nwfs's, and this is what we show in Section
\ref{Sec:props}. Because the $\ELL$-maps and $\R$-maps now carry extra structure, giving a precise
meaning to ``the pushout of an $\ELL$-map'' is a little more subtle, and showing that it exists a
little more involved: but beyond this, we find that we are able to proceed essentially as before.

The main meat of this paper, however, is Sections \ref{Sec:cgendef} and \ref{Sec:cgen:construct},
where we define and construct \emph{cofibrantly generated} \nwfs's. The definition generalises the
notion of a cofibrantly generated \wfs, a notion which describes almost every \wfs\ found in
nature; whilst the construction is both an adaptation of Quillen's \emph{small object argument} and
an example of the sort of \emph{free monoid} construction studied by Kelly in \cite{Ke80}. In fact,
we see that a cofibrantly generated \nwfs\ is, in a suitable sense, \emph{freely} generated by its
generating cofibrations. There are ramifications for the study of plain \wfs's as well, since our
method gives a recipe for the construction of functorial factorisations which are much less
redundant than Quillen's original argument, and which in many cases can be easily calculated by
hand.

\paragraph{Acknowledgements} I thank Walter Tholen for helpful feedback on an earlier draft of this
paper.

\section{Natural weak factorisation systems}
\label{Sec:nwfs} Let us start by recalling the notion of a \defn{weak factorisation system} $(\ELL,
\R)$ on a category $\C$. This is given by two classes of maps $\ELL$ and $\R$ in $\C$ which are
closed under retracts in the arrow category of $\C$, and which satisfy both
\begin{description}
\item[(lifting)] Whenever we are given a commutative square
\begin{equation}
\label{fillinsquare}
    \cd{
        A
            \ar[r]^-f
            \ar[d]_i &
        C
            \ar[d]^p \\
        B
            \ar[r]_-g &
        D
    }
\end{equation}
in $\C$, where $i \in \ELL$ and $p \in \R$, we can find a fill-in $j \colon B \to C$ such that $pj
= g$ and $ji = f$; and
\item[(factorisation)] Every map $f \colon X \to Y$ in $\C$ can be factorised as $f = pi$, where $i
\in \ELL$ and $p \in \R$.
\end{description}
In general, given maps $i \colon A \to B$ and $p \colon C \to D$, we call a square like
\eqref{fillinsquare} an \emph{$(i, p)$-lifting problem}. If every such square has a fill-in (or
``solution'') then we say that $i$ has the \emph{left lifting property} (llp) with respect to $p$
and that $p$ has the \emph{right lifting property} (rlp) with respect to $i$. Thus we can restate
the lifting axiom as: every $\R$-map has the rlp with respect to every $\ELL$-map, and vice versa.
In fact, the $\R$-maps are \emph{precisely} the maps with the rlp with respect to every $\ELL$-map,
and vice versa, so that the classes $\ELL$ and $\R$ determine each other.

One frequently requires that the factorisations for a \wfs\ should
be functorial in the following sense. Let us write $\Ar \C$ for the
arrow category of $\C$; we have two functors $\dom, \cod \colon \Ar
\C \to \C$ and a natural transformation $\kappa \colon \dom
\Rightarrow \cod$ whose component $\kappa_f$ is the map $f$. By a
\defn{functorial factorisation} $({E}, \lambda, \rho)$, we now mean a
functor ${E} \colon \Ar \C \to \C$ together with natural transformations
\[\cd{\dom \ar@2[r]^{\lambda} & {E} \ar@2[r]^-{\rho} & \cod}\]
satisfying $\kappa = \rho \cdot \lambda$. A \defn{functorial weak factorisation system} is given by
a \wfs\ $(\ELL, \R)$ together with a functorial factorisation $({E}, \lambda, \rho)$ such that each
$\lambda_f$ is in $\ELL$ and each $\rho_f$ is in $\R$. This notion is stronger than a plain \wfs,
but technically more convenient.

By strengthening a functorial \wfs\ further still, one arrives at
the notion of a natural \wfs, which as explained in the
introduction, consists of a comonad $\Ll$ and a monad $\Rr$ on $\Ar
\C$, interacting in a certain way. The definition we give is
essentially that of \cite{nwfs}, with the only novelty being the
addition of a \emph{distributive law} of the comonad over the monad.
This is a natural transformation $\Delta \colon LR \Rightarrow RL$
satisfying axioms expressing a form of compatibility between the
monad and the comonad; more precisely, it encodes a way of lifting
the comonad $\Ll$ to the category of free algebras for the monad
$\Rr$, and vice versa. We won't spell out the details here, since
equation \eqref{nwfs2} below tells us everything we need to know
about the distributive law; but the reader may like to consult
\cite{Beck} for further details.

\begin{Defn}\label{natwfs}
A \defn{natural weak factorisation system} $(\mathsf L, \mathsf R, \Delta)$ on a category $\C$ is
given by:
\begin{itemize*}
\item A comonad $\mathsf L = (L, \Phi, \Sigma)$ on $\Ar \C$,
\item A monad $\mathsf R = (R, \Lambda, \Pi)$ on $\Ar \C$, and
\item A distributive law $\Delta \colon LR \Rightarrow RL$,
\end{itemize*}
satisfying the following equalities:
\begin{align*}
\dom \cdot L & = \dom\text, & \cod \cdot L & = \dom \cdot R\text, &
\cod \cdot R & = \cod\text;\\
\dom \cdot \Phi & = 1_{\dom}\text, & \cod \cdot \Phi & = \kappa \cdot R\text, & \dom \cdot \Lambda
& = \kappa \cdot L\text, & \cod
\cdot \Lambda & = 1_{\cod}\text;\\
\text{and }\dom \cdot \Sigma & = 1_{\dom}\text, & \cod \cdot \Sigma & = \dom \cdot \Delta\text, &
\dom \cdot \Pi & = \cod \cdot \Delta\text, & \cod \cdot \Pi & = 1_{\cod}\text.
\end{align*}
\end{Defn}
We will shortly see that there is a good deal of redundancy in this definition: this is unavoidable
if we want to capture the (co)algebraic aspects of the system, but on the plus side means that we
can unravel the definition and give a much more compact description of a n.w.f.s. Firstly, we have
functors $L, R \colon \Ar \C \to \Ar \C$ satisfying $\mathord{\dom} \cdot L = \dom$,
$\mathord{\cod} \cdot L = \mathord{\dom} \cdot R$ and $\mathord{\cod} \cdot R = \cod$, which we
write as:
\[
L\left(\cd{X \ar[d]^f \\ Y}\right) = \cd{X \ar[d]^{\lambda_f} \\
{E}f\text,} \qquad
R\left(\cd{X \ar[d]^f \\ Y}\right) = \cd{{E}f \ar[d]^{\rho_f} \\
Y\text,}
\]\[
L\left(
    \cd[@1]{
        X \ar[d]_f \ar[r]^h &
        W \ar[d]^g \\
        Y \ar[r]_k &
        Z
    }
\right)  =
    \cd[@C+1em]{
        X \ar[d]_{\lambda_f} \ar[r]^h &
        W \ar[d]^{\lambda_g} \\
        {E}f \ar[r]_{{E}(h, k)} &
        {E}g\text,
    }
\qquad R\left(
    \cd[@1]{
        X \ar[d]_f \ar[r]^h &
        W \ar[d]^g \\
        Y \ar[r]_k &
        Z
    }
\right) =
    \cd[@C+1em]{
        {E}f \ar[d]_{\rho_f} \ar[r]^{{E}(h, k)} &
        {E}g \ar[d]^{\rho_g} \\
        Y \ar[r]_{k} &
        Z\text.
    }
\]
This diagram should illustrate our conventions for working with the category $\Ar \C$. When we view
a morphism $f \colon X \to Y$ of $\C$ as an object of $\Ar \C$ we will draw it \emph{vertically};
thus a morphism of $\Ar \C$ from $f$ to $g$ is a commutative square, bounded by $f$ and $g$
vertically and by two maps $h$ and $k$ horizontally, which we write as $(h, k) \colon f \to g$.

Next in the definition of \nwfs\ we have natural transformations $\Phi \colon L \Rightarrow \id$
and $\Lambda \colon \id \Rightarrow R$ satisfying $\mathord{\dom} \cdot \Phi = 1_{\dom}$,
$\mathord{\cod} \cdot \Phi = \kappa \cdot R$, $\mathord{\dom} \cdot \Lambda = \kappa \cdot L$ and
$\mathord{\cod} \cdot \Lambda = 1_{\cod}$. These conditions completely determine the components of
$\Phi$ and $\Lambda$ as being:
\[
\Phi_f =
    \cd[@1]{
        X \ar[d]_{\lambda_f} \ar[r]^{\id_X} &
        X \ar[d]^{f} \\
        {E}f \ar[r]_{\rho_f} &
        Y
    }
\quad \text{and} \quad \Lambda_f  =
    \cd[@1]{
        X \ar[d]_{f} \ar[r]^{\lambda_f} &
        {E}f \ar[d]^{\rho_f} \\
        Y \ar[r]_{\id_Y} &
        Y\text.
    }
\]
However, the existence of $\Phi$ and $\Lambda$ is not without force,
since the above squares must commute, which tells us that $f =
\rho_f \cdot \lambda_f$ for each $f \in \Ar \C$. Thus what we have
so far is precisely a functorial factorisation $({E}, \lambda,
\rho)$: we can think of a \nwfs\ as a
``functorial-factorisation-with-structure'', a viewpoint we will
espouse more comprehensively in the next section.

Continuing, we have the natural transformations $\Sigma \colon L \Rightarrow LL$, $\Delta \colon LR
\Rightarrow RL$ and $\Pi \colon RR \Rightarrow R$, satisfying $\mathord{\dom} \cdot \Sigma =
1_{\dom}$, $\cod \cdot \Sigma = \dom \cdot \Delta$, $\cod \cdot \Delta = \dom \cdot \Pi$, and
$\mathord{\cod} \cdot \Pi = 1_{\cod}$, and thus we have
\[
\Sigma_f =
    \cd[@1]{
        X \ar[d]_{\lambda_f} \ar[r]^{\id_X} &
        X \ar[d]^{\lambda_{Lf}} \\
        {E}f \ar[r]_-{\sigma_f} &
        {E}Lf\text,
    }\quad
\Delta_f  =
    \cd[@1]{
        {E}f \ar[r]^{\sigma_f} \ar[d]_{\lambda_{Rf}} &
        {E}Lf \ar[d]^{\rho_{Lf}} \\
        {E}Rf \ar[r]_{\pi_f} &
        {E}f
    }
\quad \text{and} \quad \Pi_f  =
    \cd[@1]{
        {E}Rf \ar[d]_{\rho_{Rf}} \ar[r]^{\pi_f} &
        {E}f \ar[d]^{\rho_f} \\
        Y \ar[r]_{\id_Y} &
        Y\text.
    }
\]
The intuition behind these maps is as follows. If we were dealing with a functorial \wfs, then the
$\lambda$'s would have the left lifting property with respect to the $\rho$'s, and so we would have
fill-ins for squares like this:
\[
    \cd[@1]{
        X \ar[d]_{\lambda_f} \ar[r]^-{\lambda_{Lf}} &
        {E}Lf \ar[d]^{\rho_{Lf}} \\
        {E}f \ar[r]_{\id} \ar@{.>}[ur] &
        {E}f
    }
\quad \text{and} \quad
    \cd[@1]{
        {E}f \ar[d]_{\lambda_{Rf}} \ar[r]^{\id} &
        {E}f \ar[d]^{\rho_f} \\
        {E}Rf \ar[r]_-{\rho_{Rf}} \ar@{.>}[ur] &
        Y\text.
    }
\]
This is what $\sigma_f$ and $\pi_f$ provide us with, in a uniform way. Indeed, we already know that
$\sigma_f$ and $\pi_f$ make the upper left and lower right triangles in the displayed squares
commute; that the same is true for the lower left and upper right follows from the comonad and
monad identities for $\Ll$ and $\Rr$. Explicitly, these identities assert that:
\begin{equation}
\label{nwfs1}
\begin{aligned}
  \rho_{Lf} \cdot \sigma_f & = \id_{{E}f}\text, & \pi_f \cdot \lambda_{Rf} & = \id_{{E}f}\text,\\
  {E}(1_X, \rho_f) \cdot \sigma_f & = \id_{{E}f}\text, & \pi_f \cdot {E}(\lambda_f, 1_Y)& = \id_{{E}f}\text,\\
  {E}(1_X, \sigma_f) \cdot \sigma_f & = \sigma_{Lf} \cdot \sigma_f &
  \text{and } \pi_f \cdot {E}(\pi_f, 1_Y)  & = \pi_f \cdot \pi_{Rf}
  \text.
\end{aligned}
\end{equation}
All that remains to account for are the axioms for the distributive law $\Delta$. Most of these
just repeat things we already know, and the only new equality we obtain is:
\begin{equation}
\label{nwfs2} \sigma_f \cdot \pi_f = \pi_{Lf} \cdot {E}(\sigma_f,
\pi_f) \cdot \sigma_{Rf}\text.
\end{equation}
The equations of \eqref{nwfs1} and \eqref{nwfs2} may seem rather puzzling at first; a reasonable
intuition is that they can be viewed as ensuring that every possible way of constructing a lifting
from the $\lambda$'s, $\rho$'s, $\sigma$'s and $\lambda$'s will give the same result. We can now
give the promised ``more compact'' version of the definition of \nwfs
\begin{Defn}\label{reduced}
A \defn{reduced \nwfs} $({E}, \lambda, \rho, \sigma, \pi)$ on $\C$
is given by:
\begin{itemize*}
\item A functorial factorisation $({E}, \lambda, \rho)$ on $\C$;
\item Natural transformations $\sigma \colon {E} \Rightarrow {E}L$ and $\pi \colon {E}R \Rightarrow
{E}$, where $L$ and $R$  are the unique functors $\Ar \C \to \Ar \C$
satisfying $\kappa \cdot L = \lambda$ and $\kappa \cdot R = \rho$,
\end{itemize*}
such that $\sigma \cdot \lambda = \lambda L$ and $\rho \cdot \pi = \rho R$, and such that the
equations of \eqref{nwfs1} and \eqref{nwfs2} hold.
\end{Defn}
\noindent From the preceding discussion, we see that \nwfs's on $\C$ are in bijection with reduced
natural \wfs's on $\C$ and thus we will pass between the two views without further comment.

Let us now examine the manner in which a \nwfs\ generalises a plain w.f.s. As explained in the
Introduction, we capture the ``$\ELL$-maps'' and ``$\R$-maps'' for a \nwfs\ by means of the
categories of (co)algebras for the (co)monad part of the n.w.f.s. So given given a \nwfs\ $(\Ll,
\Rr, \Delta)$, let us write $\Lmap$ for the category of coalgebras for $\Ll$, and call it the
\defn{category of $\Ll$-maps}; and similarly write $\Rmap$ for the category of algebras for $\Rr$ and call it the \defn{category of $\Rr$-maps}. Explicitly, $\Lmap$ has
\begin{itemize*}
\item \textbf{Objects} $(f, s)$ being arrows $f \colon X \to Y$
and $s \colon Y \to {E}f$ of $\C$ satisfying $s \cdot f =
\lambda_f$, $\rho_f \cdot  s = \id_Y$ and $\sigma_f \cdot s =
{E}(1_X, s) \cdot s$, and
\item \textbf{Morphisms} $(h, k) \colon (f, s) \to (g, t)$ being morphisms
$(h, k) \colon f \to g$ in $\Ar \C$ such that $t \cdot k = {E}(h,
k)\cdot s$,
\end{itemize*}
whilst $\Rmap$ has
\begin{itemize*}
\item \textbf{Objects} $(f, p)$ being arrows $f \colon X \to Y$ and $p \colon {E}f \to X$ of $\C$
satisfying $f\cdot p = \rho_f$, $p \cdot \lambda_f = \id_X$ and $p
\cdot \pi_f = p \cdot {E}(p, 1_Y)$, and
\item \textbf{Morphisms} $(h, k) \colon (f, p) \to (g, q)$ being morphisms
$(h, k) \colon f \to g$ in $\Ar \C$ such that $h \cdot p = q \cdot
{E}(h, k)$.
\end{itemize*}
We will sometimes abuse notation slightly, and write $(f, s) \colon
X \to Y$ to signify an $\Ll$-map or $\Rr$-map for which $f \colon X
\to Y$; this emphasises the idea that an $\Ll$-map or $\Rr$-map is
just a map of $\C$ with extra structure. Now, to see that these
definitions make sense, consider the case where we have a mere
functorial \wfs\ $(\ELL, \R)$ and are given an $\ELL$-map $f \colon
A \to B$. If we take the factorisation of $f$ as $f = \rho_f \cdot
\lambda_f$, then, since every $\rho$ is an $\R$-map, we will have a
solution $s$ to the lifting problem
\[
    \cd{
        A
            \ar[r]^-{\lambda_f}
            \ar[d]_{f} &
        {E}f
            \ar[d]^{\rho_f} \\
        B
            \ar@{.>}[ur]_s
            \ar@{=}[r] &
        B\text.
    }
\]
It is lifting data of this form which accompanies the $\Ll$-maps and
$\Rr$-maps for a \nwfs; moreover, this extra data is sufficient to
give us canonical solutions to \emph{all} ``$(\Ll, \Rr)$-lifting
problems''. More precisely, if we are given an $\Ll$-map $(f, s)$,
an $\Rr$-map $(g, p)$, and an $(f, g)$-lifting problem
\[
    \cd{
        A
            \ar[r]^-h
            \ar[d]_{f} &
        C
            \ar[d]^{g} \\
        B
            \ar[r]_-k &
        D\text,
    }
\]
then we have a canonical choice of lifting $j \colon B \to C$ given
by \[j := B \xrightarrow{s} Ff \xrightarrow{F(h, k)} Fg
\xrightarrow{p} C\text,\] which is now \emph{natural}, in that it is
stable under composition with morphisms of $\Lmap$ on the left and
morphisms of $\Rmap$ on the right.

The remaining ingredient in a \nwfs\ is of course
\emph{factorisation}. Because $\Lmap$ and $\Rmap$ are categories of
(co)algebras, we have adjunctions:
\[\cd{\Lmap \ar@<4pt>[r]^-{U_\Ll} \ar@{}[r]|-{\bot} & \Ar \C \ar@<4pt>[l]^-{F_\Ll}} \qquad \text{and} \qquad \cd{\Rmap \ar@<4pt>[r]^-{U_\Rr} \ar@{}[r]|-{\top} & \Ar \C\text. \ar@<4pt>[l]^-{F_\Rr}}\]
The forgetful functors $U_\Ll$ and $U_\Rr$ send an $\Ll$-map $(f, s)$ or $\Rr$-map $(f, p)$ to its
underlying $\C$-map $f$, whilst the free\footnote{We should call $F_\Ll$ ``cofree'', but we won't
labour the point.} functors $F_\Ll$ and $F_\Rr$ respectively send a map $f \colon A \to B$ of $\C$
to the $\Ll$-map $(\lambda_f, \sigma_f)$ and the $\Rr$-map $(\rho_f, \pi_f)$. Thus $F_\Ll$ and
$F_\Rr$ give us a functorial factorisation of any map of $\C$ into an $\Ll$-map followed by an
$\Rr$-map:
\[
f \colon A \to B = A \xrightarrow{(\lambda_f, \sigma_f)} {E}f
\xrightarrow{(\rho_f, \pi_f)} B\text.
\]

Now, underlying each \nwfs\ $(\Ll, \Rr, \Delta)$ is a functorial
\wfs: if we let $\ELL$ be the class of maps in $\C$ which admit some
$\Ll$-coalgebra structure, and $\R$ be the class of maps admitting
some $\Rr$-algebra structure, then the pair $(\ELL, \R)$ satisfy all
the conditions for a functorial \wfs\ except, possibly, closure of
$\ELL$ and $\R$ under retracts. So if we write $\overline \ELL$ and
$\overline \R$ for the respective retract-closures, we obtain a
functorial \wfs\ $(\overline \ELL, \overline \R)$, with the property
that the given factorisations land inside the smaller classes $\ELL$
and $\R$.

\begin{Exs}\hfill
\begin{itemize}
\item If $(\Ll, \Rr, \Delta)$ is a \nwfs\ on $\C$, then $(\Rr^\op, \Ll^\op, \Delta^\op)$ is a
\nwfs\ on $\C^\op$, and so the notion of \nwfs\ is self-dual.
\item If $(\Ll, \Rr, \Delta)$ is a \nwfs\ on $\C$, then for any object $X \in \C$ we induce a \nwfs\ of the
same name on the slice category $\C / X$ and the coslice category $X / \C$.
\item If $(\Ll, \Rr, \Delta)$ is a \nwfs\ on $\C$, then for any other category $\A$ we induce a \nwfs\
calculated pointwise on $[\A, \C]$. This stands in strong contrast to the situation with \wfs's,
where there is \emph{no} canonical lifting to functor categories.
\end{itemize}
\end{Exs}
We will not give any substantial examples now, as these will arise in due course from the theory of
\emph{cofibrantly generated} \nwfs's that we develop in Sections \ref{Sec:cgendef} and
\ref{Sec:cgen:construct}: thus the reader may like to look ahead to these examples, or to look at
those given in \cite{nwfs}.

\section{An alternative view of natural weak factorisation systems}\label{Sec:alternativeview}
We observed in the previous section that every \nwfs\ has an underlying functorial factorisation.
In this section we shall go in the other direction, and characterise \nwfs's as functorial
factorisations equipped with a \emph{bialgebra} structure.

Classically, a bialgebra is a vector space $A$ equipped with both an algebra and a coalgebra
structure, such that the coalgebra maps $\Delta \colon A \to A \otimes A$ and $\epsilon \colon A
\to k$ are algebra homomorphisms; the notion of bialgebra that we deploy to characterise \nwfs's is
a mild generalisation of this. This is an intuitively plausible idea, since both bialgebras and
\nwfs's have a ``multiplicative'' and a ``comultiplicative'' part satisfying compatibility
conditions. However, to go from \emph{plausible} to \emph{precise} will require a little work.

\subsection{Bialgebras in 2-fold monoidal categories}
One obvious way to generalise the notion of bialgebra is to restate its definition of in an
arbitrary symmetric (or braided) monoidal category; however, we will need something slightly more
general still, namely bialgebra objects where the ``algebra'' and the ``coalgebra'' parts are given
with respect to two \emph{different} monoidal structures on the same category. In order to express
the compatibility of the algebra and coalgebra parts, we first need a higher-level compatibility
between the two monoidal category structures with respect to which they are taken. This
compatibility is captured by the concept of a \emph{2-fold monoidal category}. In fact, we only
really need a \emph{strict} 2-fold monoidal category in this case, which makes the definition a
little simpler.

We recall first that a \emph{lax monoidal functor} between strict monoidal categories $\V$ and $\W$
is given by a functor $F \colon \V \to \W$ together with a natural family of
(not-necessarily-invertible) maps $m_{a, b} \colon Fa \otimes Fb \to F(a \otimes b)$ and a map $m_I
\colon I \to FI$, satisfying two coherence axioms: the first equates the two obvious ways of
getting from $Fa \otimes Fb \otimes Fc$ to $F(a \otimes b \otimes c)$, and the second says that
$m_{a, I} = m_{I, a} = \id_{Fa}$ for all $a \in \V$. One can compose lax monoidal functors to
obtain a category $\cat{StrMonCat}_{lax}$ of strict monoidal categories and lax monoidal functors;
and since $\cat{StrMonCat}_{lax}$ has finite products, we can consider monoids in it.

\begin{Defn}
A \defn{strict 2-fold monoidal category} is a monoid in $\cat{StrMonCat}_{lax}$.
\end{Defn}

If we expand this definition, a strict 2-fold monoidal category
consists of a category $\V$, two strict monoidal structures
$(\otimes, I)$ and $(\odot, \bot)$ on it, maps $m \colon \bot
\otimes \bot \to \bot$, $c \colon I \to I \odot I$ and $j \colon I
\to \bot$ making $(\bot, j, m)$ into a $\otimes$-monoid and $(I, j,
c)$ into a $\odot$-comonoid, and a natural family of maps
\[z_{A,B,C,D} \colon (A \odot B) \otimes (C \odot D) \to (A \otimes C) \odot (B \otimes D)\]
obeying six coherence laws, which equate, respectively, the two possible ways of getting from:
\begin{align*}
  (A \odot B \odot C) \otimes (A' \odot B' \odot C') & \quad \text{to} \quad (A \otimes A') \odot (B \otimes B') \odot (C \otimes C')\text,\\
  (A \odot A') \otimes (B \odot B') \otimes (C \odot C') & \quad \text{to} \quad (A \otimes B \otimes C) \odot (A' \otimes B' \otimes C')\text,\\
  (A \odot B) \otimes I & \quad \text{to} \quad (A \otimes I) \odot (B \otimes I)\text,\\
  I \otimes (A \odot B) & \quad \text{to} \quad (I \otimes A) \odot (I \otimes B)\text,\\
  (\bot \odot A) \otimes (\bot \odot B)& \quad \text{to} \quad \bot \odot (A \otimes B)\text,\\
  \text{and}\quad (A \odot \bot) \otimes (B \odot \bot)& \quad \text{to} \quad (A \otimes B) \odot \bot\text.
\end{align*}
We will write such a 2-fold monoidal category as $(\V, \otimes, I,
\odot, \bot)$. A 2-fold monoidal category is the simplest example of
an \emph{iterated} monoidal category, in the sense of \cite{iter1,
iter2}, to which we refer the reader for further examples and
applications. Note that, in one aspect, our definition is slightly
more general than those of the above-cited papers, since it does not
assume that the units $I$ and $\bot$ coincide.

Let us now see why a 2-fold monoidal category is a suitable environment for defining a notion of
bialgebra. What we want to say is the following: a bialgebra in $(\V, \otimes, I, \odot, \bot)$ is
an object $A$ together with maps $\mu \colon A \otimes A \to A$, $\eta \colon I \to A$, $\Delta
\colon A \to A \odot A$ and $\epsilon \colon A \to \bot$ such that $(A, \mu, \eta)$ is a monoid,
$(A, \Delta, \epsilon)$ is a comonoid, and such that the maps $\Delta \colon A \to A \odot A$ and
$\epsilon \colon A \to \bot$ are monoid homomorphisms.

For this last clause to make sense, we need $\otimes$-monoid
structures on $A \odot A$ and on $\bot$; and one way of obtaining
these is by lifting the $(\odot, \bot)$ monoidal structure on $\V$
to the category $\cat{Mon}_\otimes(\V)$ of $\otimes$-monoid objects
in $\V$. But this is precisely what the 2-fold monoidal structure
allows us to do. The unit for this lifted monoidal structure is the
$\otimes$-monoid $(\bot, j, m)$, whilst the tensor product of two
$\otimes$-monoids $(A, \eta^A, \mu^A)$ and $(B, \eta^B, \mu^B)$ is
given by $(A \odot B, \eta^{A \odot B}, \mu^{A \odot B})$, where
\[\eta^{A \odot B} = I \xrightarrow{c} I \odot I \xrightarrow {\eta^A \odot \eta^B} A \odot B\]
and
\[\mu^{A \odot B} = (A \odot B) \otimes (A \odot B) \xrightarrow{z_{A, B, A, B}} (A \otimes A) \odot (B \otimes B)
\xrightarrow{\mu^A \odot \mu^B} A \odot B\text.\]
This lifting
process can be seen more abstractly by noting that the operation
that assigns to each strict monoidal category the category of
monoids in it extends to a finite-product preserving functor
$\cat{Mon}(\thg) \colon \cat{StrMonCat}_{lax} \to \cat{Cat}$. Thus
monoids in $\cat{StrMonCat}_{lax}$~--~which are 2-fold monoidal
categories~--~are sent to monoids in $\cat{Cat}$~--~which are strict
monoidal categories. Regardless of how we obtain it, this lifting
allows us to define:
\begin{Defn}\label{bialgdef}
Let $(\V, \otimes, I, \odot, \bot)$ be a 2-fold monoidal category. The category $\cat{Bialg}(\V)$
of \defn{bialgebras in $\V$} is given by $\cat{Comon}_\odot(\cat{Mon}_\otimes(\V))$, the category
of $\odot$-comonoid objects in $\cat{Mon}_\otimes(\V)$.
\end{Defn}

Explicitly, such a bialgebra is given by a quintuple $(A, \eta, \mu,
\epsilon, \Delta)$ as above, such that such that $(A, \eta, \mu)$ is
a $\otimes$-monoid, $(A, \epsilon, \Delta)$ is a $\odot$-comonoid,
and such the following four diagrams commute:
\begin{equation}
\label{bialg}
\begin{gathered}
\cd{
 I \ar[r]^-\eta \ar[d]_{c} &
 A \ar[d]^\Delta \\
 I \odot I \ar[r]_-{\eta \odot \eta} &
 A \odot A\text,
} \qquad \cd{
 A \otimes A \ar[r]^-\mu \ar[d]_{\epsilon \otimes \epsilon} &
 A \ar[d]^\epsilon \\
 \bot \otimes \bot \ar[r]_-m & \bot\text,
} \qquad \cd{
 & A \ar[dr]^\epsilon \\
 I \ar[ur]^\eta \ar[rr]_j & & \bot\text,
}
\\
\cd[@C+1em]{
 A \otimes A
   \ar[rr]^\mu \ar[d]_{\Delta \otimes \Delta} & &
 A \ar[d]^\Delta \\
 (A \odot A) \otimes (A \odot A)
 \ar[r]_-{z_{A, A, A, A}} &  (A \otimes A) \odot (A \otimes A)
 \ar[r]_-{\mu \odot \mu} &
 A \odot A\text,
 }
\end{gathered}
\end{equation}
whilst a bialgebra homomorphism is a morphism of $\V$ which is simultaneously a monoid homomorphism
and a comonoid homomorphism.

\begin{Rk}\label{otherbialgview} Before returning to our pursuit of \nwfs's, we note that we can obtain the notion of
bialgebra in a 2-fold monoidal category in a dual way: it is not
only a ``$\odot$-comonoid in the category of $\otimes$-monoids'' but
also a ``$\otimes$-monoid in the category of $\odot$-comonoids''.
Indeed, we can lift the $(\otimes, I)$ monoidal structure on $\V$ to
the category $\cat{Comon}_\odot(\V)$ of $\odot$-comonoids: the unit
object is $(I, j, c)$ whilst the tensor product of two
$\odot$-comonoids $(A, \epsilon^A, \Delta^A)$ and $(B, \epsilon^B,
\Delta^B)$ is given by $(A \otimes B, \epsilon^{A \otimes B},
\Delta^{A \otimes B})$, where
\[\epsilon^{A \otimes B} = A \otimes B \xrightarrow{\epsilon^A \otimes \epsilon^B} \bot \otimes \bot \xrightarrow{m} \bot\] and
\[\Delta^{A \otimes B} = A \otimes B \xrightarrow{\Delta^A \otimes \Delta^B} (A \odot A) \otimes (B \odot B) \xrightarrow{z_{A, A, B, B}} (A \otimes B) \odot (A \otimes
B)\text.
\]And a $\otimes$-monoid in $\cat{Comon}_\odot(\V)$ is once again
a bialgebra in $\V$. To see this abstractly, observe that if $(\V,
\otimes, I, \odot, \bot)$ is a 2-fold monoidal category, then so is
$(\V^\op, \odot, \bot, \otimes, I)$; and from the explicit
definition of bialgebras given above, we can easily see that
$\cat{Bialg}(\V^\op) \cong \cat{Bialg}(\V)^\op$. So now
\begin{align*}
\cat{Bialg}(\V)^\op & \cong \cat{Bialg}(\V^\op) \\
& = \cat{Comon}_\otimes(\cat{Mon}_\odot(\V^\op)) \\
& \cong \cat{Comon}_\otimes(\cat{Comon}_\odot(\V)^\op) \\
& \cong \cat{Mon}_\otimes(\cat{Comon}_\odot(\V))^\op
\end{align*}
so that $\cat{Bialg}(\V) \cong \cat{Mon}_\otimes(\cat{Comon}_\odot(\V))$ as claimed. This second
characterisation of bialgebras will be the most useful to us when we are working with \nwfs's.
\end{Rk}

\subsection{Natural weak factorisation systems as bialgebras}\label{Sec:nwfsasbialg}
We are now ready to characterise \nwfs's on a category $\C$, which as we have already suggested,
will arise as functorial factorisations on $\C$ bearing a bialgebra structure. We now know that for
this to make sense, we need to organise functorial factorisations on $\C$ into a 2-fold monoidal
category.

Making them form a category~--~let us call it $\Ff_\C$~--~is easy
enough: objects are functorial factorisations $({E}, \lambda,
\rho)$, and morphisms $\alpha \colon ({E}, \lambda, \rho) \to ({E}',
\lambda', \rho')$ are natural transformations $\alpha \colon {E}
\Rightarrow {E}'$ making the diagram
\[
\cd{
 & \dom \ar@2[dl]_\lambda  \ar@2[dr]^{\lambda'} \\
 {E} \ar@2[rr]_\alpha \ar@2[dr]_{\rho} & &
 {E}' \ar@2[dl]^{\rho'} \\ &
 \cod
}
\]
commute. What remains is to give the two interacting monoidal structures $(\otimes, I)$ and
$(\odot, \bot)$ on $\Ff_\C$: and these arise very naturally from two ways of combining functorial
factorisations. In the first such, the tensor product $(E', \lambda', \rho') \otimes ({E}, \lambda,
\rho)$ factorises $f \colon X \to Y$ by applying ${E}$ to it, then applying ${E}'$ to the right
half of this factorisation:
\[X \xrightarrow{f} Y \quad \mapsto \quad X \xrightarrow{\lambda_f} {E}f \xrightarrow{\rho_f}
Y \quad \mapsto \quad X \xrightarrow{\lambda_f} {E}f \xrightarrow{\lambda'_{Rf}} {E}'Rf
\xrightarrow{\rho'_{Rf}} Y\text,\footnote{Note that, as in section \ref{Sec:nwfs}, we write $R$ for
the functor $\Ar \C \to \Ar \C$ corresponding to the natural transformation $\rho \colon {E}
\Rightarrow \cod$, and so on. We shall continue to do this without further note throughout the
paper.}\] and finally composing together the two ``left'' parts, $\lambda_f$ and $\lambda'_{Rf}$ to
obtain the factorisation
\[X \xrightarrow{f} Y \quad \mapsto \quad X \xrightarrow{\lambda'_{Rf} \cdot \lambda_f} {E}'Rf \xrightarrow{\rho'_{Rf}} Y\text.\]
The unit $I$ for this tensor product is also the initial object of $\Ff_\C$, namely $(\dom, 1_\dom,
\kappa)$ which factorises $f \colon X \to Y$ as \[X \xrightarrow{1_X} X \xrightarrow f Y\text.\]

The second monoidal structure is completely dual to the first; so
$({E}', \lambda', \rho') \odot ({E}, \lambda, \rho)$ factorises $f
\colon X \to Y$ by applying ${E}$ to it, then applying ${E}'$ to the
\emph{left} half of this factorisation:
\[X \xrightarrow{f} Y \quad \mapsto \quad X \xrightarrow{\lambda_f} {E}f \xrightarrow{\rho_f}
Y \quad \mapsto \quad X \xrightarrow{\lambda'_{Lf}} {E}'Lf
\xrightarrow{\rho'_{Lf}} {E}f \xrightarrow{\rho_f} Y\text,\] and
finally composing together the two ``right'' parts, $\rho_f$ and
$\rho'_{Lf}$ to obtain the factorisation
\[X \xrightarrow{f} Y \quad \mapsto \quad X \xrightarrow{\lambda'_{Lf}} {E}'Lf \xrightarrow{\rho_f \cdot \rho'_{Lf}} Y\text.\]
The unit $\bot$ for this tensor product is now the terminal object of $\Ff_\C$, namely $(\cod,
\kappa, 1_\cod)$ which factorises $f \colon X \to Y$ as \[X \xrightarrow{f} Y \xrightarrow {1_Y}
Y\text.\] One can easily check directly that these operations yield two strict monoidal structures
on $\Ff_\C$, but it will be more illuminating to see how we can deduce their existence indirectly.
Let us say that a functor $F \colon \Ar \C \to \Ar \C$ is \emph{over} $\cod \colon \Ar \C \to \C$
if $\cod \cdot F = \cod \colon \Ar \C \to \C$, and likewise, that a natural transformation $\alpha
\colon F \Rightarrow G \colon \Ar \C \to \Ar \C$ is over $\cod$ if $\cod \cdot \alpha = \id_\cod$.
Now it's easy to show that $\Ff_\C$ is isomorphic to the category $\Ff'_\C$ with:
\begin{itemize*}
\item \textbf{Objects} being pairs $(R, \Lambda)$ where $R \colon \Ar \C \to \Ar \C$ and $\Lambda
\colon \id_{\Ar \C} \Rightarrow R$ over $\cod$.
\item \textbf{Morphisms} $\Gamma \colon (R, \Lambda) \to (R', \Lambda')$ being commutative
triangles over $\cod$:
\[
   \cd{ & \id_{\Ar \C} \ar@2[dl]_\Lambda \ar@2[dr]^{\Lambda'} \\ R \ar@2[rr]_\Gamma & & R'\text.}
\]
\end{itemize*}
Now $\Ff'_\C$ has a strict monoidal structure on it:
\[
I = \cd{\id_{\Ar \C} \ar@2[d]^{\id} \\ \id_{\Ar \C}} \qquad \text{and} \qquad \cd{\id_{\Ar \C}
\ar@2[d]^{\Lambda} \\ R} \ \otimes \ \cd{\id_{\Ar \C} \ar@2[d]^{\Lambda'} \\ R'} \ = \ \cd{\id_{\Ar
\C} \ar@2[d]^{\Lambda\Lambda'} \\ RR'\text,}
\]
which transfers back to $\Ff_\C$ to give us the same-named structure
there: and so we deduce the associativity and unitality of the
latter from that of the former. Moreover, we can now easily classify
$\otimes$-monoids in $\Ff_\C$, since they correspond to
$\otimes$-monoids in $\Ff'_\C$; and giving a monoid structure on
$(R, \Lambda) \in \Ff'_\C$ is the same as giving a natural
transformation $\Pi \colon RR \Rightarrow R$ making $(R, \Lambda,
\Pi)$ into a monad over $\cod$.

We can argue dually for the $(\odot, \bot)$ monoidal structure, where the first step is now the
observation that $\Ff_\C$ is isomorphic to the category with:
\begin{itemize*}
\item \textbf{Objects} being pairs $(L, \Phi)$ where $L \colon \Ar \C \to \Ar \C$ and $\Phi
\colon L \Rightarrow \id_{\Ar \C}$ over $\dom$;
\item \textbf{Morphisms} $\Gamma \colon (L, \Phi) \to (L', \Phi')$ being commutative
triangles over $\dom$.
\end{itemize*}
Thus we have proved:
\begin{Prop}\label{monoidalstructures}
There is a strict monoidal structure $(\otimes, I)$ on $\Ff_\C$ such
that an $\otimes$-monoid structure on $({E}, \lambda, \rho)$ is the
same as an extension of the corresponding pair $(R, \Lambda)$ to a
monad over $\cod$. Dually, there is a strict monoidal structure
$(\odot, \bot)$ on $\Ff_\C$ such that a $\odot$-comonoid structure
on $({E}, \lambda, \rho)$ is the same as an extension of the
corresponding pair $(L, \Phi)$ to a comonad over $\dom$.
\end{Prop}
To relate this to \nwfs's, observe that if we take only the data
which concerns $\Rr$ in Definition \ref{natwfs} then what we have is
a monad over $\cod$ on $\Ar \C$; and likewise, taking only the data
relating to $\Ll$ gives us a comonad over $\dom$. So we can think of
$\cat{Mon}_\otimes(\Ff_\C)$ as the category of ``right halves of
\nwfs's'' and $\cat{Comon}_\odot(\Ff_\C)$ as the category of ``left
halves of \nwfs's''. Moreover, if we combine the two monoidal
structures in a simple-minded way, we nearly get enough to capture a
full \nwfs:
\begin{Prop}\label{nearlynwfs}
To give an object $({E}, \lambda, \rho)$ of $\Ff_\C$ which is
simultaneously a $\otimes$-monoid and a $\odot$-comonoid is to give:
\begin{itemize*}
\item A comonad $\mathsf L = (L, \Phi, \Sigma)$ on $\Ar \C$,
\item A monad $\mathsf R = (R, \Lambda, \Pi)$ on $\Ar \C$, and
\item A natural transformation $\Delta \colon LR \Rightarrow RL$,
\end{itemize*}
satisfying the following equalities:
\begin{align*}
\dom \cdot L & = \dom\text, & \cod \cdot L & = \dom \cdot R\text, &
\cod \cdot R & = \cod\text;\\
\dom \cdot \Phi & = 1_{\dom}\text, & \cod \cdot \Phi & = \kappa
\cdot R\text, & \dom \cdot \Lambda & = \kappa \cdot L\text, & \cod
\cdot \Lambda & = 1_{\cod}\text;\\
\text{and }\dom \cdot \Sigma & = 1_{\dom}\text, & \cod \cdot \Sigma
& = \dom \cdot \Delta\text, & \dom \cdot \Pi & = \cod \cdot
\Delta\text, & \cod \cdot \Pi & = 1_{\cod}\text.
\end{align*}
\end{Prop}
\begin{proof}
Just as in Section \ref{Sec:nwfs}, the functor $\Delta \colon LR
\Rightarrow RL$ is completely determined by the other data and the
requirements $\cod \cdot \Sigma = \dom \cdot \Delta$ and $\dom \cdot
\Pi = \cod \cdot \Delta$. Everything else follows immediately from
Proposition \ref{monoidalstructures}.
\end{proof}
Comparing this Corollary with Definition \ref{natwfs}, we see that
the only thing missing is the stipulation that $\Delta$ should be
not only a natural transformation, but also a distributive law; and
this extra requirement amounts to requiring that what we actually
have is a \emph{bialgebra} in $\Ff_\C$. For this to make sense, we
first need to show that the two monoidal structures on $\Ff_\C$
interact properly:

\begin{Prop}\label{2fold}
$(\Ff_\C, \otimes, I, \odot, \bot)$ is a strict 2-fold monoidal category.
\end{Prop}
\begin{proof}
Recall that this amounts to giving maps $m \colon \bot \otimes \bot \to \bot$, $c \colon I \to I
\odot I$ and $j \colon I \to \bot$ and a natural family of maps
\[z_{A, B, C, D} \colon ( A \odot  B) \otimes ( C \odot  D) \to ( A \otimes  C) \odot ( B \otimes  D)\]
obeying laws. Since $I$ is initial and $\bot$ is terminal in $\Ff_\C$, the maps $m, j$ and $c$ are
uniquely determined, and so we need only give the maps $z_{A, B, C, D}$, which we do directly.
Suppose that we have:
\[
A = ({E}^1, \lambda^1, \rho^1)\text, \quad B = ({E}^2, \lambda^2,
\rho^2)\text, \quad C = ({E}^3, \lambda^3, \rho^3) \quad \text{and}
\quad D = ({E}^4, \lambda^4, \rho^4)\text;
\]
then the factorisation $(A \odot B) \otimes (C \odot D)$ sends the map $f \colon X \to Y$ to
\[
X \xrightarrow{\textstyle\lambda^1_{L^2 R^{3 \odot 4} f} \cdot
\lambda^3_{L^4 f}} {E}^1L^2R^{3 \odot 4}f
\xrightarrow{\textstyle\rho^2_{R^{3 \odot 4}f} \cdot \rho^1_{L^2R^{3
\odot 4}f}} Y\text,
\]
where $R^{3 \odot 4}$ corresponds to the right-hand part of the
factorisation $({E}^3, \lambda^3, \rho^3) \odot ({E}^4, \lambda^4,
\rho^4)$. Likewise, the factorisation $(A \otimes C) \odot (B
\otimes D)$ sends $f$ to
\[
X \xrightarrow{\textstyle\lambda^1_{R^3 L^{2 \otimes 4} f} \cdot
\lambda^3_{L^{2 \otimes 4} f}} {E}^1R^3L^{2 \otimes 4}f
\xrightarrow{\textstyle\rho^2_{R^{4}f} \cdot \rho^1_{R^3L^{2 \otimes
4}f}} Y\text,
\]
where a similar meaning is attached to $L^{2 \otimes 4}$. Now, to
give $z_{A, B, C, D}$ we must give, for each such $f$, a map
${E}^1L^2R^{3 \odot 4}f \to {E}^1R^3L^{2 \otimes 4}f$, compatible
with the maps from $X$ and to $Y$ and natural in $f$. To do this,
consider the following diagram:
\[
\cd[@+1.5em]{
    {E}^3L^4f
        \ar[r]^{\lambda^2_{R^{3 \odot 4}f}}
        \ar[d]_{{E}^3(\id_X, \lambda^2_{R^4f})} &
    {E}^2R^{3 \odot 4}f
        \ar[d]^{{E}^2(\rho^3_{L^4f}, \id_Y)} \\
    {E}^3L^{2 \otimes 4}f
        \ar[r]_{\rho^3_{L^{2 \otimes 4}f}} &
    {E}^2R^{4}f\text.
 }
\]
This commutes, with both sides equal to ${E}^3L^4f \xrightarrow{\rho^3_{L^4f}} {E}^4f
\xrightarrow{\lambda^2_{R^4f}} {E}^2R^4f$. Applying $E^1$, we obtain the map
\[{E}^1\big({E}^3(\id_X, \lambda^2_{R^4f}), {E}^2(\rho^3_{L^4f}, \id_Y)\big) \colon {E}^1L^2R^{3 \odot 4}f \to {E}^1R^3L^{2 \otimes 4}f\text,\]
which we take to be the component of $z_{A, B, C, D}$ at $f$. The remaining (extensive) details are
routine.
\end{proof}

So, since our category $\Ff_\C$ bears the structure of a strict
2-fold monoidal category, we can consider bialgebras in it, in the
sense of Definition \ref{bialgdef}; and as we might hope, we have:

\begin{Prop}\label{bialgfac}
Natural weak factorisation systems on $\C$ are in bijection with bialgebras in the strict 2-fold
monoidal category $\Ff_\C$.
\end{Prop}
\begin{proof}
By Proposition \ref{nearlynwfs}, it suffices to show that, given an
object $({E}, \lambda, \rho)$ of $\Ff_\C$ which is both a
$\otimes$-monoid and a $\odot$-comonoid, the bialgebra axioms
\eqref{bialg} hold just when the equation \eqref{nwfs2} does. Now,
since $I$ is an initial object and $\bot$ is a terminal object, the
first three bialgebra axioms will always hold; whilst the fourth
holds just when the following two composites are equal for all $f
\colon X \to Y$:
\[{E}Rf \xrightarrow{\pi_f} {E}f \xrightarrow{\sigma_f} {E}Lf\]
and
\[
\cd[@C+4.5em]{ {E}Rf \ar[r]^-{\sigma_{Rf}} & {E}LRf
\ar[r]^-{{E}(\sigma_f, {E}(\sigma_f, 1_Y))} & {E}LR^{{E} \odot {E}}f
\ar[d]^{{E}({E}(1_X,
\lambda_{Rf}), {E}(\rho_{Lf}, 1_Y))} \\
{E}Lf &  {E}RLf \ar[l]^-{\pi_{Lf}} & {E}RL^{{E} \otimes {E}}f
\ar[l]^-{{E}({E}(1_X, \pi_f), \pi_f)}}
\] (where the meaning of $R^{{E} \odot {E}}$ and $L^{{E} \otimes {E}}$ is as in the proof of Proposition \ref{2fold}). Now, considering the central three maps in the latter composite, we calculate:
\begin{align*}
& {E}\big({E}(1_X, \pi_f), \pi_f\big) \cdot {E}\big({E}(1_X,
\lambda_{Rf}), {E}(\rho_{Lf}, 1_Y)\big) \cdot {E}\big(\sigma_f,
{E}(\sigma_f, 1_Y)\big)
\\ = \  & {E}\big({E}(1_X, \pi_f) \cdot {E}(1_X, \lambda_{Rf}) \cdot
\sigma_f,\, \pi_f \cdot {E}(\rho_{Lf}, 1_Y) \cdot {E}(\sigma_f,
1_Y)\big)\\ = \  &  {E}\big({E}(1_X, \pi_f \cdot \lambda_{Rf}) \cdot
\sigma_f,\, \pi_f \cdot {E}(\rho_{Lf} \cdot \sigma_f, 1_Y)\big)
\\ = \  &  {E}\big(
\sigma_f, \pi_f\big)\text.
\end{align*}
Thus the fourth bialgebra axiom holds just when $\sigma_f \cdot
\pi_f = \pi_{Lf} \cdot {E}(\sigma_f, \pi_f) \cdot \sigma_{Rf}$ for
all $f$, as required.
\end{proof}

We note in passing that we can use this characterisation theorem to
read off the correct notion of \emph{morphism} between \nwfs's on
$\C$: it is simply a map of functorial factorisations which respects
both the monad and the comonad structure. We will not make direct
use of this notion in the current paper, but it is undoubtedly
rather important, since once would expect a putative ``algebraic''
version of a \emph{Quillen model structure} to contain (amongst
other data) two \nwfs's, (trivial cofibration, fibration) and
(cofibration, trivial fibration), together with a morphism of
\nwfs's from the former to the latter.

\section{Cofibrantly generated \nwfs's: definition} \label{Sec:cgendef}
We are now ready to tackle the main topic of this paper, the definition and construction of
\emph{cofibrantly generated} \nwfs's. Recall first that a plain \wfs\ $(\ELL, \R)$ on a category
$\C$ is said to be
\defn{cofibrantly generated} if there is a set $J$ of $\ELL$-maps, called the
\defn{generating cofibrations}, such that $\R$ is precisely the class of maps with the right
lifting property with respect to each of the maps in $J$; it then follows that $\ELL$ is the class
of maps with the left lifting property with respect to each of the maps in $\R$, and so the set $J$
completely determines the factorisation system.

For example, there is a cofibrantly generated \wfs\ on the category of topological spaces whose
generating cofibrations $J$ are the inclusions $S_{n-1} \to D_n$ of the $(n-1)$-sphere into the
$n$-disc. It typifies a certain ``topological'' kind type of \wfs, where one thinks of each
generating cofibration $f \colon A \to B$ as specifying a \emph{shape} or \emph{cell} $B$ together
with the inclusion of its \emph{boundary} $A$.

When the generating cofibrations are viewed in this way, one arrives at very natural
interpretations of the two classes of the resultant w.f.s. The left class consists of retracts of
\emph{cell complexes}, which are maps $X \to Y$ obtained by a transfinite process which, starting
with $X$, iteratively picks out boundaries along which to glue in cells until arriving at $Y$. In
the example of the previous paragraph, these cell complexes directly generalise the topologist's
\emph{CW-complexes}. The right class is, of course, still determined by the right lifting property;
but given a generating cofibration $f$, we might suggestively call a lifting problem like
\begin{equation}\label{abcd}
    \cd[@-0.5em]{
        A
            \ar[r]^-{h}
            \ar[d]_f &
        C
            \ar[d]^g \\
        B
            \ar[r]_-{k} &
        D
    }
\end{equation}
a \emph{relative horn of shape $f$ in $g$}, and call a solution for this lifting problem a
\emph{filler}. Then the right class consists of those maps $g$ such that \emph{every relative horn
in $g$ has a filler}.

The purpose of this section is to develop a corresponding notion of cofibrantly generated
\emph{natural} \wfs: roughly speaking, we will say that a \nwfs\ $(\mathsf L, \mathsf R, \Delta)$
on $\C$ is cofibrantly generated by $J$ if its $\Rr$-maps are arrows $g \colon X \to Y$ in $\C$
equipped with a \emph{choice} of filler for every relative horn.

So, suppose we are given a category $\C$ and a set of maps $J$ in
it. Then we can ``algebraise'' the notion of ``having the right
lifting property with respect to $J$'': given $g \colon C \to D$ in
$\C$, we write $S_g$ for the set whose elements are $(f, g)$-lifting
problems as in \eqref{abcd} as $f$ ranges over $J$, and define
\emph{right lifting data for $g$ w.r.t.\ $J$}  to be a function
$\delta$ assigning to each lifting problem $x \in S_g$ as in
\eqref{abcd} a chosen fill-in $\delta(x) \colon B \to C$ with
$\delta(x) \cdot f = h$ and $g \cdot \delta(x) = k$. We can form
such right lifting data into a category $\ro J$, with
\begin{itemize*}
\item \textbf{Objects} being pairs $(g, \delta)$ where $g \colon C \to D$ and $\delta$ is right lifting data for $g$
with respect to $J$, and
\item \textbf{Morphisms} $(g, \delta) \to (g', \delta')$ being commutative squares
\[
    \cd[@-0.5em]{
        C
            \ar[r]^-{m}
            \ar[d]_g &
        C'
            \ar[d]^{g'} \\
        D
            \ar[r]_-{n} &
        D'
    }
\]
which commute with the right lifting data for $g$ and $g'$: that is, given an element of $S_g$, we
should get the same result from first applying $\delta$ and then postcomposing with $m$, or from
first composing on the right with the square $(m, n)$ and then applying $\delta'$.
\end{itemize*}
This category comes equipped with an obvious forgetful functor to
$\Ar \C$ which we will denote by $U_J \colon \ro J \to \Ar \C$.
\begin{Exs}\hfill\label{ex10}
\begin{itemize*}
\item Let $\C = \cat{Set}$ and let $J = \{0 \to 1\}$. Then for any map $g \colon C \to D$, we have
$S_g = D$; a typical object of $\ro J$ is a map $g \colon C \to D$ together with a map $i \colon D
\to C$ satisfying $gi = \id_D$; and a typical morphism $(g, i) \to (g', i')$ of $\ro J$ is a map
$(h, k) \colon g \to g'$ of $\Ar \C$ such that $i'k = hi$.
\item Let $\C = \cat{Set}$ and let $J = \{\inj_1 \colon 1 \to 1 + 1\}$. Then for any map $g \colon C \to D$, we have
$S_g = C \times D$; a typical object of $\ro J$ is a map $g \colon C \to D$ together with a
function $\theta \colon C \times D \to C$ satisfying $g(\theta(c, d)) = d$ for all $c, d \in D$;
and a typical morphism $(g, \theta) \to (g', \theta')$ of $\ro J$ is a map $(h, k) \colon g \to g'$
satisfying $h(\theta(c, d)) = \theta'(h(c), k(d))$ for all $c \in C$ and $d \in D$.
\item Combining the previous two, if $\C = \cat{Set}$ and $J = \{!\colon 0 \to 1, \inj_1 \colon 1 \to 1 +
1\}$, then for any map $g \colon C \to D$, we have $S_g = D + C \times D$; a typical object of $\ro
J$ is a map $g \colon C \to D$ together with a both a map $i \colon D \to C$ satisfying $gi =
\id_D$ and a map $\theta \colon C \times D \to C$ satisfying $g(\theta(c, d)) = d$ for all $c, d
\in D$.
\item Let $\C = R\text-\cat{Mod}$ be the category of modules over a commutative ring $R$, and let $J = \{0 \to R\}$. Then for any map $g \colon M \to N$, we have
$S_g = \abs N$, the underlying set of the $R$-module $N$, whilst a typical object of $\ro J$ is a
map $g \colon M \to N$ together with a mere \emph{function} $k \colon \abs N \to \abs M$ such that
$gk(n) = n$ for all $n \in N$.
\item Let $\C$ be the category of directed multigraphs, i.e., the functor category $[\cdot
\rightrightarrows \cdot, \cat{Set}]$. We write a typical object of $\C$ as \[X =
\cdl{X_a\ar@<2pt>[r]^{s} \ar@<-2pt>[r]_{t} & X_v}\] (for $a$rrows, $v$ertices, $s$ource and
$t$arget), and a typical map as $f = (f_a \colon X_a \to Y_a, f_v \colon X_v \to Y_v)$. Let $J = \{
(\bullet) \rightarrow (\bullet \to \bullet)\}$ consist of the inclusion of the graph with one
vertex as the source of the graph with one arrow and two vertices; i.e., the following map:
\[\cd{
 0 \ar@<2pt>[r]^{!} \ar@<-2pt>[r]_{!} \ar[d]_{!} & 1 \ar[d]^{\inj_1} \\
 1 \ar@<2pt>[r]^{\inj_1} \ar@<-2pt>[r]_{\inj_2} & 2\text.
}\] Now, given a map $f \colon X \to Y$ in $\C$, we have $S_{f}$
given by the pullback
\[
    \cd[@-0.5em]{
        S_{f} \pushoutcorner
            \ar[r]^-{\overline s}
            \ar[d]_{\overline f} &
        X_v
            \ar[d]^{f_v} \\
        Y_a
            \ar[r]_-{s} &
        Y_v \text;
    }
\]
thus a typical element $(a, x)$ of $S_f$ is given by an arrow $a$ in $Y$ together with a vertex $x$
of $X$ lying over its source. A typical object of $\ro J$ is given by a map $f \colon X \to Y$
equipped with a lifting of $\overline f$ through $f_a$:
\[
    \cd[@-0.5em]{ & A_x \ar[d]^{f_a} \\
        S_{f}
            \ar[r]_{\overline f}
            \ar[ur]^p &
        A_y\text,
    }
\]
which is to give, for each element $(a, x)$ of $S_f$, an arrow $p(a, x)$ of $X$ lying over $a$.
\end{itemize*}
\end{Exs}
[All of these may look like toy examples: but we have chosen them as such to give us something to
play with. If we look at something more substantial~--~$\C$ the category of simplicial sets and $J$
the set of horn inclusions, for example~--~then the most explicit we can be is that an element of
$\ro J$ as a map $g \colon C \to D$ equipped with a chosen filler for every relative horn: which is
fine but doesn't really give us anything to get our hands on.]

Now suppose that as well as a set of maps $J$, we are also given a \nwfs\ $(\Ll, \Rr, \Delta)$ on
$\C$. To say that this is cofibrantly generated by $J$ should mean that its category of $\Rr$-maps
is isomorphic to the category $\ro J$ we have just defined; however, if we take this as our
definition then we have missed out on an important subtlety: $\Rmap$ should not be isomorphic to
$\ro J$ in any old way, but in a \emph{canonical} way.

To make sense of this, we need some extra data. Observe that for a cofibrantly generated \wfs, each
map of $J$ is an $\ELL$-map: so for a cofibrantly generated \nwfs, we expect each map of $J$ to be
an $\Ll$-map. We could take this to mean that each map of $J$ admits at least one $\Ll$-map
structure, but if we are going to be consistent about our philosophy of ``algebraisation'', we
should surely take it to mean that each element of $J$ comes equipped with a \emph{chosen}
$\Ll$-map structure. Thus our additional piece of data is a factorisation
\[
 \cd[@!@-2.5em]{ J \ar[rr]^\alpha \ar[dr]_{\iota} & & \Lmap\text. \ar[dl]^{U_\Ll} \\ & \Ar \C}
\]
(Here we view $J$ as a discrete subcategory of $\Ar \C$). Using this data, we now have a
\emph{canonical} way of obtaining right lifting data w.r.t.\ $J$ from any $\Rr$-map $(g, s)$.
Indeed, $\alpha$ equips each element $f \in J$ with an $\Ll$-map structure $(f, \alpha_f)$, and so
we can solve $(f, g)$-lifting problems like \eqref{abcd} using the liftings from the \nwfs\ between
the $\Ll$-map $(f, \alpha_f)$ and the $\Rr$-map $(g, s)$. This assignation extends to a functor
$\theta$:
\[
 \cd[@!@-2.5em]{ \Rmap \ar[rr]^\theta \ar[dr]_{U_\Rr} & & \ro J\text. \ar[dl]^{U_J} \\ & \Ar \C}
\]
We now say that the \nwfs\ $(\Ll, \Rr, \Delta)$ is \defn{cofibrantly generated} by $(J, \alpha)$ if
the functor $\theta$ so defined is an isomorphism of categories: in other words, if the $\Rr$-maps
are completely determined by the lifting data that they give with respect to the generating
cofibrations $J$.

In order for this to be a sensible definition, it should be conservative over the corresponding
definition for plain \wfs's. To see this, suppose that $(\Ll, \Rr, \Delta)$ is cofibrantly
generated by $J$ and consider its underlying functorial \wfs\ $(\overline \ELL, \overline \R)$. We
recall that $\overline \ELL$ and $\overline \R$ are the respective closures under retracts of
$\ELL$, the class of maps in $\C$ admitting some $\Ll$-coalgebra structure, and $\R$, the class of
maps admitting some $\Rr$-algebra structure. In this case the category of $\Rr$-maps is isomorphic
to the category $\ro J$ of right lifting data with respect to $J$, and so a map of $\C$ lies in
$\R$ precisely when it has the right lifting property with respect to $J$. Thus the underlying
functorial \wfs\ $(\overline \ELL, \overline \R)$ is precisely the \wfs\ cofibrantly generated by
$J$.

In particular, the functorial factorisation that we construct for a cofibrantly generated \nwfs\ in
the next section gives rise to a functorial factorisation for the underlying plain \wfs: and it is
a much smaller and more tractable factorisation than one generally obtains for cofibrantly
generated \wfs's. For this reason, even the reader who feels that natural \wfs's have nothing much
to commend them over functorial \wfs's should find the following results of interest.

\section{Cofibrantly generated \nwfs's: construction}\label{Sec:cgen:construct}
\subsection{Introduction}\label{Sec:cgen:construct:intro}
Now we know what a cofibrantly generated \nwfs\ \emph{is}, we can begin to investigate the
circumstances under which we can build one. The method we use will be familiar both to topologists,
who will recognise it as a variant of Quillen's \emph{small object argument} (for a modern account
of which, see \cite{Hovey} or \cite{HH}, for example), and to category theorists, who will
recognise it as an example of the construction of the \emph{free monad on a pointed endofunctor}, a
subject treated in detail by Kelly \cite{Ke80}.

So, suppose that we are given a category $\C$ and a set $J$ of maps in it; let us work backwards
from the definition of a cofibrantly generated \nwfs\ and see if we can build one which is
generated by $J$. Our starting point is the observation that, if $\C$ is cocomplete, we can greatly
simplify the definition of our category $\ro J$. Fix a map $g \colon C \to D$ of $\C$, and consider
again the set $S_g$ of commutative squares as in \eqref{abcd}. We can view each $x \in S_g$ as a
morphism $(h_x, k_x) \colon f_x \to g$ in $\Ar \C$, and thus we can combine them into a map
\[
\spn{(h_x, k_x)}_{x \in S_g} \colon \sum_{x \in S_g} f_x \to g
\]
of $\Ar \C$; that is, a diagram
\begin{equation}
\label{singlesquare}
    \cd{
        \sum A_x
            \ar[r]^-{\left<h_x\right>}
            \ar[d]_{\sum_x f_x} &
        C
            \ar[d]^g \\
        \sum B_x
            \ar[r]_-{\left<k_x\right>} &
        D\text.
    }
\end{equation}
Now, to give right lifting data for $g$ w.r.t.\ $J$ is equivalent to
giving a diagonal fill-in for this single square. We take this
process one stage further by observing that \eqref{singlesquare}
factorises as
\begin{equation}
\label{doublesquare}
    \cd{
        \sum A_x
            \ar[r]^-{\left<h_x\right>}
            \ar[d]_{\sum_x f_x} &
        C
            \ar@{=}[r]
            \ar[d]^{\lambda^1_g} &
        C \ar[d]^g \\
        \sum B_x
            \ar[r] &
        {E}^1g
            \pullbackcorner \ar[r]_{\rho^1_g} & D\text,
    }
\end{equation}
where the left-hand square is a pushout, and that giving a fill-in
for \eqref{singlesquare} is equivalent to giving a fill-in for the
right-hand square of \eqref{doublesquare}: that is, a map $k \colon
{E}^1g \to C$ satisfying the two equalities $g \cdot k = \rho^1_g$
and $k \cdot \lambda^1_g = \id_C$. The first of these says that
there is a commutative square
\begin{equation}
\label{sugg1}
    \cd{
        {E}^1g
            \ar[r]^-{k}
            \ar[d]_{\rho^1_g} &
        C
            \ar[d]^g \\
        D
            \ar@{=}[r] &
        D\text,
    }
\end{equation}
which, if we write $R^1 \colon \Ar \C \to \Ar \C$ for the functor\footnote{We will see that this
operation really is a functor in the next section.} sending $g$ to $\rho^1_g$, corresponds to a map
$\gamma = (k, 1_D) \colon R^1g \to g$. The second equality says that
\begin{equation}
\label{sugg2}
    \cd{
        C
            \ar[d]_g
            \ar[r]^{\lambda^1_g} &
        {E}^1g
            \ar[r]^-{k}
            \ar[d]_{\rho^1_g} &
        C
            \ar[d]^g \\
        D
            \ar@{=}[r] &
        D
            \ar@{=}[r] &
        D
    } \quad = \quad
    \cd{
        C
            \ar[d]_g
            \ar@{=}[r]^{\phantom{\lambda^1_g}} &
        C
            \ar[d]^g \\
        D
            \ar@{=}[r] &
        D\text;
    }
\end{equation}
which, writing $\Lambda^1 \colon \id_{\Ar \C} \Rightarrow R^1$ for the natural
transformation\footnote{Ditto.} whose component at $g$ is $\Lambda^1_g = (\lambda^1_g, 1_D)$,
corresponds to the assertion that $\gamma \cdot \Lambda^1_g = \id_g$. Moreover, \emph{every} map
$\gamma \colon R^1g \to g$ in $\Ar \C$ satisfying this equality must arise in this way, since the
equality $\gamma \cdot \Lambda^1_g = \id_g$ forces $\gamma$ to be of the form $(k, 1_D)$. Thus we
have proven:
\begin{Prop}\label{algj}
Giving right lifting data w.r.t.\ $J$ for $g$ is equivalent to giving a map $\gamma \colon R^1g \to
g$ satisfying $\gamma \cdot \Lambda^1_g = \id_g$.
\end{Prop}
We can see this as an example of a more general concept: a \defn{pointed endofunctor} $(T, \tau)$
on a category $\K$ is a functor $T \colon \K \to \K$ together with a natural transformation $\tau
\colon \id_\K \Rightarrow T$. So $(T, \tau)$ is a ``monad without the multiplication'', and like a
monad, it gives rise to a \defn{category of algebras}, $T\text-\cat{Alg}$, with:
\begin{itemize*}
\item \textbf{Objects} being pairs $(X, x)$ where $X \in \K$ and $x \colon TX \to X$,
satisfying the unit condition $x \cdot \tau_X = \id_X$;
\item \textbf{Morphisms} $(X, x) \to (Y, y)$ being maps $f \colon X \to Y$ such that
$y \cdot Tf = f \cdot x$.
\end{itemize*}
In particular, we can consider the category $R^1$-$\cat{Alg}$ for the pointed endofunctor $(R^1,
\Lambda^1)$ above; and in this language, Proposition \ref{algj} says that objects of
$R^1\text-\cat{Alg}$ are the same as objects of $\ro J$. As one would hope, this correspondence
extends to morphisms:
\begin{Prop}\label{r1isjblob}
There is an isomorphism, commuting with the forgetful functors to $\Ar \C$, between the category
$\ro J$ of right lifting data w.r.t.\ $J$ and the category $R^1\text-\cat{Alg}$ of algebras for the
pointed endofunctor $(R^1, \Lambda^1)$ on $\Ar \C$.
\end{Prop}

Now, according to the definition in the previous section, a \nwfs\ $(\Ll, \Rr, \Delta)$ is
cofibrantly generated by $J$ if $\ro J$ is isomorphic to $\Rmap$ in a canonical way. Leaving aside
the ``in a canonical way'' part for the moment, and using the characterisation of the previous
Proposition, this means that the category $R^1\text-\cat{Alg}$ we have just defined must be
isomorphic to the category of algebras for some monad $\Rr$ on $\Ar \C$; and this monad $\Rr$ will
provide us with the right-hand side of our \nwfs

We thus are led to ask: when does such an isomorphism exist? Questions such as this are dealt with
comprehensively in \cite{Ke80}, and in this case the answer is very simple.
\begin{Prop}\label{kellyresult}
Let $(T, \tau)$ be a pointed endofunctor on a category $\K$. Then $T\text-\cat{Alg}$ is isomorphic
to the category of algebras for a monad $\Rr$ on $\K$ just when the forgetful functor $U \colon
T\text-\cat{Alg} \to \K$ has a left adjoint. In this case, $\Rr$ is called the
\defn{algebraically-free monad} on the pointed endofunctor $(T, \tau)$ and is isomorphic to the monad
generated by the adjunction $F \dashv U \colon T\text-\cat{Alg} \to \K$.
\end{Prop}
So the obvious next question is, when does this left adjoint exist? Again, \cite{Ke80} provides an
answer: if $\K$ is cocomplete and the functor $T$ is suitably ``small'', we can construct the
desired left adjoint as the colimit of a transfinite sequence. Here, ``small'' can mean something
very general, but we will only need the following two cases of it. The first is very familiar:
\begin{Defn}\label{small1}
Let $\alpha$ be a cardinal. We say that a limit ordinal $\beta$ is
\defn{$\alpha$-filtered} if, for every subset $A \subset \beta$ of cardinality $\leqslant \alpha$,
we have $\sup A < \beta$. We say that a functor $T \colon \K \to \ELL$ is \defn{$\alpha$-small} if
it preserves colimits of chains indexed by $\alpha$-filtered ordinals.
\end{Defn}
The second is a slight refinement of the first, and requires the
notion of an \defn{orthogonal factorisation system}
\cite{FreydKelly:continuous} on our category $\K$; as mentioned in
the introduction, this is given by two classes of maps $(\E, \M)$
satisfying the same axioms as a weak factorisation system, except we
strengthen the ``lifting'' property to
\begin{description}
\item[(unique lifting)] Whenever we are given a commutative square
\[
    \cd{
        A
            \ar[r]^-f
            \ar[d]_e &
        C
            \ar[d]^m \\
        B
            \ar[r]_-g &
        D
    }
\]
in $\K$, where $e \in \E$ and $m \in \M$, we can find a \emph{unique} fill-in $j \colon B \to C$
such that $mj = g$ and $je = f$.
\end{description}
Typical examples are (epi, mono) factorisations in $\cat{Set}$; and either (surjection, embedding)
or (quotient, injection) factorisations in $\cat{Top}$. Given an orthogonal factorisation system
$(\E, \M)$ on a category $\K$, we shall say that it is \defn{cowellpowered} if every object $X$ of
$\K$ possesses a mere set of isomorphism classes of $\E$-maps with domain $X$: for example, each of
the three factorisation systems just cited are cowellpowered.
\begin{Defn}\label{small2}
Let $\K$ be a category equipped with a cowellpowered factorisation system $(\E, \M)$. We say that a
functor $T \colon \K \to \ELL$ is \defn{$\alpha$-small relative to $\M$} if it preserves colimits
of chains of $\M$-maps indexed by $\alpha$-filtered ordinals.
\end{Defn}

\begin{Prop}\label{convergence} Let $(T, \tau)$ be a pointed endofunctor on a cocomplete category $\K$, such that $T$
is either $\alpha$-small or $\alpha$-small relative to $\M$ for some cowellpowered $(\E, \M)$. Then
the forgetful functor $U \colon T\text-\cat{Alg} \to \K$ has a left adjoint.
\end{Prop}
Applying this result to our pointed endofunctor $(R^1, \Lambda^1)$, we see that, as long as $R^1$
is $\alpha$-small~--~which amounts to requiring our set of generating maps $J$ to be $\alpha$-small
in a suitable sense~--~we can build the algebraically-free monad $\Rr$ on $(R^1, \Lambda^1)$.
However, we are not out of the woods yet: this approach builds a monad $\Rr$ for the right-hand
side of our putative \nwfs, but does not produce a corresponding comonad $\Ll$. Given how
intertwined the two parts of a \nwfs\ are, it may appear that we are in a somewhat hopeless
situation.

This is where the view of \nwfs's as bialgebras comes into play. The
pointed endofunctor $(R^1, \Lambda^1)$ is really another
presentation of the functorial factorisation $T = ({E}^1, \lambda^1,
\rho^1)$. This is an object of $\Ff_\C$, and in fact a \emph{pointed
object} $\tau \colon I \to T$, where the map $\tau \colon I \to T$
is the unique map from the initial object $I$. From this
perspective, building the free monad on the pointed endofunctor
$(R^1, \Lambda^1)$ is more-or-less the same thing as building the
free $\otimes$-monoid on the pointed object $(T, \tau)$ of $\Ff_\C$.

What we shall shortly see is that the functorial factorisation $T =
({E}^1, \lambda^1, \rho^1)$~--~or rather, its alternative
presentation as a pair $(L^1, \Phi^1)$~--~already admits a comonad
structure: so by Proposition \ref{monoidalstructures}, we can lift
$T$ from an object of $\Ff_\C$ to an object of
$\cat{Comon}_\odot(\Ff_\C)$. Moreover, because $\Ff_\C$ is a 2-fold
monoidal category, the $(\otimes, I)$ monoidal structure also lifts
from $\Ff_\C$ to $\cat{Comon}_\odot(\Ff_\C)$.

But now we can try to lift the free monoid construction for $(T,
\tau)$ from $\Ff_\C$ to $\cat{Comon}_\odot(\Ff_\C)$, thereby
obtaining a $\otimes$-monoid in $\cat{Comon}_\odot(\Ff_\C)$, which
is a bialgebra in $\Ff_\C$, which, by Proposition \ref{bialgfac}, is
a \nwfs\ on $\C$. Moreover, the monad for this \nwfs\ will be the
right thing~--~the algebraically-free monad $\Rr$ on $(R^1,
\Lambda^1)$~--~because the construction we used is just a lifting of
this free-monad construction.

Our plan is now as follows: first we show that our functorial factorisation $({E}^1, \lambda^1,
\rho^1)$ admits a comonad structure, and thus lifts from $\Ff_\C$ to $\cat{Comon}_\odot(\Ff_\C)$.
We then give an explicit description of the construction of the free $\otimes$-monoid on our lifted
functorial factorisation. Finally, we show that the resultant \nwfs\ really is cofibrantly
generated by $J$: which is where the ``in a canonical way'' which we laid aside earlier will be
picked back up again.

\subsection{The one-step comonad} \label{Sec:onestep} Our task in this section is to take the assignation
\begin{equation}
\label{assignation}g \colon C \to D \qquad \mapsto \qquad C
\xrightarrow{\lambda^1_g} {E}^1g \xrightarrow{\rho^1_g} D
\end{equation} of the previous section and
show that it gives us a functorial factorisation $({E}^1, \lambda^1, \rho^1)$ on $\C$ for which the
corresponding pair $(L^1, \Phi^1)$ has a natural extension to a comonad $\Ll^1 = (L^1, \Phi^1,
\Sigma^1)$. This functorial factorisation will be very familiar to readers who know Quillen's small
object argument for \wfs's: it provides the ``iterative step'' by which one transfinitely
constructs factorisations. The comonad $\Ll^1 =(L^1, \Phi^1, \Sigma^1)$ extending it will play a
similar role in the construction of natural \wfs's, and thus we christen it the ``one-step
comonad''. It turns out to have a very satisfactory universal property:
\begin{Prop}\label{goodunivprop}
$\Ll^1$ is the free ``comonad over $\dom$'' generated by $J$, in that there are bijections, natural
in $\Ll \in \cat{Comon}_\odot(\Ff_\C)$, between morphisms $\Ll^1 \to \Ll$ of
$\cat{Comon}_\odot(\Ff_\C)$ and morphisms
\[
 \cd[@!@-2.5em]{ J \ar[rr] \ar[dr]_{\iota} & & \Lmap \ar[dl]^{U_{\Ll}} \\ & \Ar \C}
\]
of $\cat{Cat} / \Ar \C$.
\end{Prop}
%
%
Both this Proposition and the construction of $\Ll^1$ which we are
about to given can be deduced from the fact that $\Ll^1$ is a
\emph{density comonad} in a certain 2-category. The notion of
density comonad embodies the idea of a comonad being ``freely
generated'' by an arrow: in this case, by the arrow $\iota \colon J
\to \Ar \C$ exhibiting $J$ as a discrete subcategory of $\Ar \C$.
Setting up the theory to explain this here would lead us too far
afield, and instead we defer this task to the Appendix. What we give
in the remainder of this section is the explicit description of what
this abstract framework yields.

So let us return to our contemplation of equation
\eqref{assignation}, which we recall arose from the following
process:
\[
    \cd{
        \sum_{x \in S_g} A_x
            \ar[r]^-{\left<h_x\right>}
            \ar[d]_{\sum_x f_x} &
        C
            \ar[d]^g \\
        \sum_{x \in S_g} B_x
            \ar[r]_-{\left<k_x\right>} &
        D
    }
\quad \rightsquigarrow \quad
    \cd{
        \sum_{x \in S_g} A_x
            \ar[r]^-{\left<h_x\right>}
            \ar[d]_{\sum_x f_x} &
        C
            \ar@{=}[r]
            \ar[d]^{\lambda^1_g} &
        C \ar[d]^g \\
        \sum_{x \in S_g} B_x
            \ar[r] &
        {E}^1g
            \pullbackcorner \ar[r]_{\rho^1_g} & D\text.
    }
\]
We want to make the assignation $g \mapsto L^1g$ into the object part of a functor, for which we
must give the value of $L^1$ on a morphism $\gamma = (h, k) \colon g \to g'$ of $\Ar \C$. We do
this by first making the assignation $g \mapsto Kg := \sum_{x \in S_g} f_x$ into the object part of
a functor, whose value on a map $\gamma \colon g \to g'$ of $\Ar \C$ is given as follows. Observe
that postcomposition with $\gamma$ induces a function
\begin{align*}
S_\gamma \colon S_g & \to S_{g'}\\
(f \xrightarrow{\delta} g) & \mapsto (f \xrightarrow{\delta} g \xrightarrow{\gamma} g')\text.
\end{align*}
And thus we take $K \gamma$ to be:
\[
K\gamma = \spn{\inj_{S_\gamma(x)}} \colon \sum_{x \in S_g} f_x \to \sum_{y \in S_{g'}} f_y\text.
\]
Thus we have a functor $K \colon \Ar \C \to \Ar \C$ for which the maps
\[
\phi_g = \spn{(h_x, k_x)}_{x \in S_g} \colon \sum_{x \in S_g} f_x = Kg \to g
\]
become the components of a natural transformation $\phi \colon K \Rightarrow \id_{\Ar \C}$. We will
now use the functor $K$ to give the value of $L^1$ on morphisms of $\Ar \C$. Indeed, given such a
morphism $\gamma = (h, k) \colon g \to g'$, we have the following diagram, whose left-hand face is
$K\gamma$ and whose top face is the domain part of a naturality square for $\phi$:
\[
\cd[@!@-3em@R-1em]{
    \sum_{x \in S_g} A_x
        \ar[dr]
        \ar[dd]_{Kg}
        \ar[rr]^{\spn{h_x}} & &
    C
        \ar[dr]^h \\ &
    \sum_{y \in S_{g'}} A_y
        \ar[dd]^{Kg'}
        \ar[rr]^{\spn{h_y}} & &
    C' \\
    \sum_{x \in S_g} B_x
        \ar[dr]
        \\ &
    \sum_{y \in S_{g'}} B_y
}
\]
Pushing out the rear face gives us $L^1g$, pushing out the front face gives us $L^1g'$, and so we
take $L^1\gamma$ to be the induced map from the rear to the front of the right-hand face. Now we
see that the diagram
\begin{equation}    \cd{
        \sum_{x \in S_g} A_x
            \ar[r]^-{\left<h_x\right>}
            \ar[d]_{\sum_x f_x} &
        C
            \ar@{=}[r]
            \ar[d]^{\lambda^1_g} &
        C \ar[d]^g \\
        \sum_{x \in S_g} B_x
            \ar[r] &
        {E}^1g
            \pullbackcorner \ar[r]_{\rho^1_g} & D\text,
    }
\end{equation}
viewed as a pair of maps $\epsilon_g \colon Kg \to L^1g$ and $\Phi^1_g \colon L^1g \to g$ in $\Ar
\C$, gives us the components of natural transformations $\epsilon \colon K \Rightarrow L^1$ and
$\Phi^1 \colon L^1 \Rightarrow \id_{\Ar \C}$, satisfying $\Phi^1 \cdot \epsilon = \phi$. In
particular, we have a copointed endofunctor $(L^1, \Phi^1)$ over $\dom$ which corresponds to a
functorial factorisation $({E}^1, \lambda^1, \rho^1)$, as claimed.

We now show that $(L^1, \Phi^1)$ can be extended to a comonad $\Ll^1 = (L^1, \Phi^1, \Sigma^1)$,
for which we must give maps $\Sigma^1_g \colon L^1g \rightarrow L^1L^1g$ over $\dom$, which we
write as:
\[    \cd{
        C
            \ar@{=}[r]
            \ar[d]_{\lambda^1_g} &
        C
            \ar[d]^{\lambda^1_{L^1g}} \\
        {E}^1g
            \ar[r]_-{\sigma^1_g} & {E}^1L^1g\text.
    }
\]
We obtain these maps as follows: we have a function
\begin{align*}
\psi_g \colon S_g & \to S_{L^1g} \\
x & \mapsto (f_x \xrightarrow{\inj_x} Kg \xrightarrow{\epsilon_g} L^1g)\text,
\end{align*}
and so can form the map
\[
\delta_g = \spn{\inj_{\psi_g(x)}} \colon \sum_{x \in S_g} f_x \to \sum_{y \in S_{L^1g}} f_y\text.
\]
Now we have the following diagram, whose left-hand face is $\delta_g$:
\[
\cd[@!@-3em@R-1.5em]{
    \sum_{x \in S_g} A_x
        \ar[dr]
        \ar[dd]_{Kg}
        \ar[rr]^{\spn{h_x}} & &
    C
        \ar@{=}[dr] \\ &
    \sum_{y \in S_{L^1g}} A_y
        \ar[dd]^{KL^1g}
        \ar[rr]^{\spn{h_y}} & &
    C \\
    \sum_{x \in S_g} B_x
        \ar[dr]
        \\ &
    \sum_{y \in S_{L^1g}} B_y
}
\]
Pushing out the rear face gives us ${E}^1g$, pushing out the front
face gives us ${E}^1L^1g$, and the induced map along the
bottom-right diagonal we take to be the value of $\sigma^1_g$.

This completes our description of the one-step comonad $\Ll^1$: but before moving on to consider
how we can use it to build a \nwfs, we should discuss what $\Ll^1$-coalgebras are. Let us write
$\Ll^1\text-\cat{Map}$ for the category of such, and call its objects $\Ll^1$-maps; as in Section
\ref{Sec:nwfs}, we write them as pairs $(f, s)$ where $f \colon X \to Y$ and $s \colon Y \to
{E}^1f$. Now, every $\Ll^1$-map will induce an $\Ll$-map in the \nwfs\ generated by $J$, and the
intuition is that they should be just those $\Ll$-maps which can be obtained using only one step's
worth of ``glueing on cells''. We will make this intuition precise in Proposition
\ref{characteriseonestep}, where we will characterise $\Ll^1$-maps as (certain) retracts of
pushouts of coproducts of the generating cofibrations; but for now the following examples should
give a good feel for what happens.

\begin{Exs}\hfill
\begin{itemize*}
\item When $\C = \cat{Set}$ and $J = \{0 \to 1\}$, we obtain the
one-step factorisation
\[
g \colon X \to Y \qquad \mapsto \qquad X \xrightarrow{\inj_1} X + Y \xrightarrow{\spn{g, \id}}
Y\text,
\]
and $\sigma^1_g = \spn{\inj_1, \inj_3} \colon X + Y \to X + (X + Y)$. In this case, an $\Ll^1$-map
$(f, s)$ is a map $f \colon X \to Y$ which is an injection: this comonad is ``property-like'' in
that any map can carry at most \emph{one} coalgebra structure. A morphism of $\Ll^1$-coalgebras
$(f, s) \to (f', s')$ is given by a map
\begin{equation}\label{eqnexample}
    \cd[@-0.5em]{
        X
            \ar[r]^h
            \ar[d]_{f} &
        X'
            \ar[d]^{f'} \\
        Y
            \ar[r]_-{k} &
        Y'
    }
\end{equation}
such that $k$ maps $Y \setminus f(X)$ into $Y' \setminus f'(X')$.
\item When $\C = \cat{Set}$ and $J = \{\inj_1 \colon 1 \to 1 + 1\}$, we obtain the
one-step factorisation
\[
g \colon X \to Y \qquad \mapsto \qquad X \xrightarrow{\inj_1} X + X \times Y \xrightarrow{\spn{g,
\pi_2}} Y\text,
\]
and $\sigma^1_g = \spn{\inj_1, \psi} \colon X + X \times Y \to X + X \times (X + X \times Y)$,
where $\psi(c, d) = (c, (c, d))$. An $\Ll^1$-map $(f, s)$ is given by an injection $f \colon X \to
Y$ together with a map $i \colon Y \setminus f(X) \to X$ (saying ``where the extra elements were
attached''). A morphism of $\Ll^1$-coalgebras $(h, k) \colon f \to f'$ is a map as in
\eqref{eqnexample} for which the following diagram commutes:
\[
    \cd[@-0.5em]{
        Y \setminus f(X)
            \ar[r]^-i
            \ar[d]_{k} &
        X
            \ar[d]^{h} \\
        Y' \setminus f'(X')
            \ar[r]_-{i'} &
        X'\text.
    }
\]
\item When $\C = \cat{Set}$ and $J = \{! \colon 0 \to 1, \inj_1 \colon 1 \to 1 + 1\}$, we obtain the
one-step factorisation
\[
g \colon X \to Y \qquad \mapsto \qquad X \xrightarrow{\inj_1} X + Y + X \times Y
\xrightarrow{\spn{g, \id_Y, \pi_2}} Y\text,
\]
and we leave the description of $\sigma^1_g$ to the reader. An $\Ll^1$-map $(f, s)$ is given by an
injection $f \colon X \to Y$, a partition of $Y \setminus f(X)$ into disjoint subsets $Y_1$ and
$Y_2$, and a function $i \colon Y_2 \to X$. The elements of $Y_1$ correspond to elements attached
via $! \colon 0 \to 1$, whilst the elements of $Y_2$ correspond to elements attached via $\inj_1
\colon 1 \to 1 + 1$, and $i$ tells us how these elements were attached.
\item When $\C = R\text-\cat{Mod}$ and $J = \{0 \to R\}$, we obtain the one-step factorisation
\[
g \colon M \to N \qquad \mapsto \qquad M \xrightarrow{\inj_1} M \oplus FN \xrightarrow{\spn{g,
\textsf{ev}}} N\text.
\]
Here, $FN$ is the free $R$-module on the underlying set of $N$ and $\textsf{ev}$ is the obvious map
from there to $N$. If we write the generators of $FN$ as $\{x_n\}_{n \in \abs N}$, then the
comultiplication map is given by
\begin{align*}
\sigma^1_g \colon M \oplus FN & \to M \oplus F(M \oplus FN)\\
(m, 0) & \mapsto (m, 0)\\
(0, x_n) & \mapsto (0, x_{x_n}) \quad \text{for $n \in \abs N$}
\end{align*}
extended linearly. An $\Ll^1$-coalgebra $(f, s)$ is a map $f \colon
M \to N$ which is a direct summand inclusion (i.e., the canonical
map $N \to f(M) \oplus N / f(M)$ is a bijection) together with a
subset $X \subset N / f(M)$ and an isomorphism $\theta \colon N /
f(M) \cong FX$. A map of $\Ll^1$-coalgebras is a morphism $(h, k)
\colon g \to g'$ such that $k(N / f(M)) \subset N' / f'(M')$, such
that $k(X) \subset X'$, and such that the diagram
\[\cd{
N / f(M) \ar[r]^k \ar[d]_{\theta} & N' / f'(M') \ar[d]^\theta \\
FX \ar[r]_{Fk} & FX' }\] commutes.
\item When $\C$ is the category of directed graphs and $J =
\{(\bullet) \to (\bullet \to \bullet)\}$, the one-step factorisation
is given by:
\[\cd{
 X_a \ar@<2pt>[d]^{t} \ar@<-2pt>[d]_{s} \ar[r]^-{f_a} & Y_a \ar@<2pt>[d]^{t} \ar@<-2pt>[d]_{s} \\
 X_v \ar[r]_{f_v}  & Y_v
} \qquad \mapsto \qquad \cd[@C+2em]{
 X_a \ar@<2pt>[d]^{t} \ar@<-2pt>[d]_{s} \ar[r]^-{\inj_1} &
 X_a + S_f \ar@<2pt>[d]^{\spn{\inj_1 \cdot t, \inj_2}} \ar@<-2pt>[d]_{\spn{\inj_1 \cdot s, \inj_1 \cdot \overline s}} \ar[r]^{\spn{f_a, \overline f}} &
 Y_a \ar@<2pt>[d]^{t} \ar@<-2pt>[d]_{s} \\
 X_v \ar[r]_-{\inj_1}  & X_v + S_f \ar[r]_{\spn{f_v, t \overline f}} & Y_v\text. }
\]
Again, we omit the description of $\sigma^1$. In this case, an $\Ll^1$-map $(f, s)$ is given by a
map $f \colon X \to Y$ such that $f_a \colon X_a \to Y_a$ and $f_v \colon X_v \to Y_v$ are
injections and such that $t\big(Y_a \setminus f_a(X_a)\big) = Y_v \setminus f_v(X_v)$. This is
another ``property-like'' comonad, and comparison with the second example is instructive: despite
the similarities, that example was not ``property-like'', and it is the extra structure borne by
the category of directed graphs relative to the category of sets which is responsible for this
difference.
\end{itemize*}
\end{Exs}

\subsection{Iterating the one-step monad}\label{iteratingtheonestep}
Now that we have defined and described the one-step comonad, we are ready to use it to build a
n.w.f.s.
\begin{Defn}\label{nwfsgenj}
Let $\C$ be a cocomplete category, $J$ a set of maps in $\C$, and $\Ll^1$ the one-step comonad
corresponding to $J$. As in the discussion at the end of Section \ref{Sec:cgen:construct:intro}, we
may view $\Ll^1$ as a pointed object $(\Ll^1, \tau)$ of $\cat{Comon} := \cat{Comon}_\odot(\Ff_\C)$,
and we write $(\Ll, \eta, \mu)$ for the free $\otimes$-monoid\footnote{For a formal definition of
which, see Definition \ref{freemonoiddefn}.} on the pointed object $(\Ll^1, \tau)$~--~if it
exists~--~and call the corresponding \nwfs\ $(\Ll, \Rr, \Delta)$ the
\defn{\nwfs\ generated by $J$}.
\end{Defn}

In terms of the $\Rr$-maps, this process is a suitably refined way
of forming the \emph{free monad} on the pointed endofunctor $(R^1,
\Lambda^1)$ which corresponds to $\Ll^1$, which is the right thing
to do because it yields a category of $\Rr$-maps which is consistent
with the requirement that our \nwfs\ should be generated by $J$.
However, we can also give a natural interpretation of what we are
doing in terms of the $\Ll$-maps.

Indeed, we shall see in the next section that $\Ll^1$-maps are closed under every possible
operation we might like them to be closed under \emph{except composition}: in which terms, we can
see the process of constructing our \nwfs\ from the one-step comonad as ``closing off $\Ll^1$ under
composition''. To make this intuition slightly less vague, we must examine in more detail the
process which assigns to an object $\Ll \in \cat{Comon}$~--~viewed as a comonad over $\dom$~--~its
category of coalgebras; in keeping with our previous notation, we will write this assignation as
$\Ll \mapsto \Ll\text-\cat{Map}$.

The first observation is that we can make this into a functor $\G \colon \cat{Comon} \to \cat{Cat}
/ \Ar \C$: indeed, a morphism $\alpha \colon \Ll \to \Ll'$ in $\cat{Comon}$ is a map of the
underlying functorial factorisations $\alpha \colon ({E}, \lambda, \rho) \to ({E}', \lambda',
\rho')$ which is also a map of comonads $\Ll \to \Ll'$; thus it
induces a morphism
\[
 \cd[@!@-2.5em]{ \Ll\text-\cat{Map} \ar[rr]^{\alpha_\ast} \ar[dr]_{U_{\Ll}} & & \Ll'\text-\cat{Map}\text. \ar[dl]^{U_{\Ll'}} \\ & \Ar \C}
\]
Explicitly, this sends an $\Ll$-algebra $(f, s)$ to the $\Ll'$-algebra $(f, \alpha_f \cdot s)$: so
$\alpha_\ast$ witnesses that ``every $\Ll$-map is an $\Ll'$-map''.\footnote{The sharp-eyed reader
will have spotted that we implicitly used this functor $\G$ in the statement of Proposition
\ref{goodunivprop}.} We can now ask how this functor $\G$ interacts with the monoidal structure on
$\cat{Comon}$. The unit is straightforward: $I \in \cat{Comon}$ corresponds to the comonad which
sends a map $f \colon X \to Y$ to $\id_X \colon X \to X$, and its category of coalgebras is
precisely the full subcategory of $\Ar \C$ whose objects are the isomorphisms.

The multiplication $\otimes$ is more interesting: from an $\Ll$-map $(f, s) \colon X \to Y$ and an
$\Ll'$-map $(g, t) \colon Y \to Z$, we can obtain an $(\Ll' \otimes \Ll)$-map structure on $gf
\colon X \to Z$. We will prove this formally in Section \ref{Sec:compositionalprops}, where it
becomes part of the statement that the functor $\G$ is \emph{lax monoidal} with respect to a
suitably defined ``compositional'' monoidal structure on $\cat{Cat} / \Ar \C$.

Given this, we can see that if we have a \emph{monoid} $(\Ll, \eta, \mu)$ in $\cat{Comon}$, then
its category of coalgebras will be closed under composition: since from a pair of $\Ll$-maps $f
\colon X \to Y$ and $g \colon Y \to Z$, we obtain an $\Ll \otimes \Ll$-map $gf \colon X \to Z$;
applying $\mu_\ast \colon (\Ll \otimes \Ll)\text-\cat{Map} \to \Lmap$ to which gives us an
$\Ll$-map structure on $gf$. In particular, forming the free monoid on a pointed object of
$\cat{Comon}$ can be seen as freely closing off its category of coalgebras under composition.

A final perspective on what we are doing, and perhaps the most convincing, comes from the
combination of Proposition \ref{goodunivprop} and Definition \ref{nwfsgenj}:
\begin{Prop}\label{goodunivprop2} If the \nwfs\ $(\Ll, \Rr, \Delta)$
generated by a set of maps $J$ exists, then it is the free \nwfs\ on $J$, in the sense that there
are bijections, natural in $(\Ll', \Rr', \Delta')$, between morphisms $(\Ll, \Rr, \Delta) \to
(\Ll', \Rr', \Delta')$ of $\cat{Bialg}(\Ff_\C)$ and morphisms
\[
 \cd[@!@-2.5em]{ J \ar[rr] \ar[dr]_{\iota} & & \Ll'\text-\cat{Map} \ar[dl]^{U_{\Ll}} \\ & \Ar \C}
\]
of $\cat{Cat} / \Ar \C$.
\end{Prop}

Before we examine the details of the free monoid construction (drawing on \cite{Ke80} once more),
let us motivate why it takes the form it does by answering the following question: if $\Ll^1$-maps
corresponds to doing one step's worth of glueing, then for what comonad $\Ll^2$ do $\Ll^2$-maps
correspond to doing \emph{two} step's worth of glueing? The obvious first guess, $\Ll^1 \otimes
\Ll^1$, turns out to be not quite right. For we observe that there are \emph{two} copies of $\Ll^1$
embedded inside $\Ll^1 \otimes \Ll^1$, via the maps:
\[\Ll^1 \xrightarrow{\tau \otimes \Ll^1} \Ll^1 \otimes \Ll^1 \xleftarrow{\Ll^1 \otimes \tau} \Ll^1\text.\]
These two embeddings correspond to taking one step's worth of glueing on cells and either
prepending or postpending it with one step's worth of doing nothing. But surely we would like to
identify these two: we do not really want to record how long we waited around before glueing some
cells on, after all. Thus, more correctly, we should take $\Ll^2$ to be the coequaliser:
\[\cd{\Ll^1 \ar@<2pt>[r]^-{\tau \otimes \Ll^1} \ar@<-2pt>[r]_-{\Ll^1 \otimes \tau} & \Ll^1 \otimes \Ll^1 \ar[r] & \Ll^2\text.}\]
This is exactly what we will do, forming each $\Ll^{\alpha^+}$ as a suitable coequaliser of $\Ll^1
\otimes \Ll^{\alpha}$ and then taking $\Ll$ to be the colimit of this sequence. This is somewhat
different from a simple-minded generalisation of the small object argument, which would correspond
to forming the colimit of a suitably long sequence of the form:
\[\Ll^1 \xrightarrow{\tau \otimes \Ll^1} \Ll^1 \otimes \Ll^1
\xrightarrow{\tau \otimes \Ll^1 \otimes \Ll^1} \Ll^1 \otimes \Ll^1 \otimes \Ll^1 \to \cdots\] in
$\cat{Comon}$. Here, the problem we have just described with respect to $\Ll^1 \otimes \Ll^1$ is
present and, indeed, drastically multiplied, giving us a plethora of different ways of glueing on
the same cells depending on how much waiting around we choose to do; and this is surely not what we
want.

\subsubsection{The theory}
In this section, we give a brief summary of the material we need from \cite{Ke80} pertaining to the
construction of a free monoid on a pointed object; except where noted, everything in this section
can be found in that paper. We will site our summary in an arbitrary cocomplete\footnote{For the
moment, we make no assumptions about the preservation of colimits in $\V$ by any of the functors $A
\otimes (\thg)$ or $(\thg) \otimes A$; in particular, we do not assume that $\V$ is either left or
right closed.} monoidal category $(\V, \otimes, I)$, because the degeneracy of the particular
example we are interested in (where the unit $I$ is also the initial object) sometimes makes it
harder to see what is going on.

We have, of course, the familiar notions of \emph{monoid} and \emph{monoid map} in $\V$, whilst, as
we have mentioned before, a \emph{pointed object} $(T, \tau)$ of $\V$ is an object $T \in \V$
equipped with a map $\tau \colon I \to T$; finally, by a \emph{map of pointed objects} $\alpha
\colon (S, \sigma) \to (T, \tau)$ we mean a map $\alpha \colon S \to T$ satisfying $\alpha \sigma =
\tau$.

\begin{Defn}\label{freemonoiddefn}
Let $(T, \tau)$ be a pointed object of $\V$: then the \defn{free monoid} on $(T, \tau)$ is a monoid
$(U, \eta, \mu)$ together with a map of pointed objects $\chi \colon (T, \tau) \to (U, \eta)$ such
that precomposition with $\chi$ induces a isomorphism, natural in $V$, between maps of monoids $(U,
\eta, \mu) \to (V, \eta', \mu')$ and maps of pointed objects $(T, \tau) \to (V, \eta')$.
\end{Defn}

Now, to build the free monoid on $(T, \tau)$, it often suffices to construct the ``free object with
an action by $T$''. To make this precise, we consider the category $T\text-\cat{Mod}$ of
\defn{modules} for a pointed object $(T, \tau)$, with
\begin{itemize*}
\item \textbf{Objects} being pairs $(X, \ x \colon T \otimes X \to X)$ in $\V$ satisfying $x \cdot (\tau \otimes X) =
1_X$;
\item \textbf{Morphisms} $f \colon (X, x) \to (Y, y)$ being maps $f \colon X \to Y$ satisfying $f \cdot x = y \cdot (T \otimes
f)$.
\end{itemize*}
Observe that if $(X, x)$ is a $T$-module and $A \in \V$, then so is $(X \otimes A, x \otimes A)$,
and that this assignation extends to a ``right action'' of the monoidal category $\V$ on
$T\text-\cat{Mod}$; that is, a functor
\begin{align*}
\mathord{\star} \colon T\text-\cat{Mod} \times \V &\to T\text-\cat{Mod} \\
\big((X, x), A\big) &\mapsto (X \otimes A, x \otimes A)
\end{align*}
satisfying the two usual laws for a right action, but weakened up to coherent isomorphism.
\begin{Prop}\label{freemonoid}
Let $(T, \tau)$ be a pointed object in $\V$. If there is a $T$-module $(X, x)$ such that $(X, x)
\star (\thg) \colon \V \to T\text-\cat{Mod}$ is left adjoint to the forgetful functor $U \colon
T\text-\cat{Mod} \to \V$, then $X$ is the underlying object of the free monoid on $(T, \tau)$.
\end{Prop}
\begin{proof}
The isomorphism $(X \otimes I, x \otimes I) = (X, x) \star I \to (X, x)$ in $T\text-\cat{Mod}$
corresponds under adjunction to a map $\eta \colon I \to X$, whilst the map $1_X \colon X \to X$
corresponds under adjunction to a map $(X, x) \star X = (X \otimes X, x \otimes X) \to (X, x)$ in
$T\text-\cat{Mod}$ underlying which is a map $\mu \colon X \otimes X \to X$ in $\V$. We refer the
reader to \cite{Ke80} for the remaining details.
\end{proof}
The hypotheses of this Proposition are equivalent to saying that, firstly, the free $T$-module on
$I$ exists and is given by $(X, x)$, and secondly, that the free $T$-module on any other $W \in \V$
can be obtained ``pointwise'' from this as $(W \otimes X, x \otimes X)$. Our route to satisfying
these hypotheses will be to attempt to construct a left adjoint to $U \colon T\text-\cat{Mod} \to
\V$ using a certain transfinite construction, which we describe in Definition \ref{freealg}. This
process may or not converge when applied to an object $X$; if it does, we say that the free
$T$-module on $X$ \emph{exists constructively}. All we need to know about this construction for the
moment is that it is obtained as the colimit of a certain transfinite sequence~--~the \emph{free
module sequence} for $X$~--~each stage of which is built using tensor products and connected
colimits of the previous stages.

\begin{Prop}
Let $(T, \tau)$ be a pointed object in $\V$ and suppose that the free $T$-module on $I$ exists
constructively and is given by $(X, x)$. If each functor $(\thg) \otimes A \colon \V \to \V$
preserves connected colimits then the forgetful functor $T\text-\cat{Alg} \to \V$ has a left
adjoint given by $(X, x) \star (\thg) \colon \V \to T\text-\cat{Alg}$.
\end{Prop}
\begin{proof}
Because the functor $(\thg) \otimes A \colon \V \to \V$ preserves connected colimits, the free
algebra sequence for $A$ is obtained, up to isomorphism, as $(\thg) \otimes A$ of the free algebra
sequence for $I$. In particular, we can take the free algebra on $A$ to be $(X \otimes A, x \otimes
A)$.
\end{proof}
Thus, under the assumption that each functor $(\thg) \otimes A \colon \V \to \V$ preserves
connected colimits, we can build the free monoid on $(T, \tau)$ whenever the free $T$-module on $I$
exists constructively; the only thing remaining is to describe what ``exists constructively''
means. We first fix some notation concerning transfinite sequences:

\begin{Defn}
Let $\kappa$ be a regular inaccessible cardinal. We write $\cat{On}$ for the well-ordered set of
ordinals smaller than $\kappa$, viewed as a posetal category. By a \emph{transfinite sequence} $X$
in $\V$, we mean a functor $X \colon \cat{On} \to \V$; we write the image of an ordinal $\alpha$ as
$X_\alpha$ and the image of the inequality $\alpha \leqslant \beta$ as $X_{\alpha, \beta} \colon
X_\alpha \to X_\beta$.
\end{Defn}
\noindent And now:
\begin{Defn}\label{freealg}
Let $(T, \tau)$ be a pointed object of $\V$, and let $A \in \V$. By the \defn{free $T$-module
sequence} for $A$, we mean the following transfinite sequence $X \colon \cat{On} \to \V$, which we
will construct simultaneously with a family of maps $\sigma_\alpha \colon T \otimes X_\alpha \to
X_{\alpha^+}$, natural in $\alpha$ and satisfying $\sigma_\alpha \cdot (\tau \otimes X_{\alpha^+})
= X_{\alpha, \alpha^+}$. Note that this last condition determines the value of the connecting maps
$X_{\alpha, \alpha^+}$ for any ordinal $\alpha \in \cat{On}$, and thus we need only give
$X_{\alpha, \beta}$ when $\beta$ is a limit ordinal.
\begin{itemize}
\item $X_0 = A$, $X_1 = T \otimes A$, $\sigma_0 = 1_{T \otimes A}$;
\item For a successor ordinal $\alpha^+$, we give $X_{\alpha^{++}}$ and $\sigma_{\alpha^+}$ by the following
coequaliser diagram:
\[\cd{ &
 X_{\alpha^+}
  \ar@/^1em/[dr]^-{\tau \otimes X_{\alpha^+}}
  \\
 T \otimes X_\alpha
  \ar@/^1em/[ur]^-{\sigma_\alpha}
  \ar@/_1em/[dr]_-{T \otimes \tau \otimes X_\alpha} & &
 T \otimes X_{\alpha^+}
  \ar[r]^-{\sigma_{\alpha^+}} &
 X_{\alpha^{++}}\text.
  \\ &
 T \otimes T \otimes X_{\alpha}
  \ar@/_1em/[ur]_-{T \otimes \sigma_{\alpha}}
}\]
\item For a non-zero limit ordinal $\gamma$, we give $X_\gamma$ by
$\colim_{\alpha < \gamma} X_\alpha$, with connecting maps $X_{\alpha, \gamma}$ given by the
injections into the colimit. We give $X_{\gamma^+}$ and $\sigma_\gamma$ by the following
coequaliser diagram:
\[
\cd{
 &  \colim X_{\alpha^+} = X_\gamma
  \ar[dr]^{\tau \otimes X_\gamma} \\
 \colim (T \otimes X_{\alpha})
  \ar[ur]^{\colim \sigma_\alpha\ }
  \ar[rr]_{\text{can}} & &
 T \otimes \colim X_{\alpha} = T \otimes X_{\gamma}
  \ar[r]_-{\sigma_{\gamma}} &
 X_{\gamma^+}
 }
\]
where ``can'' is the canonical map induced by the cocone $T \otimes X_\alpha \to T \otimes \colim
X_\alpha$.
\end{itemize}
We say that this sequence \emph{converges at $\alpha$} if $X_{\alpha, \alpha^+} \colon X_\alpha \to
X_{\alpha^+}$ is invertible for some $\alpha \in \cat{On}$; it then follows that $X_{\alpha,
\beta}$ is invertible for \emph{every} $\beta > \alpha$.
\end{Defn}
\begin{Prop}
Let $(T, \tau)$ be a pointed object of $\V$ and $A \in \V$ for which the free $T$-module sequence
for $A$ converges at $\alpha$. Then the free $T$-module on $A$ exists constructively, and is given
by $X_\alpha$ equipped with the algebra map
\[T \otimes X_\alpha \xrightarrow{\theta_\alpha} X_{\alpha^+} \xrightarrow{X_{\alpha,
\alpha^+}^{-1}}  X_\alpha\text.\] The universal map from $A$ is given by $X_{0, \alpha} \colon A
\to X_\alpha$.
\end{Prop}
To ensure that $T$-module sequences do converge,  we require, as in Proposition \ref{convergence},
some smallness assumption on $T$. Recall that we gave two such notions in Definitions \ref{small1}
and \ref{small2},  which we can reuse by stipulating that an object $T \in \V$ is
\defn{$\alpha$-small} or \defn{$\alpha$-small relative to $\M$} just when the corresponding endofunctor $T \otimes (\thg) \colon \V
\to \V$ is so. Kelly shows that if $T$ is $\alpha$-small, then every free $T$-module sequence will
converge at $\alpha$, whilst if $T$ is $\alpha$-small relative to some $\M$, then every free
$T$-module sequence will converge, though not necessarily at $\alpha$. Combining all of the above,
we have the following result:
\begin{Cor}\label{monoidresult}
Let $\V$ be a cocomplete monoidal category such that each functor $(\thg) \otimes A \colon \V \to
\V$ preserves connected colimits, and let $(T, \tau)$ be a pointed object of $\V$ such that $T$ is
either $\alpha$-small or $\alpha$-small relative to $\M$ for some cowellpowered $(\E, \M)$. Then
the free monoid on $(T, \tau)$ exists: its underlying object is given by the free $T$-module on
$I$, as constructed in Definition \ref{freealg}, and its multiplication and unit are given as in
Proposition \ref{freemonoid}.
\end{Cor}

\subsubsection{The practice} We are now ready to apply this machinery to the case of interest to
us. The result we are aiming for is the following:
\begin{Prop}\label{mainresult}
Let $\{f_j \colon A_j \to B_j\}_{j \in J}$ be a set of maps in a cocomplete category $\C$ such that
either every $\C(A_j, \thg) \colon \C \to \cat{Set}$ is $\alpha_j$-small or there is a
cowellpowered factorisation system $(\E, \M)$ on $\C$ such that every $\C(A_j, \thg) \colon \C \to
\cat{Set}$ is $\alpha_j$-small with respect to $\M$. Then the \nwfs\ generated by $J$ exists.
\end{Prop}
Our method, of course, will be to apply Corollary \ref{monoidresult}, and so we will obtain the
\nwfs\ generated by a set of maps $J$ as the colimiting value $\ELL$ of the free $\Ll^1$-module
sequence on $I$ in $\cat{Comon}$. Thus we will have:
\begin{itemize*}
\item For each $\alpha \in \cat{On}$, an object $\Ll^\alpha \in \cat{Comon}$ (where $\Ll^0 = I$ and $\Ll^1$ is itself);
\item For each $\alpha \in \cat{On}$, an ``action morphism'' $\theta_\alpha \colon \Ll^1 \otimes \Ll^\alpha \to
\Ll^{\alpha^+}$;
\item For each $\alpha < \beta \in \cat{On}$, a connecting morphism $X_{\alpha, \beta} \colon \Ll^\alpha \to
\Ll^\beta$,
\end{itemize*}
In particular, the connecting map from $\Ll^1$ into the colimiting value $\Ll$ is the ``universal
map of pointed objects'' $\chi \colon \Ll^1 \to \Ll$ of Definition \ref{freemonoiddefn}; we will
make use of this map in the next subsection.

In terms of the intuitive description given at the start of Section \ref{iteratingtheonestep}, the
$\Ll^\alpha$-maps correspond to maps given by ``at most $\alpha$ steps of glueing on cells''; the
morphisms $X_{\alpha, \beta}$ witness the fact that every $\Ll^\alpha$-map is an $\Ll^\beta$-map
when $\alpha < \beta$; and the morphisms $\theta_\alpha$ attest to the fact that anything we can do
in most $\alpha$ steps of glueing followed by a single further step of glueing, we can do in at
most $\alpha^+$ steps of glueing.

Now, although the hypotheses of Proposition \ref{mainresult} should be fairly unsurprising to
anyone used to the small object argument for plain \wfs's, the proof that they are sufficient is
surprisingly technical. Before we give it, let us see how widely these hypotheses are satisfied.
Firstly, they hold for any set of maps $J$ whatsoever in a \emph{locally presentable} category
$\C$. We recall that a category $\C$ is
\defn{locally $\kappa$-presentable} for some regular cardinal $\kappa$ if it is cocomplete and has
a set $\S$ of objects such that we have both
\begin{description}
\item[(density)] Every $X \in \C$ is a canonical colimit of elements of $\S$; and
\item[(smallness)] For each $A \in \S$, $\C(A, \thg) \colon \C \to \cat{Set}$ preserves $\kappa$-filtered
colimits,
\end{description}
and that $\C$ is \defn{locally presentable} if it is locally $\kappa$-presentable for some
$\kappa$. Now, for any object $X$ in a locally presentable category, there exists some cardinal
$\alpha$ such that $\C(X, \thg)$ is $\alpha$-small, and thus the hypotheses of the Proposition will
always hold.

Any sufficiently ``algebraic'' category is locally presentable: for example, $\cat{Set}$,
$\cat{Ab}$, $R\text-\cat{Mod}$ (modules over a ring $R$), $\cat{Ch}(R)$ (chain complexes over a
ring $R$) or any presheaf category, in particular the category $\cat{SSet}$ of simplicial sets, are
all locally presentable and in fact locally \emph{finitely} presentable, i.e., locally
$\omega$-presentable. Examples of locally presentable categories that are not locally finitely
presentable include the category $\cat{Shv}(\C)$ of sheaves for a site $\C$ and the category
$\omega\text-\cat{CPO}$ of $\omega$-complete partially ordered sets. For more on the theory of
locally presentable categories, one might refer to \cite{Borceaux} or \cite{LPC}.

There are also important examples of non-locally presentable
categories in which our Proposition can be applied: for example,
when $J$ is any set of maps whatsoever in the category $\cat{Top}$
of topological spaces or the category $\cat{Haus}$ of compact
Hausdorff topological spaces. In these cases, we need our more
refined notion of smallness: namely, relative to a cowellpowered
orthogonal factorisation system $(\E, \M)$, which for both of the
named categories, we can take to be given by $\E = $ projections and
$\M = $ subspace embeddings. We have that \emph{for any topological
space $X$ of cardinality $< \alpha$, $\cat{Top}(X, \thg)$ is
$\alpha$-small relative to the subspace embeddings}: see, for
example, \cite[Section 2.4]{Hovey} for a proof. The same holds when
$\cat{Top}$ is replaced by $\cat{Haus}$, and thus the hypotheses of
the Proposition will be satisfied, for any set $J$ of maps
whatsoever, in either of these two categories.

The rest of this section will be devoted to proving Proposition \ref{mainresult}. Our method, of
course, will be to apply Corollary \ref{monoidresult}, and so we need to check that all the
relevant hypotheses are satisfied. Apart from the smallness condition, this amounts to showing
that:
\begin{Prop}
The category $\cat{Comon}$ is cocomplete, and for each $A \in
\cat{Comon}$, the functor $(\thg) \otimes A \colon \cat{Comon} \to
\cat{Comon}$ preserves connected colimits.
\end{Prop}
\begin{proof}
We know that $\C$ is cocomplete, and thus so also is $[\Ar \C, \C]$. Now the category $\Ff_\C$ can
be obtained by slicing $[\Ar \C, \C]$ over the object $\cod$, and then coslicing this under the
object $(\kappa \colon \dom \Rightarrow \cod)$: thus $\Ff_\C$ is cocomplete. Moreover, the
forgetful functor $\cat{Comon}_\odot(\Ff_\C) \to \Ff_\C$ creates colimits, and so
$\cat{Comon}_\odot(\Ff_\C)$ is also cocomplete. For the second part, consider the composite
forgetful functor
\[U \colon \cat{Comon}_\odot(\Ff_\C) \to \Ff_\C \to [\Ar \C, \Ar \C]\text,\]
where the second arrow sends a functorial factorisation $(F, \lambda, \rho)$ to the corresponding
functor $R \colon \Ar \C \to \Ar \C$. Now, the first part creates all colimits whilst the second
creates connected colimits, and thus $U$ creates connected colimits. Moreover, it sends the
$(\otimes, I)$ monoidal structure on $\cat{Comon}$ to the monoidal structure on $[\Ar \C, \Ar \C]$
given by composition; and so, given $A \in \cat{Comon}$ with underlying object $(F, \lambda, \rho)$
in $\Ff_\C$, the following diagram commutes:
\[
\cd[@C+1em]{
    \cat{Comon}
      \ar[r]^{(\thg) \otimes A} \ar[d]_U &
    \cat{Comon} \ar[d]^U \\
    [\Ar \C, \Ar \C]
      \ar[r]_{(\thg) \cdot R} &
    [\Ar \C, \Ar \C]\text.
}
\]
But $(\thg) \cdot R$ preserves connected colimits (indeed, all colimits), and so the result
follows.
\end{proof}
We will not directly satisfy the smallness condition for $\Ll^1$, but will give sufficient
conditions for free $\Ll^1$-module sequences to converge nonetheless:
\begin{Prop}\label{reducetosmallercase}
Let $\Ll^1 \in \cat{Comon}$ be the one-step comonad generated by a set of maps $J$ in $\C$, and let
$(R^1, \Lambda^1)$ be the corresponding pointed endofunctor of $\Ar \C$. Then free $\Ll^1$-modules
exist constructively whenever the functor $R^1$ is either $\alpha$-small or $\alpha$-small with
respect to $\M$ for some cowellpowered $(\E, \M)$ on $\Ar \C$.
\end{Prop}
\begin{proof}
It suffices to show that the free $\Ll^1$-module on $I$ exists constructively, since all other free
$\Ll^1$-modules are computed ``pointwise'' from this one. We do this by considering the forgetful
functor $U \colon \cat{Comon} \to [\Ar \C, \Ar \C]$ from the proof of the previous Proposition.
Observe first that $U$ sends the pointed object $(\Ll^1, \tau)$ to the pointed endofunctor $(R^1,
\Lambda^1)$. Moreover, because $U$ both preserves the monoidal structure and creates connected
colimits, it sends the free $\Ll^1$-module sequence on $I$ to the \emph{free monad} sequence for
the pointed endofunctor $(R^1, \Lambda^1)$; we will not define this formally~--~the reader can find
it in \cite{Ke80}~--~but it should suffice if we say that it looks like Definition \ref{freealg}
with the tensor product symbols removed.

We now observe that $U$ reflects isomorphisms, so that the free $\Ll^1$-module sequence on $I$ will
converge whenever the free monad sequence for $(R^1, \Lambda^1)$ does, and so in particular
whenever $R^1$ satisfies either of the given smallness hypotheses.
\end{proof}
Note that it also follows from this proof that the underlying monad of the \nwfs\ generated from
$\Ll^1$ will be the algebraically-free monad on the pointed endofunctor $(R^1, \Lambda^1)$, which
is what we wanted. All that remains is to show that under the hypotheses of Proposition
\ref{mainresult}, the corresponding $R^1$ is suitably small.

Before we do so, we note that we can weaken the smallness
requirement on $R^1$ slightly: indeed, if we take the free monad
sequence for $R^1$~--~which is given by a chain of endofunctors of
$\Ar \C$ and natural transformations between them~--~and evaluate it
at any object $f \in \Ar \C$, then the resulting chain in $\Ar \C$
will be \emph{constant} in its codomain part, so that it suffices
for $R^1$ to be small with respect to chains of this sort. To put
this another way, observe that, since $R^1$ is a functor over
$\cod$, it restricts to functors $\res{R^1} X \colon \C / X \to \C /
X$ on each slice category of $\C$, and it is sufficient for each of
these functors to be small.

\begin{Prop}
Suppose that $\Ll^1 \in \cat{Comon}$ is generated by a set of morphisms $\{f_j \colon A_j \to
B_j\}_{j \in J}$ such that each $\C(A_j, \thg) \colon \C \to \cat{Set}$ is $\alpha_j$-small. Then
there is a regular cardinal $\alpha$ such that each functor $\res{R^1} X \colon \C / X \to \C / X$
is $\alpha$-small.
\end{Prop}
\begin{proof}
Take $\alpha$ to be a regular cardinal which is larger than all of the $\alpha_j$, so that each
$\C(A_j, \thg)  \colon \C \to \cat{Set}$ is $\alpha$-small. We aim to show that each $\res{R^1} X
\colon \C / X \to \C / X$ is $\alpha$-small. Our first observation is that, because the functor
$\Sigma_X \colon \C / X \to \C$ which forgets the projection onto $X$ creates colimits, it will
suffice to show that each composite functor $\Sigma_X \cdot \res{R^1} X \colon \C / X \to \C$ is
$\alpha$-small. Now, $\Sigma_X \cdot \res{R^1} X$ is just the restriction of $F^1 \colon \Ar \C \to
\C$ along the inclusion $i \colon \C / X \to \Ar \C$, which we can write as the composite:
\[\C / X \xrightarrow i \Ar \C \xrightarrow{L^1} \Ar \C \xrightarrow{\cod} \C\text.\]
So, since $\cod \colon \Ar \C \to \C$ preserves all colimits, we
will be done if we can prove that the composite $L^1 \cdot i \colon
\C / X \to \Ar \C$ is $\alpha$-small. To do this, we first show that
the composite $K \cdot i \colon \C / X \to \Ar \C$ is
$\alpha$-small, where we recall from Section \ref{Sec:onestep} that
$K \colon \Ar \C \to \Ar \C$ is the functor defined by \[Kg =
\sum_{x \in S_g} f_x\text.\] So let $\beta$ be an $\alpha$-filtered
ordinal and let $g_{(\thg)} \colon \beta \to \C / X$: we must show
that $K \cdot i$ preserves the colimit of this sequence, or
equivalently, since $i$ preserves connected colimits, that $K$
preserves the colimit of $i \cdot g_{(\thg)} \colon \beta \to \Ar
\C$. To do this, observe first that for any $f \colon A \to B$ and
$g \colon C \to X$ in $\C$, we have
\[\Ar \C(f, g) \cong \sum_{k \in \C(B, X)} \C / X(\Sigma_k
f, g)\] where $\Sigma_k \colon \C / B \to \C / X$ is the functor
given by postcomposition with $k$. Thus for any $f \colon A \to B$
such that that the functor $\C(A, \thg) \colon \C \to \cat{Set}$ is
$\alpha$-small, we have:
\begin{align*}
  \Ar \C(f, \colim_i g_i)
  & \cong \sum_{k \in \C(B, X)} \C / X(\Sigma_k f, \colim_i g_i) \\
  & \cong \sum_{k \in \C(B, X)} \colim_i \C / X(\Sigma_k f, g_i) \\
  & \cong \colim_i \sum_{k \in \C(B, X)} \C / X(\Sigma_k f, g_i) \\
  & = \colim_i \Ar \C (f, g_i)\text,
\end{align*}
where the step from the first to the second line follows from the
smallness of $\C(A, \thg)$ and the fact that the forgetful functor
$\C / X \to \C$ creates colimits. Since for each $f_j \in J$, we
know that $\C(A_j, \thg)$ is $\alpha$-small, we now deduce:
\begin{align*}
S_{\colim_i g_i} & = \sum_{j \in J} \Ar \C(f_j, \colim_i g_i) \\ & \cong \sum_{j \in J} \colim_i \Ar \C(f_j, g_i) \\
& \cong \colim_i \sum_{j \in J} \Ar \C(f_j, g_i) = \colim_i S_{g_i}
\end{align*}
And so we have:
\begin{align*}
K(\colim_i g_i) & = \sum_{x \in S_{\colim_i g_i}} \!\!\! f_x \\ & \cong \sum_{x \in \colim_i S_{g_i}} \!\!\! f_x \\
& \cong \colim_i \sum_{x \in S_{g_i}} \!\!\! f_x = \colim_i
Kg_i\text.
\end{align*}
Thus $K \cdot i \colon \C / X \to \Ar \C$ is $\alpha$-small; it
remains to deduce that the same is true of $L^1 \cdot i$. What we
will in fact prove is the stronger statement that $L^1$ preserves
any colimit which $K$ does. We first recall that $L^1g$ is obtained
from $Kg$ by pushing out along the domain of the canonical map
$\phi_g \colon Kg \to g$ given by
\[\phi_g = \spn{(h_x, k_x)}_{x \in S_g} \colon \sum_{x \in S_g} f_x \to g\text,\]
which we shall write as $L^1g = (h_g)_\ast (Kg)$. Now suppose that
$K$ preserves the colimit of some $F \colon \A \to \Ar \C$. Then on
the one hand, we have
\[\colim\big(L^1(Fa)\big) = \colim\big((h_{Fa})_\ast(KFa)\big) \cong (\colim
h_{Fa})_\ast (\colim KFa)\] since colimits commute with pushouts. On
the other, we have
\[L^1 (\colim Fa) = (h_{\colim Fa})_\ast (K \colim Fa)\text.\]
But we have the following commutative diagram:
\[\cd[@C+1.5em]{
 \colim KFa \ar[r]^-{\colim
\phi_{Fa}} \ar[d]_{\text{can}} &
  \colim Fa \ar@{=}[d] \\
 K \colim Fa \ar[r]_-{\phi_{\colim Fa}} &
 \colim Fa\text,
}\] and thus, since can is an isomorphism and hence already a
pushout, we deduce that
\begin{align*}
\colim\big(L^1(Fa)\big) & = (\colim h_{Fa})_\ast( \colim KFa) \\
& \cong (h_{\colim Fa})_\ast(\dom\ \text{can})_\ast(\colim KFa) \\
& \cong (h_{\colim Fa})_\ast(K \colim Fa) \\
& = L^1 (\colim Fa)
\end{align*}
as desired.
\end{proof}

In order to state the corresponding result for our more refined form of smallness, we first need
the following straightforward fact: if a category $\C$ comes equipped with a cowellpowered
orthogonal factorisation system $(\E, \M)$ then there is a corresponding cowellpowered
factorisation system of the same name on each slice category $\C / X$, whose $\E$-maps and
$\M$-maps are precisely those which become so when one forgets the projection down to $X$.

\begin{Prop}
Let $(\E, \M)$ be a cowellpowered orthogonal factorisation system on
$\C$, and let $\Ll^1 \in \cat{Comon}$ be generated as before by a
set of morphisms $\{f_j \colon A_j \to B_j\}_{j \in J}$ such that
each $\C(A_j, \thg) \colon \C \to \cat{Set}$ is $\alpha_j$-small
with respect to $\M$. Then there is a regular cardinal $\alpha$ such
that each functor $\res{R^1} X \colon \C / X \to \C / X$ is
$\alpha$-small with respect to $\M$.
\end{Prop}

The proof of the following is identical to the proof of the previous Proposition; and with it we
have the final ingredient to complete the proof of Proposition \ref{mainresult}.

\begin{Exs}\hfill
\begin{itemize*}
\item When $\C = \cat{Set}$ and $J = \{0 \to 1\}$, the \nwfs\ generated by $J$ is
the same as the one-step factorisation $\Ll^1$ we constructed before. Essentially, this is because
there is no \emph{need} for more than one step's worth of ``glueing on cells'', because such
glueings do not create any new boundaries into which further cells can be glued. So the underlying
functorial factorisation is
\[
g \colon X \to Y \qquad \mapsto \qquad X \xrightarrow{\inj_1} X + Y \xrightarrow{\spn{g, \id}}
Y\text;
\]
as before, and we have seen that the $\Ll$-coalgebras are precisely the injections. For the monad
part $\Rr$, we have $\pi_g = \spn{\inj_1, \inj_2, \inj_2} \colon X + Y + Y \to X + Y$, but we do
not need to describe the $\Rr$-algebras, for this example or any of the following ones, because
they are precisely the elements of $\ro J$ that we described in Examples \ref{ex10}: which is as we
would hope, since this was the whole point of setting up all this machinery!
\item When $\C = \cat{Set}$ and $J = \{\inj_1 \colon 1 \to 1 + 1\}$, the \nwfs\ generated by $J$
has the underlying functorial factorisation
\[
g \colon X \to Y \qquad \mapsto \qquad X \xrightarrow{\lambda_g} X \times Y^\ast
\xrightarrow{\rho_g} Y\text,
\]
where $Y^\ast$ is the free monoid on $Y$, whose elements are (possibly empty) lists $(y_1, \dots,
y_n)$ of elements of $Y$; we thus write elements of $X \times Y^\ast$ as $(x, y_1, \dots, y_n)$ for
some $n \geqslant 0$. Now $\lambda_g$ sends $x$ to $(x)$ and $\rho_g$ is given by:
\begin{align*}
\rho_g(x) &= g(x)\\
\rho_g(x, y_1, \dots, y_n) &= y_n \qquad \text{for all $n > 0$.}
\end{align*}
The map $\sigma_g \colon X \times Y^\ast \to X \times (X \times Y^\ast)^\ast$ sends $(x, y_1,
\dots, y_n)$ to the element
\[\big(x, (x, y_1), (x, y_1, y_2), \dots, (x, y_1, \dots, y_n)\big)\]
of $X \times (X \times Y^\ast)^\ast$, whilst the map $\pi_g \colon X
\times Y^\ast \times Y^\ast \to X \times Y^\ast$ is given by
\[\pi_g(x, (y_1, \dots, y_k), (y_{k+1}, \dots, y_n)) = (x, y_1,
\dots, y_n)\text.\] An $\Ll$-map is given by an injection $f \colon X \to Y$, together with a
partition of $Y$ into disjoint subsets $Y_0, Y_1, Y_2, \dots$, where $Y_0 = f(X)$, and for each $n
\geqslant 0$ a map $i_n \colon Y_{n+1} \to Y_n$. We can view such maps as specifying an $X$-indexed
family of well-founded trees, whose roots are labelled by the elements of $X$ and whose other nodes
are labelled by the elements of $Y \setminus f(X)$. A morphism of $\Ll$-maps is given by a map $(h,
k) \colon f \to f'$ which respects the partitions of $Y$ and $Y'$ and commutes with the attaching
maps $i_n$ and $i'_n$; in terms of trees, this amounts to giving a function $h \colon X \to X'$
together with a $X$-indexed family of height-preserving morphisms from the tree labeled by $x$ to
the tree labelled by $h(x)$.

\item When $\C = \cat{Set}$ and $J = \{! \colon 0 \to 1, \inj_1 \colon 1 \to 1 + 1\}$,
the \nwfs\ generated by $J$ has functorial factorisation
\[
g \colon X \to Y \qquad \mapsto \qquad X \xrightarrow{\lambda_g} (X
+ Y) \times Y^\ast \xrightarrow{\rho_g} Y\text;
\]
where with the same conventions as before, $\lambda_g$ sends $x$ to
$(x)$ whilst $\rho_g$ is given by:
\begin{align*}
\rho_g(x) &= g(x) \qquad \text{for $x \in X$,}\\
\rho_g(y) &= y \qquad \text{for $y \in Y$,}\\
\rho_g(\star, y_1, \dots, y_n) &= y_n \qquad \text{for all $n > 0$ and $\star \in X+Y$.}
\end{align*}
An $\Ll$-map is given by an injection $f \colon X \to Y$, a
partition of $Y$ into disjoint subsets $Y_0, Z_0, Y_1, Z_1, \dots$,
where $Y_0 = f(X)$, and for each $n \geqslant 0$, functions $i_n
\colon Y_{n+1} \to Y_n$ and $j_n \colon Z_{n+1} \to Z_n$. In terms
of trees, we can see this as specifying a family of trees, such that
every node is labelled by an element of $Y$ and such that elements
of $f(X)$ only ever label the roots of trees. We can give a
description of morphisms of $\Ll$-maps in similar terms.
\item When $\C = R\text-\cat{Mod}$ and $J = \{0 \to R\}$, we are in
the same situation as in the first example: the \nwfs\ generated by $J$ coincides with the one-step
comonad, and for the same reason~--~that we can always glue on all the cells we want in only one
step.
\item When $\C$ is the category of directed graphs and $J =
\{(\bullet) \to (\bullet \to \bullet)\}$, the \nwfs\ generated by $J$ has the following functorial
factorisation. Given a map $g \colon X \to Y$, the directed graph $Eg$ has vertices of two sorts:
\begin{itemize*}
\item Vertices $x \in X_v$, and
\item Sequences $(x, b_1, y_1, b_2, y_2, \dots, b_n, y_n)$, where $x$ is a vertex of $X$, each $y_i$ is a vertex of
$Y$ and each $b_i$ is an arrow of $Y$, satisfying $s(b_i) = y_{i-1}$ and $t(b_i) = y_i$ for each $i
\geqslant 0$ (with the convention that $y_0 = f(x)$),
\end{itemize*}
and arrows of two sorts:
\begin{itemize*}
\item Arrows $a \in X_v$, and
\item Sequences $(x, b_1, y_1, b_2, y_2, \dots, y_{n-1}, b_n)$ as above, but omitting the final $y_n$.
\end{itemize*}
The source and target of an arrow $a \in X_v$ are given by $s(a)$ and $t(a)$, whilst the source and
target of an arrow $(x, b_1, y_1, \dots, y_{n-1}, b_n)$ are given by:
\begin{align*}
s(x, b_1, y_1, \dots, y_{n-1}, b_n) &= (x, b_1, y_1, \dots, b_{n-1}, y_{n-1})\\
t(x, b_1, y_1, \dots, y_{n-1}, b_n) &= \big(x, b_1, y_1, \dots, b_n, t(b_n)\big)\text.
\end{align*}
The map $\lambda_g \colon X \to Eg$ is the obvious inclusion, whilst the map $\rho_g \colon Eg \to
Y$ is given by
\begin{align*}
\rho_g(x) &= g(x) \text{ for $x$ a vertex of $X$;}\\
\rho_g(x, b_1, y_1, b_2, y_2, \dots, b_n, y_n) &= y_n\text;\\
\rho_g(a) &= g(a) \text{ for $a$ an arrow of $X$, and}\\
\rho_g(x, b_1, y_1, \dots, y_{n-1}, b_n) &= b_n\text.
\end{align*}
Skipping over the description of $\sigma$ and $\pi$, which the reader should be able to figure out
by now, we observe that once again the comonad $\Ll$ is ``property-like.'' This time, a map $f
\colon X \to Y$ is an $\Ll$-map just when both $f_a$ and $f_v$ are injections and we can partition
$Y_a$ into sets $A_0, A_1, A_2, \dots$ and $Y_v$ into sets $V_0, V_1, \dots$ satisfying the
following properties:
\begin{itemize*}
\item $A_0 = f_a(X_a)$ and $V_0 = f_v(X_v)$;
\item $s(A_i) \subset V_{i-1}$ for $i \geqslant 1$, and
\item $t(A_i) = V_{i}$ for $i \geqslant 1$.
\end{itemize*}
\end{itemize*}
\end{Exs}

\subsection{Cofibrant generation}\label{Sec:cofibrantgeneration}
We have one final loose end to tie up in this section: we must show that the \nwfs\ $(\Ll, \Rr,
\Delta)$ generated by a set of maps $J$ is in fact \emph{cofibrantly generated} by $J$, in the
sense of Section \ref{Sec:cgendef}.

In order for this to make sense, we need to specify one addition
piece of data, namely how we want to view our generating
cofibrations $J$ as $\Ll$-maps. The map by which we do this is, in
the language of Proposition \ref{goodunivprop2}, the universal map
exhibiting $(\Ll, \Rr, \Delta)$ as the free \nwfs\ on $J$. To
construct it explicitly, we first lift through the category of
$\Ll^1$-maps:
\begin{Prop}\label{canonicalstructure}
Every generating cofibration $f \in J$ carries a canonical structure $\Ll^1$-map structure; in
other words, we have a lifting
\[
 \cd[@!@-2.5em]{ J \ar[rr]^{\alpha} \ar[dr]_{\iota} & & \Ll^1\text-\cat{Map}\text. \ar[dl]^{U_{\Ll^1}} \\ & \Ar \C}
\]
\end{Prop}
\begin{proof}
If $f$ is a generating cofibration, then we have the element $i =
(\id_f \colon f \to f)$ in $S_f$, and so can make $f$ into an
$\Ll^1$-coalgebra $(f, \alpha_f)$ by taking $\alpha_f$ to be the
codomain part of the morphism
\[f \xrightarrow{\inj_i} Kf \xrightarrow{\epsilon_f} L^1f\]
of $\Ar \C$. The ``canonicity'' of this lifting amounts to the fact
that it is \emph{another} universal map, this time the one
exhibiting $\Ll^1$ as the free ``comonad over $\dom$'' on $J$ in the
sense of Proposition \ref{goodunivprop}.
\end{proof}

We now use the fact that ``every $\Ll^1$-map is an $\Ll$-map'': more
formally, we obtain from the universal map $\chi \colon \Ll^1 \to
\Ll$ of Definition \ref{freemonoiddefn} a functor $\Ll^1$-$\cat{Map}
\to \Lmap$, and hence a lifting of $\beta \colon J \to \Lmap$ given
by the following composite:
\[
 \cd[@!@-2.5em@C+3em]{ J \ar[r]^-{\alpha} \ar[dr]_{\iota} & \Ll^1\text-\cat{Map} \ar[r]^{\chi_\ast} \ar[d]^{U_{\Ll^1}} &
 \Lmap \ar[dl]^{U_{\Ll}} \\ & \Ar \C\text.}
\]
Concretely, if we write $\chi_f \colon E^1f \to Ef$ for the
underlying components of $\chi \colon \Ll^1 \to \Ll$, then we have
$\beta(f) = (f, \chi_f \cdot \alpha_f)$. Our goal now is the
following result:
\begin{Prop}
Let $(\Ll, \Rr, \Delta)$ be the \nwfs\ generated by a set of maps $J$. Then $(\Ll, \Rr, \Delta)$ is
cofibrantly generated by $(J, \beta)$.
\end{Prop}
\begin{proof}
Recall that this means that there is a canonical isomorphism between the category $\ro J$ of right
lifting data with respect to $J$ and the category $R^1\text-\cat{Alg}$ of algebras for the pointed
endofunctor $(R^1, \Lambda^1)$ corresponding to $\Ll^1$. In order to show this, we must examine the
relationship between the categories $R^1\text-\cat{Alg}$ and $\Rmap$. First we note that underlying
the universal map of pointed objects $\chi \colon \Ll^1 \to \Ll$ is a map of pointed endofunctors
$(R^1, \Lambda^1) \to (R, \Lambda)$, which induces a functor
\begin{align*}
\chi^\ast \colon \Rmap & \to R^1\text-\cat{Alg}\\
(f \colon X \to Y, s \colon Ef \to X) & \mapsto (f \colon X \to Y, s \cdot \chi_f \colon E^1f \to
X)\text.
\end{align*}
But by the proof of Proposition \ref{reducetosmallercase}, we know that $\Rr$ is the
algebraically-free monad on $(R^1, \Lambda^1)$, and thus that $\Rmap$ is isomorphic to
$R^1\text-\cat{Alg}$; which stated more carefully says that the functor $\chi^\ast$ is an
isomorphism of categories. We now have the following situation:
\[
 \cd[@!@-2.5em]{ \Rmap \ar[rr]^\theta \ar[dr]_{\chi^\ast} & & \ro J\text. \ar@{<-}[dl]^{\psi} \\ & R^1\text-\cat{Alg}}
\]
where $\psi$ is the isomorphism of categories of Proposition \ref{r1isjblob} and where $\theta$ is
the canonical map of Section \ref{Sec:cgendef}. We know that both the diagonal arrows are
isomorphisms, and want to conclude that the horizontal arrow is an isomorphism: so if we can show
that the diagram commutes, we will be done, and we can do this by direct calculation. First, the
upper side: if we are given an $\Rr$-map $(g \colon C \to D, s \colon Eg \to C)$ and an element $x
\in S_g$, that is an $(f, g)$-lifting problem
\[
    \cd{
        A
            \ar[r]^-h
            \ar[d]_{f} &
        C
            \ar[d]^{g} \\
        B
            \ar[r]_-k &
        D\text,
    }
\]
for some $f \in J$, then $\theta(g, s)$ solves it by taking the codomain part of the morphism
\[f \xrightarrow{\inj_i} Kf \xrightarrow{\epsilon_f} L^1f\text,\]
which is a map $B \to E^1f$, and composing it with the morphism
\[E^1f \xrightarrow{\chi_f} Ef \xrightarrow{E(h, k)} Eg \xrightarrow{s} C\]
to obtain a morphism $j \colon B \to C$. For the lower side, if we are given the same $\Rr$-map
$(g, s)$ and lifting problem $x \in S_g$, then $\psi\chi^\ast(g, s)$ solves it by taking the
codomain part of the morphism
\[f \xrightarrow{\inj_x} Kg \xrightarrow{\epsilon_g} L^1g\text,\]
which is a map $B \to E^1f$, and composing it with the morphism
\[E^1g \xrightarrow{\chi_g} Eg \xrightarrow{s} C\text.\] But since the following diagram commutes:
\[\cd{
    & Kf \ar[r]^{\epsilon_f} \ar[d]^{K(h, k)} & L^1f \ar[d]^{L^1(h, k)} \\
    f \ar[ur]^{\inj_i} \ar[r]_{\inj_x} & Kg \ar[r]_{\epsilon_g} & L^1g\text,
}\] this latter is the same as taking the codomain part of $\epsilon_f \cdot \inj_i$ and composing
it with the morphism
\[E^1f \xrightarrow{E^1(h, k)} E^1g \xrightarrow{\chi_g} Eg \xrightarrow{s} C\text,\]
which by naturality of $\chi$ is the same as $s \cdot E(h, k) \cdot \chi_f$. Thus we have $\theta =
\psi \chi^\ast$ as claimed, and so $\theta$ is an isomorphism as desired.
\end{proof}

\section{Properties of $\Ll$-maps and $\Rr$-maps for \nwfs's}\label{Sec:props}
In this section, we broaden our attention from cofibrantly generated \nwfs's to \nwfs's in general.
Our concern will be to enumerate the closure properties that the categories of $\Ll$-maps and
$\Rr$-maps for a \nwfs\ have. As we go along, we will apply our results in the cofibrantly
generated case and see to what extent they allow us to give a characterisation of the $\Ll$-maps.
The answer is not wholly satisfactory: we do not achieve such a neat result as we have for plain
\wfs's, but we can come close.

In fact, most of the properties we are about to exhibit do not even require a full \nwfs, but only
\emph{half} of one: either a comonad over $\dom$ or a monad over $\cod$. The particular example we
should bear in mind is the one-step comonad $\Ll^1$ for a cofibrantly generated \nwfs, and for this
reason it is the comonad case that we will consider, though of course everything we do can be
straightforwardly dualised.

\subsection{Basic properties}\label{Sec:basicprops}
Suppose we are given a comonad $\Ll = (L, \Phi, \Sigma)$ over $\cod$ on $\Ar \C$, to which we apply
our usual conventions, writing $L(f \colon X \to Y)$ as $\lambda_f \colon X \to Ef$, and so on. We
begin simply:

\begin{Prop}\label{closedundercolimits}
The category $\Lmap$ is \defn{closed under colimits} in $\Ar \C$, in that the forgetful functor
$U_\Ll \colon \Lmap \to \Ar \C$ creates colimits.
\end{Prop}
\begin{proof}
Because $\Lmap$ is the category of coalgebras for a comonad on $\Ar \C$.
\end{proof}
\begin{Prop}
The category $\Lmap$ \defn{contains the isomorphisms}, in that every isomorphism $f \colon X \to Y$
in $\C$ can be equipped with a unique $\Ll$-coalgebra structure.
\end{Prop}
\begin{proof}
Given an isomorphism $f \colon X \to Y$, we make it into an $\Ll$-coalgebra $(f, s)$ by taking $s =
\lambda_f \cdot f^{-1}$; easy verification shows that this satisfies the coalgebra axioms.
Conversely, if $s \colon Y \to Ef$ makes $(f, s)$ into an $\Ll$-coalgebra, then from the coalgebra
axiom $s \cdot f = \lambda_f$ we deduce that $s = \lambda_f \cdot f^{-1}$.
\end{proof}
A more abstract view of this last Proposition is available: viewing $\Ll$ as an object of the
category $\cat{Comon} = \cat{Comon}_\odot(\Ff_\C)$, as in the previous section, we know that there
is a unique map $\tau \colon I \to \Ll$ from the initial object inducing a functor
$I\text-\cat{Map} \to \Lmap$. So ``every $I$-map is an $\Ll$-map'', and the $I$-maps are precisely
the isomorphisms.
\begin{Cor}\label{tech}
Suppose we are given $\Ll$-coalgebras $(f, s)$ and $(g, t)$, where $f \colon X \to Y$ an
isomorphism. Then every map $(h, k) \colon f \to g$ of $\Ar \C$ lifts to a map of $\Ll$-coalgebras
$(f, s) \to (g, t)$.
\end{Cor}
\begin{proof}
We must verify that $tk = E(h, k)s$. By the previous Proposition, $s = \lambda_f f^{-1}$ so that
$E(h, k)s = E(h, k)\lambda_f f^{-1} = \lambda_g hf^{-1} = tghf^{-1} = tkff^{-1} = tk$ as required.
\end{proof}
Though this last result might seem somewhat technical, when combined with Proposition
\ref{closedundercolimits} it already implies that $\Ll$-maps are closed under pushout. To give a
precise meaning to ``closure under pushout'', we need the notion of a \emph{cocartesian lifting}.

Suppose we are given an arbitrary functor $U \colon \A \to \C$, an arrow $f \colon C \to D$ in $\C$
and an object $X \in \A$ lying over $C$. Then a \defn{cocartesian lifting of $f$ at $X$} is a
universal way of turning $X$ into an object $f_\ast X$ lying over $C'$: it is given by such an
object together with a ``push forward'' map $\vec f \colon X \to f_\ast X$ in $\A$ lying over $f$:
\[
\cd{
    X \ar[rr]^{\vec f} & {} \ar@{.>}[d(0.25)];[d(0.75)]_{U} & f_\ast X \\
    C \ar[rr]_{f} & & D\text.
}\] The universality of this lifting amounts to saying that whenever
we are given a map $g \colon X \to Z$ in $\A$ whose underlying map
in $\C$ factors through $f \colon C \to D$, there is a unique
lifting to a factorisation of $g$ through $\vec f$ in $\A$. If every
cocartesian lifting exists we call $U \colon \A \to \C$ an
\defn{opfibration}; in this case we can think of $U$ as manifesting
$\A$ as a category ``indexed over $\C$''. There is no encyclopaedic
reference dealing with (op)fibrations, but one might consult
\cite{Borceaux}, for example.
\begin{Ex} For the domain functor $\dom \colon \Ar \C \to \C$, a cocartesian lifting of a map $h \colon A \to C$ of
$\C$ at an object $f \colon A \to B$ of $\Ar \C$ is given by a \emph{pushout} of $f$ along $h$:
\[
\cd{
    A \ar[rr]^{h} \ar[d]_f & & C \ar[d]^{g} \\
    B \ar[rr]_{k} & {} \ar@{.>}[d(0.25)];[d(0.75)]^{\cod} & D \pullbackcorner \\
    A \ar[rr]_{h} & & C\text.
}\] In particular, if $\C$ has all pushouts, then $\dom \colon \Ar \C \to \C$ is an opfibration.
\end{Ex}
Suppose that as well as $U \colon \A \to \C$, we have a further
functor $U' \colon \B \to \C$ together with a morphism $F \colon \A
\to \B$ commuting with the projections to $\C$. Then we say that $F$
\defn{creates cocartesian liftings over $\C$} if, for every
cocartesian map $g$ of $\B$ together with a lifting of its source to
$\A$, there exists a unique cocartesian map $\overline g$ of $\A$
satisfying $F\overline g = g$.
\begin{Ex}
Suppose we are given a class of maps $\E$ in a category $\C$ which are \emph{stable under pushout},
in the weak sense that whenever a pushout of an $\E$-map exists, it is another $\E$-map; then,
viewing $\E$ as a full subcategory of $\Ar \C$, we can express this by saying that the inclusion
functor $\E \to \Ar \C$ creates cocartesian liftings over $\C$.
\end{Ex}
The generalisation of this last example to the present situation is now immediate:
\begin{Prop}\label{closedunderpushout}
The category of $\Ll$-maps is \defn{stable under pushout} along arbitrary maps of $\C$, in the
sense that in the diagram
\[ \cd[@!@-2.5em]{ \Lmap \ar[dr]_{\dom \cdot U_{\Ll}} \ar[rr]^{U_{\Ll}}  & & \Ar \C \ar[dl]^{\dom} \\ & \C\text,}\]
the horizontal arrow creates cocartesian liftings over $\C$.
\end{Prop}
\noindent In other words, we must show that given a pushout square
\[
\cd[@C-1em]{
    A \ar[rr]^{h} \ar[d]_f & & C \ar[d]^{g} \\
    B \ar[rr]_{k} & {} & D \pullbackcorner \\
}\] in $\C$, together with an $\Ll$-map structure $(f, s)$ on $f$, there is a unique lift of $g$ to
an $\Ll$-map $(g, t)$ with respect to which $(h, k)$ is a cartesian arrow.
\begin{proof}
The argument we give here is due to \cite{IK}. We have the following situation:
\[
\cd[@!0@+1em]{
    & (\id_A, \star) \ar[dr]^{(\id_A, f)} \ar[dl]_{(h, h)} & &  & & \id_A \ar[dr]^{(\id_A, f)} \ar[dl]_{(h, h)}
    \\ (\id_C, \star) & & (f, s) \ar@{.>}[r(0.7)];[r(1.5)]^{U_\Ll} & & \id_C \ar[dr]_{(\id_C, g)} & & f \ar[dl]^{(h, k)}
    \\ & & & & & g\text, }
\]
where $\star$ represents the unique $\Ll$-coalgebra structures on the isomorphisms $\id_A$ and
$\id_C$, and we are applying Corollary \ref{tech} to deduce the existence of the arrows on the
left-hand side. Now, the right-hand diagram is a pushout in $\Ar \C$ and so, since $U_\Ll$ creates
colimits, there has a unique lifting of it to a pushout diagram in $\Lmap$. Thus we have a lifting
of $g$ to some $(g, t) \in \Lmap$ together with an arrow $(f, s) \to (g, t)$; and the universal
property of pushout says that this arrow is cocartesian.
\end{proof}
It is quite useful to have a concrete description of the induced $\Ll$-algebra structure on $g$ in
the situation of the previous proof, and a short calculation shows that it is given by the induced
map $t \colon D \to Eg$ in the following diagram:
\[
\cd{
 A
   \ar[d]_f
   \ar[r]^h &
 C
    \ar[d]^g    \ar@/^8pt/[ddr]^{\lambda_g}
\\
 B
    \ar[r]_k    \ar@/_8pt/[drr]_{E(h, k) \cdot s} &
 D \pullbackcorner
    \ar@{.>}[dr]^t
 \\
 & & Eg\text.
}
\]

The final thing we wish to consider in this section is the issue of \emph{closure under retracts}.
In the case of a plain \wfs, the class of $\ELL$-maps is closed under all retracts. The same is not
true here: indeed, if we are given an $\Ll$-map $(g, s)$ and a retract diagram
\[\cd{
    U \ar@{ >->}[r]^-{i_1} \ar[d]_f &
    V \ar@{->>}[r]^-{p_1} \ar[d]^g &
    U \ar[d]^f \\
    X \ar@{ >->}[r]_-{i_2} &
    Y \ar@{->>}[r]_-{p_2} &
    X
} \]
 in $\Ar \C$ (so $p_1i_1 = \id_U$ and $p_2i_2 = \id_X$) and want to make $f$ into an
$\Ll$-map, then the only logical way to do so is via the map
\begin{equation}\label{retractdiagram} r := X \xrightarrow{i_1} Y \xrightarrow{s} {E}g
\xrightarrow{{E}(p_1, p_2)} {E}f\text,\end{equation} and, although this satisfies the first two
axioms for an $\Ll$-map, it need not satisfy the third. However, category theory provides us with
conditions under which the above procedure will work: roughly speaking, if we can find a further
$\Ll$-map $(h, t)$ which ``measures'' $(g, s)$ in a suitable sense. Explicitly, we have:
\begin{Defn}
A \defn{contractible equaliser} in $\Ar \C$ is a diagram
\[\cd{f \ar@{ >->}[r]^i & g \ar@{->>}@/^2ex/[l]^p \ar@{ >->}@<-0.6ex>[r]_j \ar@<0.6ex>[r]^k & h \ar@{->>}@/^3.8ex/[l]^q}\]
such that $(i, p)$ and $(j, q)$ are retracts satisfying $ji = ki$ and $qk = ip$. A
\defn{contractible pair} in $\Ll$-$\cat{Map}$ is given by a contractible equaliser in $\Ar \C$
together with a lifting of $k$ and $j$ to morphisms of $\Ll$-$\cat{Map}$:
\[\cd{(g, s) \ar@<-0.6ex>[r]_j \ar@<0.6ex>[r]^k & (h, t)\text.}\]
\end{Defn}
In this richer setting, where the retract $(i, p)$ forms part of a contractible pair,
\eqref{retractdiagram} \emph{does} give a valid $\Ll$-map structure on $f$, and moreover one for
which $(i_1, i_2)$~--~though not necessarily $(p_1, p_2)$~--~becomes a morphism of
$\Ll$-$\cat{Map}$. This is a consequence of the standard result that \emph{the forgetful functor
$U_{\Ll} \colon \Ll\text-\cat{Map} \to \Ar \C$ creates equalisers for contractible pairs}; this is
part of the proof of Beck's monadicity theorem which can be found in any good book on category
theory: \cite{TTT}, for example. In such a situation, we shall call $(f, r)$ a \defn{retract
equaliser} of $(g, s)$. To summarise, we have that:
\begin{Prop}\label{closedunderretracts}
The category of $\Ll$-maps is \defn{closed under retract equalisers}, in that the forgetful functor
$U_{\Ll} \colon \Ll\text-\cat{Map} \to \Ar \C$ creates equalisers for contractible pairs.
\end{Prop}
We can now use this to give the analogue for \nwfs's of the so-called ``retract argument'' for
plain \wfs's.
\begin{Prop}\label{retractargument}
Every  $\Ll$-map $(f, r)$ is a retract equaliser of a cofree one: that is, one of the form
$(\lambda_f, \sigma_f)$.
\end{Prop}
\begin{proof}
This is another standard part of the monadicity theorem. The retract in question is
\[\cd{
    X \ar@{=}[r] \ar[d]_f &
    X \ar@{=}[r] \ar[d]^{\lambda_f} &
    X \ar[d]^f \\
    Y \ar@{ >->}[r]_-{r} &
    {E}f \ar@{->>}[r]_-{\rho_f} &
    Y
} \] which appears in the following contractible pair:
\[
\cd{F_{\Ll}f \ar@<-0.6ex>[r]_-{Lr} \ar@<0.6ex>[r]^-{\Sigma_f} & F_{\Ll}Lf} \quad \dashrightarrow
\quad \cd[@C+0.5em]{f \ar@{
>->}[r]^-{r} & Lf \ar@{->>}@/^3ex/[l]^-{\Phi_f} \ar@{ >->}@<-0.6ex>[r]_-{Lr} \ar@<0.6ex>[r]^-{\Sigma_f} & LLf\text. \ar@{->>}@/^4.5ex/[l]^-{\Phi_{Lf}}}\qedhere\]
\end{proof}

\subsection{Characterising $\Ll^1$-maps for cofibrantly generated \nwfs's}
We now apply the results we have just proved to give a characterisation of the $\Ll^1$-maps for the
one-step comonad generated by a set of maps in a cocomplete category $\C$. We know from Proposition
\ref{canonicalstructure} that every generating cofibration is an $\Ll^1$-map; we know from
Proposition \ref{closedundercolimits} that every coproduct of $\Ll^1$-maps is an $\Ll^1$-map; and
we know from Proposition \ref{closedunderpushout} that every pushout of an $\Ll^1$-map is an
$\Ll^1$-map. So in particular, \emph{every pushout of a coproduct of generating cofibrations is an
$\Ll^1$-map}.

There now arises a very natural question: observe that the functor $L^1 \colon \Ar \C \to \Ar \C$
sends a map $g$ to a pushout of a coproduct of generating cofibrations, so that if we equip each of
these generating cofibrations with its canonical $\Ll^1$-map structure we obtain a lifting of $L^1$
through $\Ll^1\text-\cat{Map}$:
\[
 \cd[@!@-2.5em]{ \Ar \C \ar[rr]^{\overline{L^1}} \ar[dr]_{L^1} & & \Ll^1\text-\cat{Map} \ar[dl]^{U_{\Ll^1}} \\ & \Ar \C}
\qquad \text{where} \qquad \overline{L^1} g = \spn{h_x}_\ast \big(\sum_{x \in S_g} (f_s,
\alpha_{f_s})\big)\text.
\]
Now, we have another lifting of $L^1$ through $\Ll\text-\cat{Map}$~--~namely, the free functor
$F_\Ll \colon \Ar \C \to \Ll\text-\cat{Map}$~--~and since, morally, both of these liftings do
exactly the same thing, we might wonder if they are naturally isomorphic. In fact, the result is
even stronger:
\begin{Prop}\label{characterisel1maps}
With the notation of the previous discussion, we have $\overline{L^1} = F_{\Ll^1}$.
\end{Prop}
\begin{proof}
Either a somewhat fiddly calculation, or, using the more abstract language of the Appendix, a
straightforward manipulation with universal properties; it can be found as Proposition II.4.2 of
\cite{Dbc}.
\end{proof}
This result implies, in particular, that \emph{every cofree
$\Ll^1$-map~--~}i.e., one of the form $(\lambda^1_g,
\sigma^1_g)$~--~\emph{is a pushout of a coproduct of generating
cofibrations}. This is almost a complete characterisation of the
$\Ll^1$-maps, and the final step, as in the case of plain \wfs's, is
to apply the retract argument. From Proposition
\ref{retractargument} we now deduce:
\begin{Prop}\label{characteriseonestep} The coalgebras for the one-step comonad $\Ll^1$ generated
by a set $J$ are precisely the retract equalisers of pushouts of coproducts of generating
cofibrations.
\end{Prop}

\subsection{Compositional properties}\label{Sec:compositionalprops}
There is one further property which we expect of the $\Ll$-maps and $\Rr$-maps for a \nwfs: namely,
that they should be closed under composition, and even \emph{transfinite} composition. Importantly,
these properties do \emph{not} hold in general for a mere comonad over $\dom$ or monad over $\cod$;
nonetheless, we can still frame our results in these broader settings. As before, we prefer the
comonadic version, but note that everything we do in this section applies equally well on the
monadic side.

To express the notion of being closed under composition, we will use a suitable monoidal structure
$(\I, \bullet)$ on the category $\cat{Cat} / \Ar \C$. The idea is that, if we view objects of
$\cat{Cat} / \Ar \C$ as being abstract categories of ``structured maps of $\C$'', then
$\bullet$-monoid structures on an object should correspond to ways of \emph{composing} these
structured maps. The monoidal structure in question has unit $\I$ given by:
\[
\I = (\C \xrightarrow{\id_{(\thg)}} \Ar \C)\text;\] whilst the tensor product of $U_\A \colon \A
\to \Ar \C$ and $U_\A \colon \B \to \Ar \C$ is given as follows. First we form the diagram
\[
\cd{
 \A \bullet \B
    \ar[rr]
    \ar[dd] \pushoutcorner \ar@{.>}[dr]|{U_{\A, \B}} & &
 \A
    \ar[d]^{U_\A} \\ & \C^{\b 3} \ar[r] \ar[d] \pushoutcorner &
 \Ar \C
    \ar[d]^{\dom} \\
 \B
    \ar[r]_{U_\B} &
 \Ar \C
    \ar[r]_{\cod} &
 \C
}
\]
in which both squares are pullbacks. $\C^{\b 3}$ is the functor category $[\cdot \to \cdot \to
\cdot, \C]$, whose objects are pairs of composable arrows in $\C$, and so we can consider the
evident
 ``composition'' functor $\text{cmp} \colon \C^{\b 3} \to \Ar \C$. We now define the
tensor product of $(\A, U_\A)$ and $(\B, U_\B)$ to be $(\A \bullet \B, U_{\A \bullet \B})$, where
\[U_{\A \bullet \B} = \A \bullet \B \xrightarrow{U_{\A, \B}} \C^{\b 3} \xrightarrow{\text{cmp}} \Ar \C\text.\]
So a typical element of $(\A \bullet \B)$ has the form $(\sigma, \tau)$, where $\sigma \in \B$ lies
over $f \colon X \to Y$ and $\tau \in \A$ lies over $g \colon Y \to Z$, whilst the projection onto
$\Ar \C$ is given by $U_{\A \bullet \B}(\sigma, \tau) = gf$.

The result we aiming for is:
\begin{Prop}\label{closedundercomp}
The category of $\Ll$-maps for a \nwfs\ is \defn{closed under composition}, in that $U_\Ll \colon
\Lmap \to \Ar \C$ is a $\bullet$-monoid in $\cat{Cat} / \Ar \C$.
\end{Prop}
\noindent Note that this statement says not only that $\Ll$-maps are closed under composition, but
also that this composition is associative and unital. As foreshadowed in Section
\ref{iteratingtheonestep}, our method for proving it will be to show that the functor $\G \colon
\cat{Comon} \to \cat{Cat} / \Ar \C$ which sends a comonad over $\dom$ to its category of coalgebras
is \emph{lax monoidal}. In particular, $\G$ sends monoids to monoids, so that if we have a
$\otimes$-monoid $(\Ll, \eta, \mu)$ in $\cat{Comon}$~--~i.e., a \nwfs\ on $\C$~--~then we induce a
$\bullet$-monoid structure on $\Lmap$ in $\cat{Cat} / \Ar \C$. Thus we will have proved the
previous Proposition if we can prove:
\begin{Prop}\label{GLax}
$\G$ is a lax monoidal functor $(\cat{Comon}, \otimes, I) \to (\cat{Cat} / \Ar \C, \mathord
\bullet, \I)$.
\end{Prop}
\noindent In order to do this, we need to provide an explicit description of the monoidal structure
$(\otimes, I)$ on $\cat{Comon}$. The unit $I$, of course, we understand: it is the comonad on $\Ar
\C$ which sends a map $f \colon X \to Y$ to $\id_X \colon X \to X$. The tensor product $\Ll^2
\otimes \Ll^1$ of two comonads is somewhat trickier to describe, primarily for notational reasons.
Its underlying functorial factorisation is:
\[X \xrightarrow{f} Y \quad \mapsto \quad X \xrightarrow{\lambda^2_{R^1f} \cdot \lambda^1_f} {E}^2R^1f \xrightarrow{\rho^2_{R^1f}} Y\text;\]
and to make it into a comonad over $\dom$, we must give maps
\[\sigma^{2 \otimes 1}_f \colon E^2R^1f \to E^2R^1L^{2 \otimes 1}f\text,\]
where we write $L^{2 \otimes 1}$ for the functor sending $f$ to $\lambda^2_{R^1f} \cdot
\lambda^1_f$. If we extract the relevant data from Remark \ref{otherbialgview} and Proposition
\ref{2fold}, we find that $\sigma^{2 \otimes 1}_f$ is given by:
\[ E^2R^1f \xrightarrow{\sigma^2_{R^1f}} E^2L^2R^1f \xrightarrow{E^2(\sigma^1_f, 1)} E^2(\lambda^2_{R^1f} \cdot \rho^1_{L^1f}) \xrightarrow{E^2\big(E^1(1, \lambda^2_{R^1f}), 1\big)} E^2R^1L^{2 \otimes 1}f\text.\]
With this in place, we can now prove Proposition \ref{GLax}. For $\G$ to be a lax functor, we first
need to provide a unit comparison map
\[ \cd[@!@-2.5em]{ \C \ar[dr]_{\id_{(\thg)}} \ar[rr]^{m_I}  & & I\text-\cat{Map} \ar[dl]^{U_I} \\ & \Ar \C\text.}\]
But $I$-$\cat{Map}$ is just the full subcategory of $\Ar \C$ whose objects are the isomorphisms,
and so we can take the obvious factorisation of $\id_{(\thg)}$ through $I$-$\cat{Map}$. We also
need to provide comparison maps
\[ \cd[@!0@+3.5em]{ \Ll^2\text-\cat{Map} \bullet \Ll^1\text-\cat{Map} \ar[dr]_{U_{\Ll^2 \bullet \Ll^1}} \ar[rr]^{m_{\Ll^2, \Ll^1}}  & & (\Ll^2 \otimes \Ll^1)\text-\cat{Map} \ar[dl]^{U_{\Ll^2 \otimes \Ll^1}} \\ & \Ar \C\text.}\]
On objects $m_{\Ll^1, \Ll^2}$ is given by specifying, for every $\Ll^1$-map $(f, s) \colon X \to Y$
and $\Ll^2$-map $(g, t) \colon Y \to Z$, an $(\Ll^2 \otimes \Ll^1)$-map structure on $gf \colon X
\to Z$. We take this to be the following map $u \colon Z \to E^2R^1(gf)$:
\[Z \xrightarrow{t} E^2g \xrightarrow{E^2(s, 1)} E^2(g \cdot \rho^1_f) \xrightarrow{E^2(E^1(1, g), 1)} E^2R^1(gf)\text.\]
Verifying the coalgebra axioms is routine. On morphisms, we have no real choice: a morphism on the
left hand side is a triple $(h, k, l)$ where $(h, k) \colon (f, s) \to (f', s')$ is a morphism
$\Ll^1$-$\cat{Map}$ and $(k, l) \colon (g, t) \to (g', t')$ is a morphism of $\Ll^2$-$\cat{Map}$,
and we are forced to send this to the morphism $(h, l) \colon (gf, u) \to (g'f', u')$ of $(\Ll^2
\otimes \Ll^1)$-$\cat{Map}$. Of course, we have to check that this morphism is a coalgebra
morphism; but this is again routine, as are the remaining details of the proof, which amount to
nothing more than checking a large number of coherence axioms. \hfill $\Box$

\vskip0.5\baselineskip Thus we conclude that the category of $\Ll$-maps for a \nwfs\ is closed
under composition; explicitly, if we are given two $\Ll$-maps $(f, s) \colon X \to Y$ and $(g, t)
\colon Y \to Z$, then their composite is given by $(gf, u) \colon X \to Z$, where $u \colon Z \to
E(gf)$ is given by:
\begin{equation}\label{compositionrule}
Z \xrightarrow{t} Eg \xrightarrow{E(s, 1)} E(g \cdot \rho_f) \xrightarrow{E(E(1, g), 1)} ER(gf)
\xrightarrow{\pi_{gf}} E(gf)\text.\end{equation} We have entirely dual results for the $\Rr$-maps,
but it might be worth spelling these out, since a little care is required. In this case, we
consider the category $\cat{Mon} := \cat{Mon}_\otimes(\Ff_C)$ of $\otimes$-monoids in $\Ff_\C$,
which we can view as the category of ``monads over $\cod$''. Once more we have a functor into
$\cat{Cat} / \Ar \C$, assigning to each such monad its category of algebras, but it is now a
\emph{contravariant} functor \[\H \colon \cat{Mon}^\op \to \cat{Cat} / \Ar \C\text.\] We can show
as before that $\H$ is lax monoidal, where the monoidal structure on the left-hand side is now the
$(\odot, \bot)$ monoidal structure. In particular, $\H$ sends $\odot$-monoids in $\cat{Mon}^\op$,
which are $\odot$-\emph{co}monoids in $\cat{Mon}$, to $\bullet$-monoids in $\cat{Cat} / \Ar \C$;
and since $\odot$-comonoids in $\cat{Mon}$ are precisely \nwfs's, we deduce that the category of
$\Rr$-maps for a \nwfs\ is closed under composition.

We turn now from finite to transfinite composition of maps. In general, if $\gamma$ is an ordinal
and $X$ is a $\gamma$-indexed chain in some category $\D$~--~in other words, a functor $X \colon
\gamma \to \D$~--~then the \defn{transfinite composite} of $X$ exists just when the colimit of $X$
does, and is given by the injection of $X_0$ into the colimit. Thus to say that a category admits
transfinite composition is simply to say that it admits colimits of chains.

In order to apply this notion to the $\Ll$-maps for a \nwfs, we need to form them into a category
which is different from the category $\Lmap$: namely, the category whose objects are those of $\C$
and whose set of morphisms from $X$ to $Y$ is the set of $\Ll$-maps $(f, s) \colon X \to Y$. By
Proposition \ref{closedundercomp}, this does gives us a category, whose composition law is given by
\eqref{compositionrule} and whose identity at $X$ is given by the unique lifting of $\id_X$ to an
$\Ll$-map. We shall denote this category by $\C_\Ll$ (and correspondingly $\C_\Rr$); observe that
whereas $\Lmap$ had a forgetful functor to $\Ar \C$, the category $\C_\Ll$ has a forgetful functor
$U \colon \C_\Ll \to \C$.
\begin{Prop}\label{closedtransfinite}
The category of $\Ll$-maps for a \nwfs\ is \defn{closed under transfinite composition}, in that the
forgetful functor $U \colon \C_\Ll \to \C$ creates colimits of chains.
\end{Prop}
\noindent The proof is surprisingly complex. First we need two lemmas:
\begin{Lem}
Colimits in $\Lmap$ commute with composition in the sense that the functors $i \colon \C \to \Lmap$
and $m \colon \Lmap \bullet \Lmap \to \Lmap$ exhibiting $\Lmap$ as a $\bullet$-monoid preserve
colimits strictly.
\end{Lem}
\begin{proof}
This is clear for $i \colon \C \to \Lmap$, which sends $X \in \C$ to the unique $\Ll$-map structure
$(\id_X, \star)$ on $\id_X$. To show the same for $m$, suppose that we are given a diagram $F
\colon \A \to \Lmap \mathbin{\bullet} \Lmap$ such that $\colim F$ exists. In particular, this
implies that the colimit of the underlying diagram $U_{\Ll \bullet \Ll} F \colon \A \to \Ar \C$
exists. But this is also the underlying diagram of $mF \colon \A \to \Lmap$, and so because $U_\Ll$
creates colimits, $mF$ also has a colimit, which moreover has the same underlying object in $\Ar
\C$ as $m \colim F$. All that remains to do is to check that the $\Ll$-coalgebra structures on
$\colim mF$ and $m \colim F$ agree. So suppose that each $F(a)$ is given by an $\Ll$-map $(f_a,
s_a) \colon X_a \to Y_a$ and an $\Ll$-map $(g_a, t_a) \colon Y_a \to Z_a$, and consider the
following diagram:
\[\cd[@+1em@C+1em]{
    \overrightarrow{Z_a} \ar[r]^-{\overrightarrow{t_a}} &
    \overrightarrow{Eg_a} \ar[r]^-{\overrightarrow{E(s_a, 1)}} \ar[d]_{\text{can}} &
    \overrightarrow{E(g_a \cdot \rho_{f_a})} \ar[r]^-{\overrightarrow{E(E(1, g_a), 1)}} \ar[d]_{\text{can}} &
    \overrightarrow{ER(g_a f_a)} \ar[r]^-{\overrightarrow{\pi_{g_af_a}}} \ar[d]^{\text{can}} &
    \overrightarrow{E(g_a f_a)} \ar[dd]^{\text{can}} \\
    &
    E \overrightarrow{g_a} \ar[r]_-{E(\overrightarrow{s_a}, 1)} &
    E(\overrightarrow{g_a \cdot \rho_{f_a}}) \ar[d]_{E(\text{can}, 1)} \ar[r]_{E(\overrightarrow{E(1, g_a)}, 1)} &
    E \overrightarrow{R(g_af_a)} \ar[d]^{\text{can}} \\
    &
    &
    E(\overrightarrow{g_a} \cdot \rho_{\overrightarrow{f_a}}) \ar[r]_{E(E(1, \overrightarrow{g_a}), 1)} &
    ER(\overrightarrow{g_a f_a}) \ar[r]_{\pi_{\overrightarrow{g_af_a}}} &
    E(\overrightarrow{g_a f_a})\text.
}\] Here, $\overrightarrow{Z_a}$ is shorthand for $\colim_a Z_a$, and so on. Each small square
commutes, and thus the two extremal routes around the big diagram are the same: but the upper of
these is the $\Ll$-coalgebra structure on $\colim mF$, whilst the lower is the coalgebra structure
on $m \colim F$.
\end{proof}
\begin{Lem}
Suppose that $(f, s) \colon X \to Y$ and $(g, t) \colon Y \to Z$ are $\Ll$-maps, and that $(gf, u)$
is their composite according to \eqref{compositionrule}. Then the morphism $(1_X, g) \colon f \to
gf$ of $\Ar \C$ is a morphism of $\Ll$-maps $(1_X, g) \colon (f, s) \to (gf, u)$.
\end{Lem}
\begin{proof}
We must show that $ug = E(1, g)s$, and so calculate
\begin{align*}
ug & = \pi_{gf} \cdot E(E(1, g), 1) \cdot E(s, 1) \cdot t \cdot g\\
& = \pi_{gf} \cdot E(E(1, g), 1) \cdot E(s, 1) \cdot \lambda_g\\
& = \pi_{gf} \cdot E(E(1, g), 1) \cdot \lambda_{g \cdot \rho_f} \cdot s\\
& = \pi_{gf} \cdot \lambda_{R(gf)} \cdot E(1, g) \cdot s\\
& = E(1, g) \cdot s
\end{align*}
as required.
\end{proof}
\begin{proof}[Proof of Proposition \ref{closedtransfinite}]
Suppose that we are given a $\gamma$-chain $X \colon \gamma \to \C_\Ll$: so for each $\alpha <
\gamma$ we have an object $X_\alpha$ of $\C$ and for each $\alpha < \beta < \gamma$ we have an
$\Ll$-map $(f_{\alpha, \beta}, s_{\alpha, \beta}) \colon X_\alpha \to X_\beta$. Suppose also that
we have a colimit for the underlying chain $UX \colon \gamma \to \C$; so we have a colimiting
object $Y$ of $\C$ and maps $g_\alpha \colon X_\alpha \to Y$ commuting with the $f_{\alpha,
\beta}$'s. We must show that this can be lifted to a colimit for $X$, for which we must equip each
$g_\alpha$ with an $\Ll$-map structure $(g_\alpha, t_\alpha)$ compatible with the $s_{\alpha,
\beta}$'s.

So given $\delta < \gamma$, we equip $g_\delta$ with an $\Ll$-map structure by considering the
chain $Z^\delta \colon \gamma \to \Lmap$ given as follows. For $\alpha \leqslant \delta$, we take
$Z^\delta$ to be constant with value $(\id_{X_\delta}, \star)$, and for $\alpha > \delta$ we take
$Z^\delta_\alpha = (f_{\delta, \alpha}, s_{\delta, \alpha})$ with connecting maps given by
$Z^\delta_{\alpha, \beta} = (\id_{X_0}, f_{\alpha, \beta})$: by the second of our two lemmas this
connecting map is a valid morphism of $\Ll$-maps.

Now, from the given colimit for $UX$ we obtain a colimit for the underlying chain $U_\Ll \cdot
Z^\delta \colon \gamma \to \Ar \C$, namely the object $g_\delta \colon X_\delta \to Y$ of $\Ar \C$;
and because the forgetful functor $U_\Ll$ creates colimits, we obtain from this a colimit for the
chain $Z^\delta \colon \gamma \to \Lmap$, which gives us the required $\Ll$-map structure
$(g_\delta, t_\delta)$ on $g_\delta$.

It remains to show that these $\Ll$-map structures $(g_\alpha, t_\alpha)$ are compatible with the
$s_{\alpha, \beta}$'s; explicitly, we need to show that $(g_\beta, t_\beta) \circ (f_{\alpha,
\beta}, s_{\alpha, \beta}) = (g_\alpha, t_\alpha)$. But by our first lemma, we know that
precomposing $\colim Z^\beta$ with $(f_{\alpha, \beta}, s_{\alpha, \beta})$ will give us the same
result as precomposing every element of $Z^\beta$ with $(f_{\alpha, \beta}, s_{\alpha, \beta})$ and
taking the colimit of the resultant chain $Z'$. But it is easy to see that $\colim Z'$ is precisely
$\colim Z^\alpha$ and so the result follows.
\end{proof}

\subsection{Characterising $\Ll$-maps for cofibrantly generated \nwfs's}
We would now like to use the results of the previous section to give a characterisation of the
$\Ll$-maps for a cofibrantly generated n.w.f.s. In one direction, this is straightforward: by
Section \ref{Sec:cofibrantgeneration} we know that every generating cofibration is an $\Ll$-map,
and by Propositions \ref{closedundercolimits}, \ref{closedunderpushout}, \ref{closedunderretracts}
and \ref{closedtransfinite}, we know that the $\Ll$-maps are closed under colimits, pushouts,
retract equalisers and transfinite composition. So in particular, \emph{every retract equaliser of
a transfinite composition of pushouts of coproducts of generating cofibrations is an $\Ll$-map}.

What is harder to come by is a result in the other direction, saying that \emph{every} $\Ll$-map is
of this form. We would like to mimic the argument we gave for the $\Ll^1$-maps, where we first
characterised the ``cofree'' $\Ll^1$-maps~--~that is, those of the form $(\lambda^1_f,
\sigma^1_f)$~--~and then applied Proposition \ref{retractargument}, the ``generalised retract
argument'', to produce a characterisation of an arbitrary $\Ll^1$-map.

The problem is that we have been unable to find a simple characterisation of the ``cofree''
$\Ll$-maps. The reason for this is the way in which these cofree maps are constructed: we perform a
transfinite induction, at each stage of which we first glue some extra cells on, and then
coequalise away the ones that we shouldn't have added because they were there already. If we could
show that this was equivalent to simply glueing slightly fewer cells on in the first place, then we
could conclude that every cofree map was just given by transfinitely glueing on cells, and then our
characterisation result would follow easily.

However, we can see no way of proving this statement, and so the
best we can do for the moment is to refer back to the plain \wfs\
underlying our cofibrantly generated n.w.f.s. We know that this is a
cofibrantly generated \wfs, and there \emph{is} a characterisation
of its left class of maps: they are precisely the retracts of
transfinite composites of pushouts of coproducts of generating
cofibrations.

\begin{Prop}
Every retract equaliser of a transfinite composition of pushouts of coproducts of generating
cofibrations is an $\Ll$-map, whilst the underlying morphism in $\C$ of any $\Ll$-map is a retract
of a transfinite composition of pushouts of coproducts of generating cofibrations.
\end{Prop}

\section*{Appendix: universality of the one-step comonad}
The purpose of this appendix is to expand upon the abstract
description of the one-step comonad $\Ll^1$ generated by a set of
maps $J$ which we hinted at in Section \ref{Sec:onestep}, and to use
it to give a proof of Proposition \ref{goodunivprop}, which, we
recall, told us that that $\Ll^1$ is \emph{freely} generated by $J$.

We will do this by using a certain amount of the theory of
\emph{2-categories}. Now, a 2-category is a category which has not
only objects $X, Y, Z, \dots$ and morphisms $f, g, h, \dots$ but
also 2-cells $\alpha \colon f \Rightarrow g$ between these 1-cells
which can be composed together in various ways, subject to axioms
which make any multiple composites we might form unambiguous; for a
good introduction to the subject the reader might refer to
\cite{Elements}. The ur-2-category is $\cat{Cat}$ whose objects,
1-cells and 2-cells are respectively, (large) categories, functors
and natural transformations; and with this in mind, one can crudely
think of two-dimensional category theory as being \emph{abstract
category theory}, in the same way that topos theory is
\emph{abstract set theory} and model category theory is
\emph{abstract homotopy theory}. We will be using it to ``do
category theory'' in one particular 2-category, which is a close
relative of $\cat{Cat}$:
\begin{Defn}
For $\C$ a category, the 2-category $\cat{Cat} / \C$ of ``categories over $\C$'' has:
\begin{itemize*}
\item \textbf{Objects} $(\A, U_\A)$ being categories $\A$ together with a functor $U_\A \colon \A \to
\C$;
\item \textbf{1-cells} $F \colon (\A, U_\A) \to (B, U_\B)$ being functors $H \colon \A \to \B$
such that
\[
\cd[@C-1em@R+1em]{
    \A \ar[rr]^F \ar[dr]_{U_\A} & & \B \ar[dl]^{U_\B}
    \\ & \C
}
\]
commutes;
\item \textbf{2-cells} $\alpha \colon F \Rightarrow G $ being natural transformations $\alpha
\colon F \Rightarrow G$ such that $U_\B \cdot \alpha = \id_{U_\A}$; or diagramatically
\[
\cd[@C-1em@R+1em]{
    \A \ar@/^0.9em/[rr]^F \ar@/_0.9em/[rr]_G \dtwocell{rr}{\alpha} \ar[dr]_{U_\A} & & \B \ar[dl]^{U_\B}
    \\ & \C\text.
}
\]
\end{itemize*}
\end{Defn}
The way one ``does category theory'' in a 2-category is by recognising that concepts which we are
familiar with in $\cat{Cat}$ are definable purely in terms of diagrams of objects, 1-cells and
2-cells and so can be defined in an arbitrary 2-category. For example, we can define a
\defn{comonad} in a 2-category
$\K$ to be an object $X$, together with a 1-cell $t \colon X \to X$
and a pair of 2-cells $\epsilon \colon t \Rightarrow \id_X$ and
$\Delta \colon t \Rightarrow tt$ making the following two diagrams
commute:
\[
\cd{
 & t \ar@2[dl]_{\Delta} \ar@2[dr]^{\Delta} \ar@2[d]^{\id_t} \\
 t & tt \ar@2[l]^{t\epsilon} \ar@2[r]_{\epsilon t} & t\text,
} \qquad \cd{
 t \ar@2[r]^{\Delta} \ar@2[d]_{\Delta} &
 tt \ar@2[d]^{t \Delta} \\
 tt \ar@2[r]_{\Delta t} &
 ttt\text.
}
\]
In $\cat{Cat}$, this reduces to the usual notion of comonad; whilst a comonad in the 2-category
$\cat{Cat} / \C$ consists of:
\begin{itemize*}
\item An object $(U_\A \colon \A \to \C)$;
\item A functor $T \colon \A \to \A$ satisfying $U_\A \cdot T = T$, and
\item Natural transformations $\epsilon \colon T \Rightarrow \id_\A$ and $\Delta \colon T
\Rightarrow TT$ satisfying $U_\A \cdot T = \id_{U_\A}$,
\end{itemize*}
and such that $(T, \epsilon, \Delta)$ is a comonad in the usual sense. The relevance of this notion
becomes manifest when we observe that a comonad on the object $(\dom \colon \Ar \C \to \C)$ is
precisely what we have been calling a ``comonad over $\dom$'': that is, a functorial factorisation
$(E, \lambda, \rho)$ together with a comonad structure on the associated pair $(L, \Phi)$.

As a further instance of this process of abstraction, consider the
notion of a \emph{left Kan extension}. We can define an analogous
notion, which is usually known just as \emph{left extension}, in an
arbitrary 2-category $\K$: given 1-cells $f \colon X \to Y$ and $g
\colon X \to Z$, we say that a pair $(h \colon Y \to Z, \theta
\colon g \Rightarrow hf)$, as in
\[\cd{
    & X \ar[dr]^g \ar[dl]_f \ltwocell{d}{\theta} \\
  Y \ar[rr]_h & & Z
}\] exhibits $h$ as the\footnote{Note that a left extension, if it exists, is unique up to unique
isomorphism, and so we can speak with justifiable looseness of \emph{the} left extension.}
\textbf{left extension} of $g$ along $f$ if, given any diagram
\[\cd{
    & X \ar[dr]^g \ar[dl]_f \ltwocell{d}{\psi} \\
  Y \ar[rr]_k & & Z
}\] there is a unique 2-cell $\phi \colon h \Rightarrow k$ such that
\[
\cd{
    & X \ar[dr]^g \ar[dl]_f \ltwocell{d}{\theta} \dtwocell[1.3]{d}{\phi} \\
  Y \ar[rr]^h \ar@/_2.2em/[rr]_{k} & & Z
} \qquad = \qquad \cd{
    & X \ar[dr]^g \ar[dl]_f \ltwocell{d}{\psi} \\
  Y \ar[rr]_k & & Z\text.
}\] In the 2-category $\cat{Cat}$, the notion of left extension is the familiar notion of left Kan
extension\footnote{Though some authors reserve the name left Kan extension for a slightly stronger
notion: see the discussion in Section 4.3 of \cite{Kellybook}.}; and what we will be interested in
is left extensions in the 2-category $\cat{Cat} / \C$. The reason for this is that we can use left
extensions to construct comonads:
\begin{Prop}
Let $\K$ be a 2-category, and suppose that $f \colon X \to Y$ is a 1-cell of $\K$ such that the
left extension of $f$ along itself exists, and is given by:
\[\cd{
    & X \ar[dr]^f \ar[dl]_f \ltwocell{d}{\theta} \\
  Y \ar[rr]_t & & Y\text.
}\] Then $t$ can be made into a comonad for which $\theta \colon f
\Rightarrow tf$ is a ``coaction'' of $t$ on $f$.
\end{Prop}
\begin{proof} This is essentially an exercise in the
universal property of the left extension. Since we have a diagram
\[\cd{
    & X \ar[dr]^f \ar[dl]_f \ltwocell{d}{\id_{f}} \\
  Y \ar[rr]_{\id_Y} & & Y\text,
}\] we induce, by the universal property of left extension, a 2-cell
$\epsilon \colon t \Rightarrow \id_{Y}$ such that $\epsilon \cdot
\theta = \id_f$. Moreover, we have a diagram
\[\cd[@C-1.5em]{
   & & X \ar[drr]^f \ar[dll]_f \ar[d]^{f} & &  \\
  Y \ar[rr]_t & & Y \ar[rr]_t \ltwocell[0.35]{ull}{\theta}  \ltwocell[0.35]{urr}{\theta}  & & Y\text,
}\] and so, again by the universal property, we induce a 2-cell
$\Delta \colon t \Rightarrow hh$ such that $\Delta \cdot \theta =
t\theta \cdot \theta$. Applying the universal property again (the
``uniqueness'' part this time), we deduce that $(t, \epsilon,
\Delta)$ satisfies the axioms for a comonad in our 2-category $\K$.
Finally, to say that $\theta \colon f \Rightarrow tf$ is a
``coaction'' of $t$ on $f$ is to say that the following diagrams
commute:
\[
\cd{
 f \ar@2[d]_{\theta} \ar@2[dr]^{\id_f} \\
 tf \ar@2[r]_{\epsilon f} & f\text,
} \qquad \cd{
 f \ar@2[r]^{\theta} \ar@2[d]_{\theta} &
 tf \ar@2[d]^{t \theta} \\
 tf \ar@2[r]_{\Delta f} &
 ttf\text,
}
\]
which follows from yet another application of the universal property.
\end{proof}
This comonad is known as the \emph{density comonad} of $f \colon X
\to Y$: the most comprehensive source of information on density
comonads~--~or rather the dual \emph{codensity monads}~--~is also
the place where they were first introduced, namely the thesis of
Dubuc \cite{Dbc}. One important universal property demonstrated by
Dubuc is that the coaction $\theta \colon f \Rightarrow tf$ induced
by a density comonad  is the \emph{universal} coaction on $f$. To
make this precise, we first define a \defn{morphism of comonads}
$\alpha \colon (t, \epsilon, \Delta) \to (t', \epsilon', \Delta')$
on $Y$ to be given by a 2-cell $\alpha \colon t \Rightarrow t'$ in
$\K$ which is compatible with the comonad structures in that the two
diagrams
\[
\cd{
 t \ar@2[d]_{\alpha} \ar@2[dr]^{\epsilon} \\
 t' \ar@2[r]_{\epsilon'} & \id_t\text,
} \qquad \cd{
 t \ar@2[r]^{\Delta} \ar@2[d]_{\alpha} &
 tt \ar@2[d]^{\alpha \alpha} \\
 t' \ar@2[r]_{\Delta'} &
 t't'
}
\]
commute. Now the universality of the coaction $\theta \colon f
\Rightarrow tf$ given above can be expressed by saying that, for any
other comonad $(t', \epsilon', \Delta')$ on $Y$ and coaction
$\theta' \colon f \Rightarrow t'f$, there is a unique comonad
morphism $\alpha \colon (t, \epsilon, \Delta) \to (t', \epsilon',
\Delta')$ for which $\theta' = \alpha \theta$.

Let us now see what this universal property says in the 2-category
$\cat{Cat} / \C$ we are interested in. Observe first that to give a
coaction
\[\cd{
    & (\A, U_\A) \ar[dr]^F \ar[dl]_F \ltwocell{d}{\theta} \\
  (\B, U_\B) \ar[rr]_T & & (\B, U_\B)\text.
}\] of a comonad $\mathsf T = (T, \epsilon, \Delta)$ on a 1-cell $F$
is to give, for each $A \in \A$, a morphism $\theta_A \colon FA \to
TFA$ in $\B$ making $FA$ into a $T$-coalgebra in a manner that is
natural in morphisms of $\A$. We deduce\footnote{The reader with
some knowledge of two-dimensional category theory will observe the
formal theory of monads \cite{FTM} raising its head here.} that
coactions of $\mathsf T$ on $F$ are in bijection with liftings of $F
\colon \A \to \B$ through the category of coalgebras for $\mathsf
T$:
\[
\cd[@!C@C-2em@R+1em]{
    \A \ar[rr]^{\overline F} \ar[dr]_{F} & & \mathsf T\text-\cat{Coalg} \ar[dl]^{U_{\mathsf T}}
    \\ & \B\text,
}
\]
in which terms the above universal property can be restated as:
\begin{Prop}
Let $F \colon (\A, U_\A) \to (\B, U_\B)$ be a 1-cell of $\cat{Cat} /
\C$ which admits a codensity monad $\mathsf T = (T, \epsilon,
\Delta)$. Then there is a bijection, natural in $\mathsf T'$,
between morphisms of comonads $\mathsf T \to \mathsf T'$ on $(\B,
U_\B)$ and liftings
\[
\cd[@!C@C-2em@R+1em]{
    \A \ar[rr]^{\overline F} \ar[dr]_{F} & & \mathsf T'\text-\cat{Coalg} \ar[dl]^{U_{\mathsf T'}}
    \\ & \B\text,
}
\]
of $F$ through $\mathsf T'\text-\cat{Coalg}$.
\end{Prop}
In particular, let us consider the case where $F$ is the following
1-cell of $\cat{Cat} / \C$:
\[
\cd[@C-1em@R+1em]{
    J \ar[rr]^{\iota} \ar[dr]_{\dom \cdot \iota} & & \Ar \C \ar[dl]^{\dom}
    \\ & \C\text.
}
\]
Suppose for a moment that this 1-cell admits a left extension along
itself; now if we suggestively write the corresponding comonad on
$(\Ar \C, \dom)$ as $\Ll^1$, then the above Proposition reduces to
precisely Proposition \ref{goodunivprop}. Therefore we will have
proved this latter Proposition, and given our promised abstract
description of the comonad $\Ll^1$, if we can show that the left
extension of $\iota$ along itself exists. To do this, we use the
fact that the 2-category $\cat{Cat} / \C$ that we are working in is
very closely related to $\cat{Cat}$, so that under suitable
circumstances, we can build left extensions in $\cat{Cat} / \C$ from
left extensions in $\cat{Cat}$.
\begin{Prop}\label{liftslice}
Let $F$ and $G$ be 1-cells
\[\cd[@C-2em]{
    & (\A, U_\A) \ar[dr]^G \ar[dl]_F  \\
  (\B, U_\B) & & (\D, U_\D)
}\] in $\cat{Cat} / \C$, where $\A$ is a small category, $\D$ has colimits preserved by $U_\D$, and
$U_\D \colon \D \to \C$ is an opfibration. Then the left extension of $G$ along $F$ exists.
\end{Prop}
This is not the most general result possible: it suffices that $\D$ be
\emph{fibre-cocomplete}~--~so that every fibre category of $\D$ has colimits which every
cocartesian map preserves. However, the proof becomes much more intricate if we do so, and so we
shall content ourselves with this slightly weaker result.
\begin{proof}
As outlined above, we will construct the required left extension using left Kan extensions in
$\cat{Cat}$ together with the opfibration structure on $U_\D$. We begin by taking the left Kan
extension $(K, \psi)$ of $G$ along $F$, which we can do because $\A$ is small and $\D$ is
cocomplete:
\[\cd{
    & \A \ar[dr]^G \ar[dl]_F \ltwocell{d}{\psi} \\
  \B \ar[rr]_K & & \D\text.
}\] The problem is that this Kan extension does not respect the functors down to $\C$: there is no
reason for us to have $U_\D K = U_\B$. However, if we can produce a natural family of maps
\[\phi_b \colon U_\D Kb \to U_\B b\]
then we can use the opfibration structure on $U_\D$ to ``correct'' the element $Kb$ to an element
$Hb = (\phi_b)_\ast(Kb)$ which lies in the right fibre: in other words, we obtain a new functor $H
\colon \B \to \D$ that \emph{does} satisfy $U_\D H = U_\B$ as required. So all we need now is a
suitable natural transformation $\phi \colon U_\D K \Rightarrow U_\B$. To get this, observe that
because $U_\D$ preserves colimits, the diagram
\[\cd{
    & \A \ar[dr]^{U_\D G} \ar[dl]_F \ltwocell{d}{U_\D\psi} \\
  \B \ar[rr]_{U_\D K} & & \C
}\] is also a left Kan extension. Now, because $U_\B F = U_\D G = U_\A$, we have the diagram
\[\cd{
    & \A \ar[dr]^{U_\D G} \ar[dl]_F \ltwocell{d}{\id_{U_\A}} \\
  \B \ar[rr]_{U_\B} & & \C\text,
}\] and so by the universal property of Kan extension, we induce a natural transformation $\phi
\colon U_\D K \Rightarrow U_\B$ satisfying $F\phi \cdot U_\D \psi = \id_{U_\A}$. Now we can use
$\phi$ to ``correct'' $K$ as outlined above. Formally, we say that because $U_\D$ is an
opfibration, so is the functor $[\B, U_\D] \colon [\B, \D] \to [\B, \C]$, and so the natural
transformation $\phi$~--~seen as a morphism of $[\B, \C]$~--~induces a functor
\[\phi_\ast \colon [\B, \D]_{U_\D K} \to [\B, \D]_{U_\B}\text.\]
Applying this to $K \in [\B, \D]_{U_\D K}$ yields our ``corrected'' functor $H = \phi_\ast K \colon
\B \to \D$ satisfying $U_\D H = U_\B$, together with a natural transformation $\chi = \vec \phi
\colon K \Rightarrow H$ satisfying $U_\D \chi = \phi$. Since $U_\D H = U_\B$, we have a 1-cell $H
\colon (\B, U_\B) \to (\D, U_\D)$ in $\cat{Cat} / \C$; whilst if we define the natural
transformation $\theta \colon G \Rightarrow HF$ to be the composite
\[
\cd{
    & \A \ar[dr]^G \ar[dl]_F \ltwocell{d}{\psi} \dtwocell[1.3]{d}{\chi} \\
  \B \ar[rr]|K \ar@/_2.2em/[rr]_{H}  & & \D
}
\]
then we have $U_\D \theta = U_\D \chi F \cdot U_\D \psi = \phi F \cdot U_\D \psi = \id_{U_\A}$ so
that $\theta$ gives a 2-cell of $\cat{Cat} / \C$, which exhibits $H$ as the left extension of $G$
along $F$ in $\cat{Cat} / \C$:
\[\cd[@C-2em]{
    & (\A, U_\A) \ar[dr]^G \ar[dl]_F \ltwocell{d}{\theta} \\
  (\B, U_\B) \ar[rr]_{H} & & (\D, U_\D)\text.
}\] Verification of the universal property is left as an exercise to the reader.
\end{proof}
\begin{Cor}\label{Kanthing}
Let $\C$ be a cocomplete category, let $J$ be a set of maps in $\C$ and let $\iota \colon J \to \Ar
\C$ be the inclusion of the discrete subcategory $J$ into $\Ar \C$. Then the left extension
\[\cd[@C-2em]{
    & (J, \dom \cdot \iota) \ar[dr]^\iota \ar[dl]_\iota \ltwocell{d}{\theta} \\
  (\Ar \C, \dom) \ar[rr]_{L^1} & & (\Ar \C, \dom)
}\] exists in $\cat{Cat} / \C$, and thus generates a comonad over
$\dom$ $\Ll^1$ satisfying Proposition \ref{goodunivprop}.
\end{Cor}
This completes our abstract description of the comonad $\Ll^1$. Our
final task is to calculate explicitly what our machinery gives us,
and check that it tallies with the comonad we gave in Section
\ref{Sec:onestep}. Since the description of the comultiplication
that we gave there was secretly derived from the results of this
Appendix, the only thing we have to check is that our machinery
gives the right underlying functorial factorisation. Now, the first
thing the construction of Proposition \ref{liftslice} tells us to do
is to form the following left Kan extension:
\[\cd{
    & J \ar[dr]^\iota \ar[dl]_\iota \ltwocell{d}{\psi} \\
  \Ar \C \ar[rr]_K & & \Ar \C\text,
}\] the value of which at an object $g \colon C \to D \in \Ar \C$ is given by the colimit
\[Kg = \int^{(f, \gamma) \in \iota \downarrow g} \iota(f)\text.\]
Because $J$ is discrete, the indexing category $\iota \downarrow g$ degenerates to a set, an
element of which is a pair $(f, \gamma)$ where $f \in J$ and $\gamma \colon f \rightarrow g$ in
$\Ar \C$. So
 $\iota \downarrow g$ is the set $S_g$ that we considered before, and
\[Kg = \sum_{x \in S_g} f_x\text.\]
Having formed $K$, the next step is to ``correct'' it to a functor over $\dom$, which we do by
pushing out along the components of a natural transformation $\phi \colon \dom \cdot K \Rightarrow
\dom \colon \Ar \C \to \C$ obtained from the universal property of Kan extension. We find that in
this case, $\phi_g$ is the map:
\[\phi_g = \spn{h_x}_{x \in S_g} \colon \sum_x A_x \to C\]
(where $g \colon C \to D$ as before). Therefore $L^1g$ is given by the right-hand map in the
pushout diagram
\[
    \cd{
        \sum A_x
            \ar[r]^-{\left<h_x\right>}
            \ar[d]_{\sum_x f_x} &
        C
            \ar[d]^{\lambda^1_g}\\
        \sum B_x
            \ar[r] &
        E^1g\text, \pullbackcorner
    }
\]
which is precisely what we had in Section \ref{Sec:onestep}.

\bibliography{biblio}

\end{document}